\documentclass{article}

\usepackage{xspace,amsmath,amssymb,amsthm,amsfonts,epsfig,syntonly,dsfont,pifont,natbib}
\usepackage{cite,bm,color,url,textcomp}
\usepackage{algorithmicx,algorithm,algpseudocode}
\usepackage{epstopdf}
\usepackage{empheq,float}
\usepackage{graphicx,subfigure,balance}
\usepackage{multirow}
\usepackage{float}
\usepackage[dvipsnames]{xcolor}

\newtheorem{assumption}{Assumption}

\newtheorem{lemma}{Lemma}

\newtheorem{theorem}{Theorem}
\newtheorem{definition}{Definition}
\theoremstyle{definition}

\DeclareMathOperator*{\argmax}{arg\,max}
\DeclareMathOperator*{\argmin}{arg\,min}

\graphicspath{{./figs_modified/}}
\usepackage{wrapfig}

\usepackage{times}
\usepackage[margin=1.08in]{geometry}

\usepackage[utf8]{inputenc} 
\usepackage[T1]{fontenc}    
\usepackage{hyperref}       
\usepackage{url}            
\usepackage{booktabs}       
\usepackage{nicefrac}       
\usepackage{microtype}      

\usepackage{scrextend}

\title{\LARGE{Enhancing Parameter-Free Frank Wolfe \\with an Extra Subproblem}}

\author{ 
	Bingcong Li\textsuperscript{\rm 1} ~~
	Lingda Wang\textsuperscript{\rm 2} ~~
	Georgios B. Giannakis\textsuperscript{\rm 1} ~~
	Zhizhen Zhao\textsuperscript{\rm 2}  	\vspace{0.1cm} \\	 
	 \textsuperscript{\rm 1}\textit{University of Minnesota} \\
	 \textsuperscript{\rm 2}\textit{University of Illinois at Urbana-Champaign}\\
	 \texttt{\{lixx5599, georgios\}@umn.edu} \\
	 \texttt{\{lingdaw2, zhizhenz\}@illinois.edu} 
	 }

\begin{document}

\maketitle
\begin{abstract}
Aiming at convex optimization under structural constraints, this work introduces and analyzes a variant of the Frank Wolfe (FW) algorithm termed ExtraFW. The distinct feature of ExtraFW is the pair of gradients leveraged per iteration, thanks to which the decision variable is updated in a prediction-correction (PC) format. Relying on no problem dependent parameters in the step sizes, the convergence rate of ExtraFW for general convex problems is shown to be ${\cal O}(\frac{1}{k})$, which is optimal in the sense of matching the lower bound on the number of solved FW subproblems. However, the merit of ExtraFW is its faster rate ${\cal O}\big(\frac{1}{k^2} \big)$ on a class of machine learning problems. Compared with other parameter-free FW variants that have faster rates on the same problems, ExtraFW has improved rates and fine-grained analysis thanks to its PC update. Numerical tests on binary classification with different sparsity-promoting constraints demonstrate that the empirical performance of ExtraFW is significantly better than FW, and even faster than Nesterov's accelerated gradient on certain datasets. For matrix completion, ExtraFW enjoys smaller optimality gap, and lower rank than FW.
\end{abstract}

\section{Introduction}

The present work deals with efficient algorithms for solving the optimization problem
\begin{align}\label{eq.prob}
	\min_{\mathbf{x} \in {\cal X}} f(\mathbf{x})	
\end{align}
where $f$ is a smooth convex function, while the constraint set ${\cal X} \subset \mathbb{R}^d$ is assumed to be convex and compact, and  $d$ is the dimension of the variable $\mathbf{x}$. Throughout we denote by $\mathbf{x}^* \in {\cal X}$ a minimizer of \eqref{eq.prob}. For many machine learning and signal processing problems, the constraint set ${\cal X}$ can be structural but it is difficult or expensive to project onto. Examples include matrix completion in recommender systems \citep{freund2017} and image reconstruction \citep{harchaoui2015}, whose constraint sets are nuclear norm ball and total-variation norm ball, respectively. The applicability of projected gradient descent (GD) \citep{nesterov2004} and Nesterov's accelerated gradient (NAG) \citep{allen2014,nesterov2015} is thus limited by the computational barriers of projection, especially as $d$ grows large.

An alternative to GD for solving \eqref{eq.prob} is the Frank Wolfe (FW) algorithm \citep{frank1956,jaggi2013,lacoste2015}, also known as the `conditional gradient' method. FW circumvents the projection in GD by solving a subproblem with a \textit{linear} loss per iteration.
For a structural ${\cal X}$, such as the constraint sets mentioned earlier, it is possible to solve the subproblem
either in closed form or through low-complexity numerical methods \citep{jaggi2013,garber2015}, which saves computational cost relative to projection. In addition to matrix completion and image reconstruction, FW has been appreciated in several applications including structural SVM \citep{lacoste2013block}, video colocation \citep{joulin2014}, optimal transport \citep{luise2019sinkhorn}, and submodular optimization \citep{mokhtari2018}, to name a few.

Although FW has well documented merits, it exhibits slower convergence when compared to NAG. Specifically, FW satisfies $f(\mathbf{x}_k)-f(\mathbf{x}^*) = {\cal O}(\frac{1}{k})$, where the subscript $k$ is iteration index. This convergence slowdown is confirmed by the lower bound, which indicates that the number of FW subproblems to solve in order to ensure $f(\mathbf{x}_k) - f(\mathbf{x}^*) \leq \epsilon$, is no less than ${\cal O}\big(\frac{1}{\epsilon}\big)$ \citep{lan2013complexity,jaggi2013}. Thus, FW is a lower-bound-matching algorithm, in general. However, improved FW type algorithms are possible either in empirical performance, or, in speedup rates for certain subclasses of problems. Next, we deal with these improved rates paying attention to whether implementation requires knowing parameters such as the smoothness constant or the diameter of ${\cal X}$.

\textit{Parameter-dependent FW with faster rates.} This class of algorithms utilizes parameters that are obtained for different instances of $f$ and ${\cal X}$. Depending on the needed parameters, these algorithms are further classified into: i) line search based FW; ii) shorter step size aided FW; and iii) conditional gradient sliding (CGS). Line search based FW relies on $f(\mathbf{x})$ evaluations, which renders inefficiency when acquisition of function values is costly. The vanilla FW with line search converges with rate ${\cal O}(\frac{1}{k})$ on general problems \citep{jaggi2013}. Jointly leveraging line search and `away steps,' variants of FW converge linearly for strongly convex problems when ${\cal X}$ is a polytope~\citep{guelat1986,lacoste2015}; see also \citep{pedregosa2018,braun2018blended}. To improve the memory efficiency of away steps, a variant is further developed in \citep{garber2016linear}. Shorter step sizes refer to those used in \citep{levitin1966,garber2015}, where the step size is obtained by minimizing a one-dimensional quadratic function over $[0,1]$. Shorter step sizes require the smoothness parameter, which needs to be estimated for different loss functions. If ${\cal X}$ is strongly convex, and the optimal solution is at the boundary of ${\cal X}$, it is known that FW converges linearly~\citep{levitin1966}. For uniformly (and thus strongly) convex sets, faster rates are attained given that the optimal solution is at the boundary of ${\cal X}$~\citep{kerdreux2020}. When both $f$ and ${\cal X}$ are strongly convex, FW with shorter step size converges at a rate of ${\cal O}(\frac{1}{k^2})$, regardless of where the optimal solution resides~\citep{garber2015}. The last category is CGS, where both smoothness parameter and the diameter of ${\cal X}$ are necessary. In CGS, the subproblem of the original NAG that relies on projection is replaced by gradient sliding that solves a sequence of FW subproblems. A faster rate ${\cal O}(\frac{1}{k^2})$ is obtained at the price of: i) requiring at most ${\cal O}(k)$ FW subproblems in the $k$th iteration; and ii) an inefficient implementation since the NAG subproblem has to be solved up to a certain accuracy.

\textit{Parameter-free FW.} The advantage of a parameter-free algorithm is its efficient implementation. Since no parameter is involved, there is no concern on the quality of parameter estimation. This also saves time and effort because the step sizes do not need tuning. Although implementation efficiency is ensured, theoretical guarantees are challenging to obtain. This is because $f(\mathbf{x}_{k+1}) \leq f(\mathbf{x}_k)$ cannot be guaranteed without line search or shorter step sizes. Faster rates for parameter-free FW are rather limited in number, and most of existing parameter-free FW approaches rely on diminishing step sizes at the order of ${\cal O}(\frac{1}{k})$. For example, the behavior of FW when $k$ is large and ${\cal X}$ is a polytope is investigated under strong assumptions on $f(\mathbf{x})$ to be twice differentiable and locally strongly convex around $\mathbf{x}^*$~\citep{bach2020}. AFW \citep{li2020} replaces the subproblem of NAG by a single FW subproblem, where constraint-specific faster rates are developed. Taking an active $\ell_2$ norm ball constraint as an example, AFW guarantees a rate of ${\cal O}\big( \frac{\ln k}{k^2} \big)$. A natural question is whether the $\ln k$ in the numerator can be eliminated. In addition, although the implementation involves no parameter, the analysis of AFW relies on the value $\max_{\mathbf{x} \in {\cal X}} f(\mathbf{x})$. 

Aiming at parameter-free FW with faster rates (on certain constraints) that can bypass the limitations of AFW, the present work deals with the design and analysis of ExtraFW. The `extra' in its name refers to the pair of gradients involved per iteration, whose merit is to enable a `prediction-correction' (PC) type of update. Though the idea of using two gradients to perform PC updates originates from projection-based algorithms, such as ExtraGradient \citep{korpelevich1976} and Mirror-Prox \citep{nemirovski2004,diakonikolas2017,kavis2019}, leveraging PC updates in FW type algorithms for faster rates is novel.  

Our contributions are summarized as follows.

\begin{enumerate}
	\item[\textbullet] A new parameter-free FW variant, ExtraFW, is studied in this work. The distinct feature of ExtraFW is the adoption of two gradient evaluations per iteration to update the decision variable in a prediction-correction (PC) manner. 
	
\item[\textbullet] It is shown that ExtraFW convergences with a rate of ${\cal O}(\frac{1}{k})$ for general problems. And for constraint sets including active $\ell_1$, $\ell_2$ and $n$-support norm balls, ExtraFW guarantees an accelerated rate ${\cal O}(\frac{1}{k^2})$. 

\item[\textbullet] Unlike most of faster rates in FW literatures, ExtraFW is parameter-free, so that no problem dependent parameter is required. Compared with another parameter-free algorithm with faster rates, AFW \citep{li2020}, introducing PC update in ExtraFW leads to several advantages: i) the convergence rate is improved by a factor of ${\cal O}(\ln k)$ on an $\ell_2$ norm ball constraint; and ii) the analysis does not rely on the maximum value of $f(\mathbf{x})$ over ${\cal X}$. 

\item[\textbullet] The efficiency of ExtraFW is corroborated on two benchmark machine learning tasks. The faster rate ${\cal O}(\frac{1}{k^2})$ is achieved on binary classification, evidenced by the possible improvement of ExtraFW over NAG on multiple sparsity-promoting constraint sets. For matrix completion, ExtraFW improves over AFW and FW in both optimality error and the rank of the solution.

\end{enumerate}

\textbf{Notation}. Bold lowercase (uppercase) letters denote vectors (matrices); $\| \mathbf{x}\|$ stands for a norm of $\mathbf{x}$, with its dual norm written as $\|\mathbf{x} \|_*$; and $\langle \mathbf{x}, \mathbf{y} \rangle$ denotes the inner product of $\mathbf{x}$ and $\mathbf{y}$. We also define $x \wedge y:= \min\{ x, y\}$. 

\section{Preliminaries}

This section reviews FW and AFW in order to illustrate the proposed algorithm in a principled manner. We first pinpoint the class of problems to focus on. 

\begin{assumption}\label{as.1}
	(Lipschitz Continuous Gradient.)
	The function $f: {\cal X} \rightarrow \mathbb{R}$ has $L$-Lipchitz continuous gradients; that is, $\|\nabla f(\mathbf{x}) \!-\! \nabla f(\mathbf{y}) \|_* \leq L \| \mathbf{x}-\mathbf{y} \|, \forall\, \mathbf{x}, \mathbf{y} \in {\cal X}$.
\end{assumption} 

\begin{assumption}\label{as.2}
	(Convex Objective Function.)
	The function $f: {\cal X} \rightarrow \mathbb{R}$ is convex; that is, $f(\mathbf{y}) - f(\mathbf{x}) \geq \langle \nabla f(\mathbf{x}), \mathbf{y} - \mathbf{x} \rangle, \forall\, \mathbf{x}, \mathbf{y} \in {\cal X}$.
\end{assumption} 
\begin{assumption}\label{as.3}
	(Constraint Set.)
	The constraint set $\cal X$ is convex and compact with diameter $D$, that is, $\| \mathbf{x} - \mathbf{y} \| \leq D, \forall \mathbf{x}, \mathbf{y} \in {\cal X}$.
\end{assumption} 
Assumptions \ref{as.1} -- \ref{as.3} are standard for FW type algorithms, and will be taken to hold true throughout. A blackbox optimization paradigm is considered in this work, where the objective function and constraint set can be accessed through oracles only. In particular, the first-order oracle (FO) and the linear minimization oracle (LMO) are needed.
\begin{definition}
	(FO.) The first-order oracle takes $\mathbf{x} \in {\cal X}$ as an input and returns its gradient $\nabla f(\mathbf{x})$.
\end{definition}
\begin{definition}
	(LMO.) The linear minimization oracle takes a vector $\mathbf{g} \in \mathbb{R}^d$ as an input and returns a minimizer of $\min_{\mathbf{x}\in {\cal X}} \langle \mathbf{g}, \mathbf{x} \rangle$.
\end{definition}
Except for gradients, problem dependent parameters such as function value, smoothness constant $L$, and constraint diameter $D$ are not provided by FO and LMO. Hence, algorithms relying only on FO and LMO are parameter-free. Next, we recap FW and AFW with parameter-free step sizes to gain more insights for the proposed algorithm.

\textbf{FW recap.}
FW is summarized in Alg. \ref{alg.fw}. 
A subproblem with a linear loss, referred to also as an \textit{FW step}, is solved per iteration via LMO. The FW step can be explained as finding a minimizer over ${\cal X}$ for the following supporting hyperplane of $f(\mathbf{x})$, 
\begin{align}\label{eq.fw_lb}
	f(\mathbf{x}_k) + \langle \nabla f(\mathbf{x}_k), \mathbf{x} - \mathbf{x}_k \rangle.
\end{align}
Note that \eqref{eq.fw_lb} is also a lower bound for $f(\mathbf{x})$ due to convexity.
Upon obtaining $\mathbf{v}_{k+1}$ by minimizing \eqref{eq.fw_lb} over ${\cal X}$, $\mathbf{x}_{k+1}$ is updated as a convex combination of $\mathbf{v}_{k+1}$ and $\mathbf{x}_k$ to eliminate the projection. 
The parameter-free step size is usually chosen as $\delta_k = \frac{2}{k+2}$. As for convergence, FW guarantees $f(\mathbf{x}_k) - f(\mathbf{x}^*) = {\cal O}(\frac{LD^2}{k})$. 
 
\begin{algorithm}[t]
    \caption{FW \citep{frank1956}}\label{alg.fw}
    \begin{algorithmic}[1]
    	\State \textbf{Initialize:} $\mathbf{x}_0\in {\cal X}$
    	\For {$k=0,1,\dots,K-1$}
    		\State $\mathbf{v}_{k+1} = \argmin_{\mathbf{x} \in \cal X} \langle \nabla f(\mathbf{x}_k), \mathbf{x} \rangle$ 
			\State $\mathbf{x}_{k+1} = (1-\delta_k) \mathbf{x}_k + \delta_k \mathbf{v}_{k+1}$ 
		\EndFor
		\State \textbf{Return:} $\mathbf{x}_K$
	\end{algorithmic}
\end{algorithm}

\textbf{AFW recap.} As an FW variant, AFW in Alg. \ref{alg.afw} relies on Nesterov momentum type update, that is, it uses an auxiliary variable $\mathbf{y}_k$ to estimate $\mathbf{x}_{k+1}$ and calculates the gradient $\nabla f(\mathbf{y}_k)$. If one writes $\mathbf{g}_{k+1}$ explicitly, $\mathbf{v}_{k+1}$ can be equivalently described as a minimizer over ${\cal X}$ of the hyperplane
\begin{align}\label{eq.afw_lb}
	\sum_{\tau=0}^k w_k^\tau \big[ f(\mathbf{y}_\tau) + \langle \nabla f(\mathbf{y}_\tau), \mathbf{x} -  \mathbf{y}_\tau \rangle \big]
\end{align}
where $w_k^\tau = \delta_\tau  \prod_{j = \tau+1}^k   (1-\delta_j)$ and $\sum_{\tau=0}^k w_k^\tau \approx 1$ (the sum depends on the choice of $\delta_0$). Note that $f(\mathbf{y}_\tau) + \langle \nabla f(\mathbf{y}_\tau), \mathbf{x} - \mathbf{y}_\tau \rangle$ is a supporting hyperplane of $f(\mathbf{x})$ at $\mathbf{y}_\tau$, hence \eqref{eq.afw_lb} is a lower bound for $f(\mathbf{x})$ constructed through a weighted average of supporting hyperplanes at $\{\mathbf{y}_\tau\}$. AFW converges at ${\cal O}\big(\frac{LD^2}{k}\big)$ on general problems. When the constraint set is an active $\ell_2$ norm ball, AFW has a faster rate ${\cal O}\big( \frac{LD^2}{k} \wedge \frac{T LD^2 \ln k}{k^2}\big)$, where $T$ depends on $D$. Writing this rate compactly as ${\cal O}\big(\frac{T LD^2 \ln k}{k^2}\big)$, it is observed that AFW achieves acceleration with the price of a worse dependence on other parameters hidden in $T$. However, even for the $k$-dependence, AFW is ${\cal O}(\ln k)$ times slower compared with other momentum based algorithms such as NAG. This slowdown is because that the lower bound \eqref{eq.afw_lb} is constructed based on $\{\mathbf{y}_k\}$, which are estimated $\{\mathbf{x}_{k+1}\}$. We will show that relying on a lower bound constructed using $\{\mathbf{x}_{k+1}\}$ directly, it is possible to avoid this ${\cal O}(\ln k)$ slowdown.

%

\begin{algorithm}[t]
    \caption{AFW \citep{li2020}}\label{alg.afw}
    \begin{algorithmic}[1]
    	\State \textbf{Initialize:} $\mathbf{x}_0\in {\cal X}$, $\mathbf{g}_0 = \mathbf{0}$, $\mathbf{v}_0 = \mathbf{x}_0$
    	\For {$k=0,1,\dots,K-1$}
    		\State $\mathbf{y}_k = (1-\delta_k) \mathbf{x}_k + \delta_k \mathbf{v}_k$
    		\State $\mathbf{g}_{k+1} = (1-\delta_k) \mathbf{g}_k + \delta_k \nabla f(\mathbf{y}_k)$
    		\State $\mathbf{v}_{k+1} = \argmin_{\mathbf{x} \in \cal X} \langle  \mathbf{g}_{k+1}, \mathbf{x} \rangle$
			\State $\mathbf{x}_{k+1} = (1-\delta_k) \mathbf{x}_k + \delta_k \mathbf{v}_{k+1}$ 
		\EndFor
		\State \textbf{Return:} $\mathbf{x}_K$
	\end{algorithmic}
\end{algorithm}


\section{ExtraFW}\label{sec.exfw}

This section introduces the main algorithm, ExtraFW, and establishes its constraint dependent faster rates. 

\subsection{Algorithm Design}
\begin{algorithm}[t]
    \caption{ExtraFW}\label{alg.exfw}
    \begin{algorithmic}[1]
    	\State \textbf{Initialize:} $\mathbf{x}_0 $, $\mathbf{g}_0 = \mathbf{0}$, and $\mathbf{v}_0 = \mathbf{x}_0$
			\For {$k=0,1,\dots,K-1$}
				\State $\mathbf{y}_k =  (1-\delta_k) \mathbf{x}_k + \delta_k \mathbf{v}_k $ \Comment{prediction}
				\State $\hat{\mathbf{g}}_{k+1} = (1-\delta_k) \mathbf{g}_k + \delta_k \nabla f(\mathbf{y}_k)$ 
				\State $\hat{\mathbf{v}}_{k+1} = \argmin_{\mathbf{v} \in {\cal X}} \langle  \hat{\mathbf{g}}_{k+1}, \mathbf{v}  \rangle $	
				\State $\mathbf{x}_{k+1} =  (1-\delta_k) \mathbf{x}_k + \delta_k \hat{\mathbf{v}}_{k+1} $ \Comment{correction}
				\State $\mathbf{g}_{k+1} =  (1-\delta_k) \mathbf{g}_k +  \delta_k \nabla f(\mathbf{x}_{k+1}) $ 
				\State $\mathbf{v}_{k+1} = \argmin_{\mathbf{v} \in {\cal X}} \langle  \mathbf{g}_{k+1}, \mathbf{v}  \rangle $		\Comment{extra FW step}		
			\EndFor
		\State \textbf{Return:} $\mathbf{x}_K$
	\end{algorithmic}
\end{algorithm}

ExtraFW is summarized in Alg. \ref{alg.exfw}. Different from the vanilla FW and AFW, two FW steps (Lines 5 and 8 of Alg. \ref{alg.exfw}) are required per iteration. Compared with other algorithms relying on two gradient evaluations, such as Mirror-Prox \citep{diakonikolas2017,kavis2019}, ExtraFW reduces the computational burden of the projection. In addition, as an FW variant, ExtraFW can capture the properties such as sparsity or low rank promoted by the constraints more effectively through the update than those projection based algorithms. Detailed elaboration can be found in Section \ref{sec.numerical} and Appendix \ref{apdx.numerical}. To facilitate comparison with FW and AFW, ExtraFW is explained through constructing lower bounds of $f(\mathbf{x})$ in a ``prediction-correction'' manner. The merits of the PC update compared with AFW are: i) the elimination of $\max_{\mathbf{x} \in {\cal X}} f(\mathbf{x})$ in analysis; and ii) it improves the convergence rate on certain class of problems as we will see later.

\textbf{Lower bound prediction.} Similar to AFW, the auxiliary variable $\mathbf{y}_k$ in Line 3 of Alg.~\ref{alg.exfw} can be viewed as an estimate of $\mathbf{x}_{k+1}$. The first gradient is evaluated at $\mathbf{y}_k$, and is incorporated into $\hat{\mathbf{g}}_{k+1}$, which is an estimate of the weighted average of $\{ \nabla f(\mathbf{x})_\tau \}_{\tau=1}^{k+1}$. By expanding $\hat{\mathbf{g}}_{k+1}$, one can verify that $\hat{\mathbf{v}}_{k+1}$ can be obtained equivalently through minimizing the following weighted sum, 
\begin{align}\label{eq.exfw_lb_est}
	\sum_{\tau=0}^{k-1}  w_k^\tau \Big[ f(\mathbf{x}_{\tau+1}) + \langle \nabla f(\mathbf{x}_{\tau+1}), \mathbf{x} - \mathbf{x}_{\tau+1} \rangle \Big] + \delta_k \Big[ f(\mathbf{y}_k) +  \big\langle \nabla f(\mathbf{y}_k) , \mathbf{x}  - \mathbf{y}_k \big\rangle  \Big], 
\end{align}
where $w_\tau =\delta_\tau \prod_{j = \tau+1}^k (1-\delta_j) $ and $\sum_{\tau=0}^{k-1} w_\tau + \delta_k \approx 1$. Note that each term inside square brackets forms a supporting hyperplane of $f(\mathbf{x})$, hence \eqref{eq.exfw_lb_est} is an (approximated) lower bound of $f(\mathbf{x})$ because of convexity. As a prediction to $f(\mathbf{x}_{k+1}) +  \langle \nabla f(\mathbf{x}_{k+1}) , \mathbf{x}  - \mathbf{x}_{k+1} \rangle$, the last bracket in \eqref{eq.exfw_lb_est} will be corrected once $\mathbf{x}_{k+1}$ is obtained.

\textbf{Lower bound correction.} The gradient $\nabla f(\mathbf{x}_{k+1})$ is used to obtain a weighted averaged gradients  $\mathbf{g}_{k+1}$. By unrolling $\mathbf{g}_{k+1}$, one can find that $\mathbf{v}_{k+1}$ is a minimizer of the following (approximated) lower bound of $f(\mathbf{x})$
\begin{align}\label{eq.exfw_lb_correct}
	\sum_{\tau=0}^{k-1} & w_k^\tau \Big[ f(\mathbf{x}_{\tau+1}) + \big\langle \nabla f(\mathbf{x}_{\tau+1}), \mathbf{x} - \mathbf{x}_{\tau+1} \big\rangle \Big] + \delta_k \Big[ f(\mathbf{x}_{k+1}) +  \big\langle \nabla f(\mathbf{x}_{k+1}) , \mathbf{x}  - \mathbf{x}_{k+1} \big\rangle  \Big].
\end{align}
Comparing \eqref{eq.exfw_lb_est} and \eqref{eq.exfw_lb_correct}, we deduce that the terms in the last bracket of \eqref{eq.exfw_lb_est} are corrected to the true supporting hyperplane of $f(\mathbf{x})$ at $\mathbf{x}_{k+1}$. In sum, the FW steps in ExtraFW rely on lower bounds of $f(\mathbf{x})$ constructed in a weighted average manner similar to AFW. However, the key difference is that ExtraFW leverages the supporting hyperplanes at true variables $\{ \mathbf{x}_k\}$ rather than the auxiliary ones $\{ \mathbf{y}_k\}$ in AFW through a ``correction'' effected by \eqref{eq.exfw_lb_correct}. In the following subsections, we will show that the PC update in ExtraFW performs no worse than FW or AFW on general problems, while harnessing its own analytical merits on certain constraint sets.

\subsection{Convergence of ExtraFW}\label{sec.conv_gen}
We investigate the convergence of ExtraFW by considering the general case first. The analysis relies on the notion of estimate sequence (ES) introduced in \citep{nesterov2004}. An ES ``estimates'' $f$ using a sequence of surrogate functions $\{\Phi_k(\mathbf{x}) \}$ that are analytically tractable (e.g., being quadratic or linear). ES is formalized in the following definition.
\begin{definition}\label{def.es}
	A tuple $\big( \{\Phi_k(\mathbf{x}) \}_{k=0}^\infty,\{ \lambda_k \}_{k=0}^\infty \big)$ is called an estimate sequence of function $f(\mathbf{x})$ if $\lim_{k\rightarrow \infty} \lambda_k = 0$ and for any $\mathbf{x}\in {\cal X }$ we have $
		\Phi_k (\mathbf{x}) \leq (1 - \lambda_k)	 f(\mathbf{x}) + \lambda_k \Phi_0 (\mathbf{x})$.
\end{definition}

The construction of ES varies for different algorithms (see e.g., \citep{kulunchakov2019,nesterov2004,lin2015,li2019bb}). However, the reason to rely on the ES based analysis is similar, as summarized in the following lemma. 
\begin{lemma}\label{lemma.es_convergence}
	For $\big( \{\Phi_k(\mathbf{x}) \}_{k=0}^\infty,\{ \lambda_k \}_{k=0}^\infty \big)$ satisfying the definition of ES, if $f(\mathbf{x}_k) \leq \min_{\mathbf{x} \in {\cal X}} \Phi_k(\mathbf{x}) + \xi_k, \forall k$, it is true that
	\begin{align*}
		f(\mathbf{x}_k) - f(\mathbf{x}^*) \leq \lambda_k \big( \Phi_0 (\mathbf{x}^*) - f(\mathbf{x}^*) \big) + \xi_k , \forall\, k.
	\end{align*}
\end{lemma}

As shown in Lemma \ref{lemma.es_convergence}, $\lambda_k$ and $\xi_k$ jointly characterize the convergence rate of $f(\mathbf{x}_k)$. (Consider $\lambda_k = {\cal O}(\frac{1}{k})$ and $\xi_k = {\cal O}(\frac{1}{k})$ for an example.) Keeping Lemma \ref{lemma.es_convergence} in mind, we construct \textit{two} sequences of \textit{linear} surrogate functions for analyzing ExtraFW, which highlight the differences of our analysis with existing ES based approaches
\begin{subequations}\label{eq.phi}
\begin{align}
	\Phi_0(\mathbf{x}) & = \hat{\Phi}_0 (\mathbf{x}) \equiv f(\mathbf{x}_0) \\
	\hat{\Phi}_{k+1}(\mathbf{x}) &= (1 - \delta_k) \Phi_k(\mathbf{x}) + \delta_k \big[ f(\mathbf{y}_k)  +  \langle \nabla f(\mathbf{y}_k)   , \mathbf{x}  - \mathbf{y}_k \rangle  \big] , ~\forall k \geq 0 \\
	\Phi_{k+1}(\mathbf{x}) &= (1 - \delta_k) \Phi_k(\mathbf{x})  +  \langle \nabla f(\mathbf{x}_{k+1}) , \mathbf{x}  - \mathbf{x}_{k+1} \rangle  \big] , ~\forall k \geq 0.  
\end{align}
\end{subequations}

Clearly, both $\Phi_k(\mathbf{x})$ and $\hat{\Phi}_k(\mathbf{x})$ are linear in $\mathbf{x}$, in contrast to the quadratic surrogate functions adopted for analyzing NAG \citep{nesterov2004}. Such linear surrogate functions are constructed specifically for FW type algorithms taking advantage of the compact and convex constraint set. Next we show that \eqref{eq.phi} and proper $\{\lambda_k\}$ form two different ES of $f$. 

\begin{lemma}\label{lemma.es}
	If we choose $\lambda_0 = 1$, $\delta_k \in (0,1)$, and $\lambda_{k+1} = (1-\delta_k) \lambda_k ~\forall k \geq 0$, both $\big( \{\Phi_k(\mathbf{x}) \}_{k=0}^\infty,$$ \{ \lambda_k \}_{k=0}^\infty \big)$ and $\big( \{\hat{\Phi}_k(\mathbf{x}) \}_{k=0}^\infty, \{ \lambda_k \}_{k=0}^\infty \big)$ satisfy the definition of ES.
\end{lemma}

The key reason behind the construction of surrogate functions in \eqref{eq.phi} is that they are closely linked with the lower bounds \eqref{eq.exfw_lb_est} and \eqref{eq.exfw_lb_correct} used in the FW steps, as stated in the next lemma.

\begin{lemma}\label{lemma.vstar}
	Let $\mathbf{g}_0 = \mathbf{0}$, then it is true that $ \mathbf{v}_k = \argmin_{\mathbf{x} \in {\cal X}} \Phi_k (\mathbf{x})$ and $\hat{\mathbf{v}}_k = \argmin_{\mathbf{x} \in {\cal X}} \hat{\Phi}_k (\mathbf{x})$.
\end{lemma}
After relating the surrogate functions in \eqref{eq.phi} with ExtraFW, exploiting the analytical merits of the surrogate functions $\Phi_k(\mathbf{x})$ and $\hat{\Phi}_k(\mathbf{x})$, including being linear, next we show that $f(\mathbf{x}_k) \leq \min_{\mathbf{x} \in {\cal X}} \Phi_k(\mathbf{x}) + \xi_k, \forall k$, which is the premise of Lemma \ref{lemma.es_convergence}.

\begin{lemma}\label{lemma_phistar}
	Let $\xi_0 = 0$ and other parameters chosen the same as previous lemmas. Denote $\Phi_k^* := \Phi_k(\mathbf{v}_k)$ as the minimum value of $\Phi_k(\mathbf{x})$ over ${\cal X}$ (cf. Lemma \ref{lemma.vstar}), then ExtraFW guarantees that for any $k\geq 0$
	\begin{align*}
		f(\mathbf{x}_k )\leq \Phi_k^* + \xi_k, ~\text{with}~ \xi_{k+1} = (1-\delta_k) \xi_k + \frac{3 L D^2}{2} \delta_k^2. 
	\end{align*}
\end{lemma}
Based on Lemma \ref{lemma_phistar}, the value of $f(\mathbf{x}_k )$ and $\Phi_k^*$ can be used to derive the stopping criterion if one does not want to preset the iteration number $K$. Further discussions are provided in Appendix \ref{apdx.stop} due to space limitation. Now we are ready to apply Lemma \ref{lemma.es_convergence} to establish the convergence of ExtraFW. 

\begin{theorem}\label{thm.general}
Suppose that Assumptions \ref{as.1}, \ref{as.2} and \ref{as.3} are satisfied. Choosing $\delta_k = \frac{2}{k+3}$, and $\mathbf{g}_0 = \mathbf{0}$, ExtraFW in Alg. \ref{alg.exfw} guarantees	
	\begin{align*}
		f(\mathbf{x}_k) - f(\mathbf{x}^*) = {\cal O} \bigg( \frac{LD^2}{k} \bigg), \forall k.
	\end{align*}
\end{theorem}

This convergence rate of ExtraFW has the same order as AFW and FW. In addition, Theorem \ref{thm.general} translates into ${\cal O}(\frac{LD^2}{\epsilon})$ queries of LMO to ensure $f(\mathbf{x}_k) - f(\mathbf{x}^*) \leq \epsilon$, which matches to the lower bound \citep{lan2013complexity,jaggi2013}.

\textbf{The obstacle for faster rates.} As shown in the detailed proof, one needs to guarantee that either $\| \mathbf{v}_k - \hat{\mathbf{v}}_{k+1} \|^2 $ or $\| \mathbf{v}_{k+1} - \hat{\mathbf{v}}_{k+1} \|^2 $ is small enough to obtain a faster rate than Theorem \ref{thm.general}. This is difficult in general because there could be multiple $\mathbf{v}_k$ and $\hat{\mathbf{v}}_k$ solving the FW steps. A simple example is to consider the $i$th entry  $[\mathbf{g}_k]_i=0$. The $i$th entry $[\mathbf{v}_k]_i$ can then be chosen arbitrarily as long as $\mathbf{v}_k \in {\cal X}$. The non-uniqueness of $\mathbf{v}_k$ prevents one from ensuring a small upper bound of $\|\mathbf{v}_k -\hat{\mathbf{v}}_{k+1}\|^2,~\forall~\mathbf{v}_k$. In spite of this, we will show that together with the structure on ${\cal X}$, ExtraFW can attain faster rates.

\subsection{Acceleration of ExtraFW}\label{sec.acc}
In this subsection, we provide constraint-dependent accelerated rates of ExtraFW when ${\cal X}$ is some norm ball. Even for projection based algorithms, most of faster rates are obtained with step sizes depending on $L$ \citep{nemirovski2004,diakonikolas2017}. Thus, faster rates for parameter-free algorithms are challenging to establish. An extra assumption is needed in this subsection.

\begin{assumption}\label{as.4}
	The constraint is active, i.e., $\| \nabla f(\mathbf{x}^*)\|_2 \geq G >0$.
\end{assumption}
It is natural to rely on the position of the optimal solution in FW type algorithms for analysis, and one can see this assumption also in \citep{levitin1966,dunn1979,li2020,kerdreux2020}. For a number of machine learning tasks, Assumption \ref{as.4} is rather mild. Relying on Lagrangian duality, it can be seen that problem \eqref{eq.prob} with a norm ball constraint is equivalent to the regularized formulation $\min_{\mathbf{x}} f(\mathbf{x})+ \gamma g(\mathbf{x})$, where $\gamma \geq 0$ is the Lagrange multiplier, and $g(\mathbf{x})$ denotes some norm. In view of this, Assumption \ref{as.4} simply implies that $\gamma > 0$ in the equivalent regularized formulation, that is, the norm ball constraint plays the role of a regularizer. Given the prevalence of the regularized formulation in machine learning, it is worth investigating its equivalent constrained form \eqref{eq.prob} under Assumption \ref{as.4}.

Technically, the need behind Assumption \ref{as.4} can be exemplified through a one-dimensional problem. Consider minimizing $f(x) = x^2$ over ${\cal X}= \{x| x \in [-1,1]\}$. We clearly have $x^* = 0$ for which the constraint is inactive at the optimal solution. Recall a faster rate of ExtraFW requires $\| \hat{v}_{k+1} - v_{k+1}\|_2$ to be small. When $x_k$ is close to $x^* = 0$, it can happen that $\hat{g}_{k+1}>0$ and $g_{k+1}<0$, leading to $\hat{v}_{k+1} = -1$ and $v_{k+1} = 1$. The faster rate is prevented by pushing $v_{k+1}$ and $\hat{v}_{k+1}$ further apart from each other. 

Next, we consider different instances of norm ball constraints as examples to the acceleration of ExtraFW. For simplicity of exposition, the intuition and technical details are discussed using an $\ell_2$ norm ball constraint in the main test. Detailed analysis for $\ell_1$ and $n$-support norm ball \citep{argyriou2012} constraints are provided in Appendix.

\textbf{$\ell_2$ norm ball constraint.} Consider ${\cal X}:= \{ \mathbf{x}| \| \mathbf{x} \|_2 \leq \frac{D}{2} \}$. In this case, $\mathbf{v}_{k+1}$ and $\hat{\mathbf{v}}_{k+1}$ admit closed-form solutions, taking $\mathbf{v}_{k+1}$ as an example,
\begin{align}\label{eq.acc_opt_v}
	\mathbf{v}_{k+1} = \argmin_{\mathbf{x}\in{\cal X}} \langle \mathbf{g}_{k+1}, \mathbf{x} \rangle = - \frac{D}{2 \| \mathbf{g}_{k+1} \|_2} \mathbf{g}_{k+1}.
\end{align}
We assume that when using $\mathbf{g}_{k+1}$ as the input to the LMO, the returned vector is given by \eqref{eq.acc_opt_v}. This is reasonable since it is what we usually implemented in practice. Though rarely happen, one can choose $\mathbf{v}_{k+1} = \hat{\mathbf{v}}_{k+1}$ to proceed if $\mathbf{g}_{k+1} = \mathbf{0}$. Similarly, we can simply set $ \hat{\mathbf{v}}_{k+1} = \mathbf{v}_k$ if $\hat{\mathbf{g}}_{k+1} = \mathbf{0}$. The uniqueness of $\mathbf{v}_{k+1}$ is ensured by its closed-form solution, wiping out the obstacle for a faster rate.


\begin{theorem}\label{thm.acc}
	Suppose that Assumptions \ref{as.1}, \ref{as.2}, \ref{as.3} and \ref{as.4} are satisfied, and ${\cal X}$ is an $\ell_2$ norm ball. Choosing $\delta_k = \frac{2}{k+3}$, and $\mathbf{g}_0 = \mathbf{0}$, ExtraFW in Alg. \ref{alg.exfw} guarantees
	\begin{align*}
		f(\mathbf{x}_k) - f(\mathbf{x}^*) = {\cal O} \bigg( \frac{LD^2}{k} \wedge \frac{L D^2 T }{k^2} \bigg), \forall k
	\end{align*}
	where $T$ is a constant depending only on $L$, $G$, and $D$.
\end{theorem}
Theorem \ref{thm.acc} admits a couple of interpretations. By writing the rate compactly, ExtraFW achieves accelerated rate ${\cal O}\big(  \frac{TLD^2}{k^2} \big), \forall k$ with a worse dependence on $D$ compared to the vanilla FW. Or alternatively, the ``asymptotic'' performance at $k \geq T$ is strictly improved over the vanilla FW. It is worth mentioning 
that the choices of $\delta_k$ and $\mathbf{g}_0$ are not changed compared to Theorem \ref{thm.general} so that the parameter-free implementation is the same regardless whether accelerated. In other words, prior knowledge on whether Assumption \ref{as.4} holds is not needed in practice. Compared with CGS, ExtraFW sacrifices the $D$ dependence in the convergence rate to trade for i) the nonnecessity of the knowledge of $L$ and $D$, and ii) ensuring two FW subproblems per iteration (whereas at most ${\cal O}(k)$ subproblems are needed in CGS). When comparing with AFW \citep{li2020}, the convergence rate of ExtraFW is improved by a factor of ${\cal O}(\ln k)$, and the analysis does not rely on the constant $M:= \max_{\mathbf{x}\in {\cal X}} f(\mathbf{x})$. 

\textbf{$\ell_1$ norm ball constraint.} For the sparsity-promoting constraint ${\cal X}:= \{ \mathbf{x}| \| \mathbf{x} \|_1 \leq R \}$, the FW steps can be solved in closed form too. Taking $\mathbf{v}_{k+1}$ as an example, we have 
\begin{align}\label{eq.fw_l1}
	\mathbf{v}_{k+1} = R  \cdot [0, & \ldots, 0,  -{\rm sgn} [\mathbf{g}_{k+1}]_i ,0,\ldots, 0]^\top  ~~  {\rm with} ~~i = \argmax_j | [\mathbf{g}_{k+1}]_j |.
\end{align}
We show in Theorem \ref{thm.l1} (see Appendix \ref{apdx.ell1}) that when Assumption \ref{as.4} holds and the set $\argmax_j \big{|}[ \nabla f(\mathbf{x}^*) ]_j \big{|}$ has cardinality $1$, a faster rate ${\cal O}(\frac{T_1 LD^2}{k^2})$ can be obtained with the constant $T_1$ depending on $L$, $G$, and $D$. 
The additional assumption here is known as \textit{strict complementarity}, and has been adopted also in, e.g.,\citep{ding2020spectral,garber2020}.

\textbf{$n$-support norm ball constraint.} The $n$-support norm ball is a tighter relaxation of a sparsity prompting $\ell_0$ norm ball combined with an $\ell_2$ norm penalty compared with the ElasticNet \citep{zou2005}. It is defined as ${\cal X}:= {\rm conv} \{\mathbf{x}| \|\mathbf{x} \|_0 \leq n, \| \mathbf{x} \|_2 \leq R \}$, where ${\rm conv}\{\cdot\}$ denotes the convex hull \citep{argyriou2012}. The closed-form solution of $\mathbf{v}_{k+1}$ is given by \citep{liu2016efficient}
\begin{align}\label{eq.fw_opt_v}
	\mathbf{v}_{k+1} = - \frac{R}{\| {\rm top}_n (\mathbf{g}_{k+1}) \|_2} {\rm top}_n (\mathbf{g}_{k+1})
\end{align}
where ${\rm top}_n (\mathbf{g})$ denotes the truncated version of $\mathbf{g}$ with its top $n$ (in magnitude) entries preserved. A faster rate ${\cal O}(\frac{T_2 LD^2}{k^2})$ is guaranteed by ExtraFW under Assumption \ref{as.4}, and a condition similar to strict complementarity (see Theorem \ref{thm.nsupp} in the Appendix \ref{apdx.nsupp}). Again, the constant $T_2$ here depends on $L$, $G$, and $D$.

\begin{figure*}[t]
	\vspace{-0.1cm}
	\centering
	\begin{tabular}{cccc}
		\hspace{-0.2cm}
		\includegraphics[width=.24\textwidth]{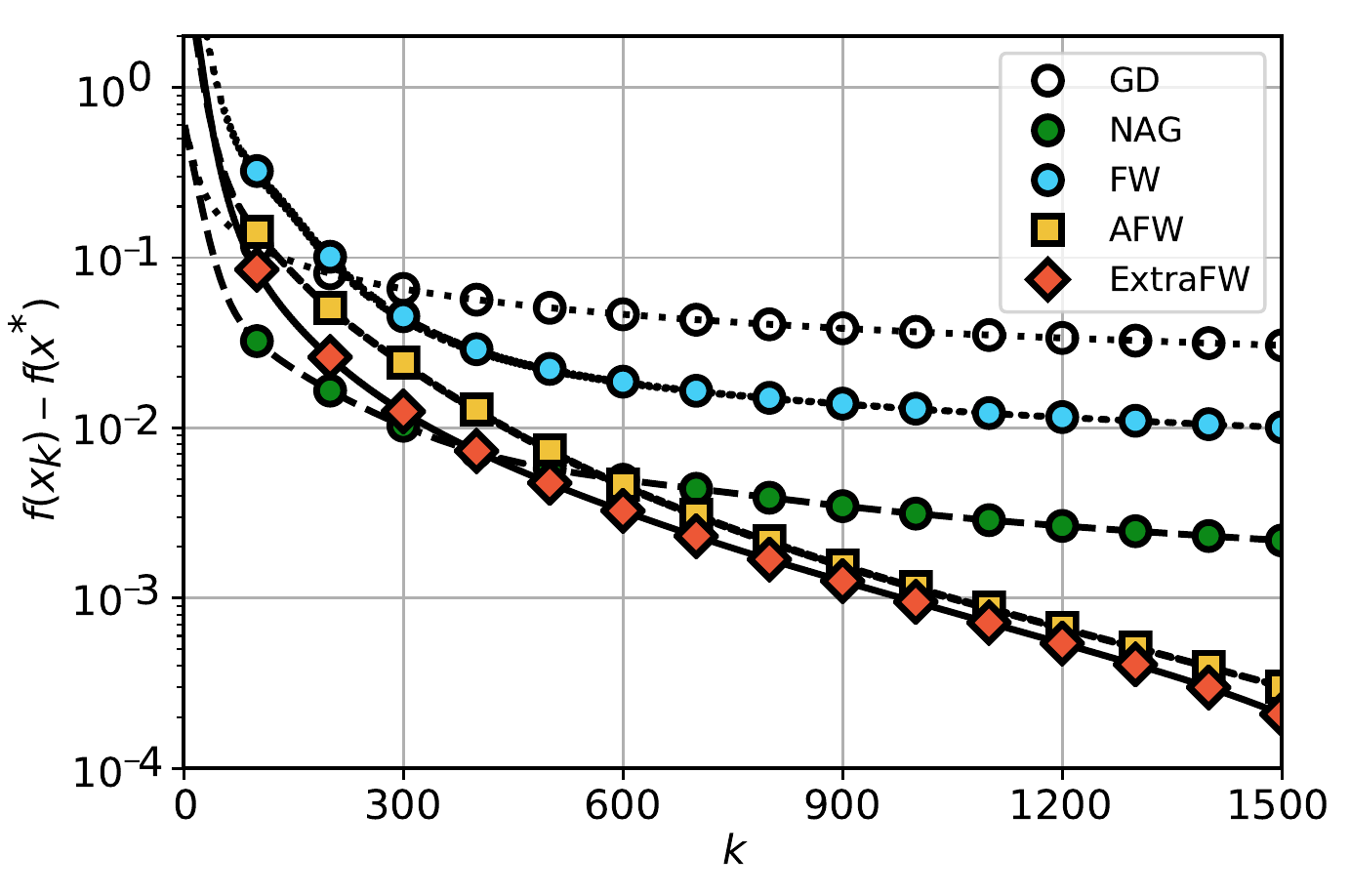}&
		\hspace{-0.4cm}
		\includegraphics[width=.24\textwidth]{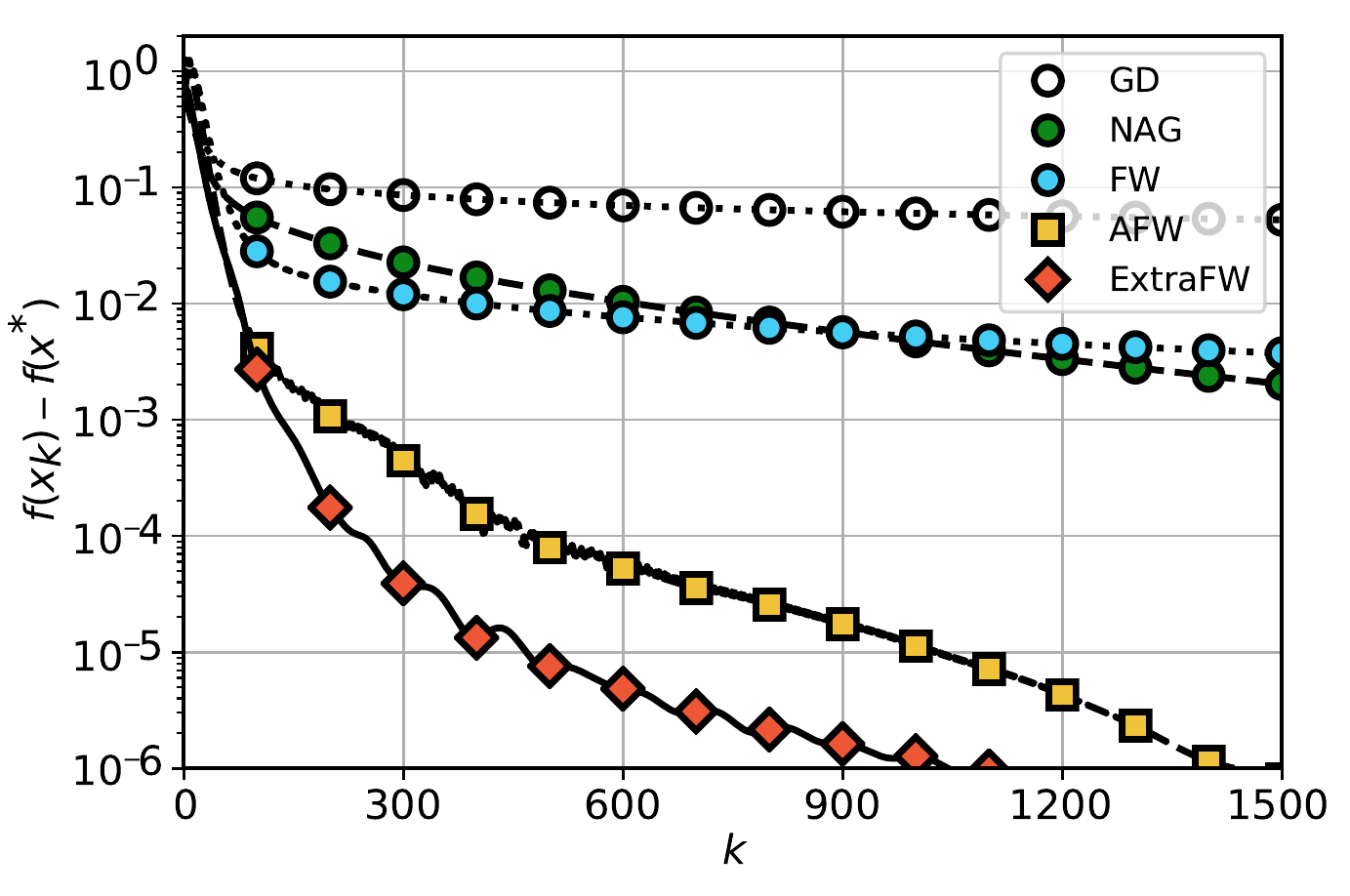}&
		\hspace{-0.4cm}
		\includegraphics[width=.24\textwidth]{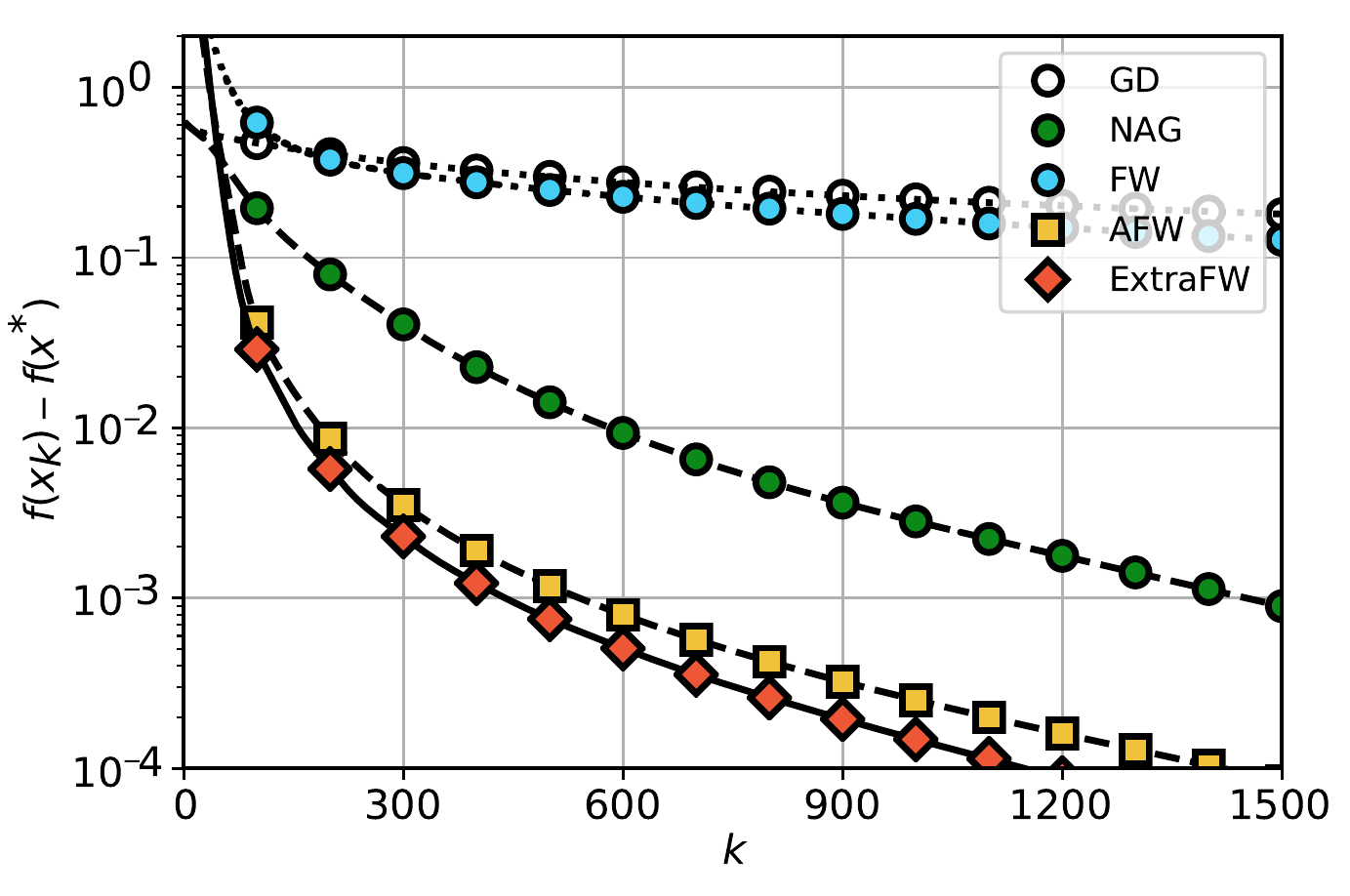}&
		\hspace{-0.4cm}
		\includegraphics[width=.24\textwidth]{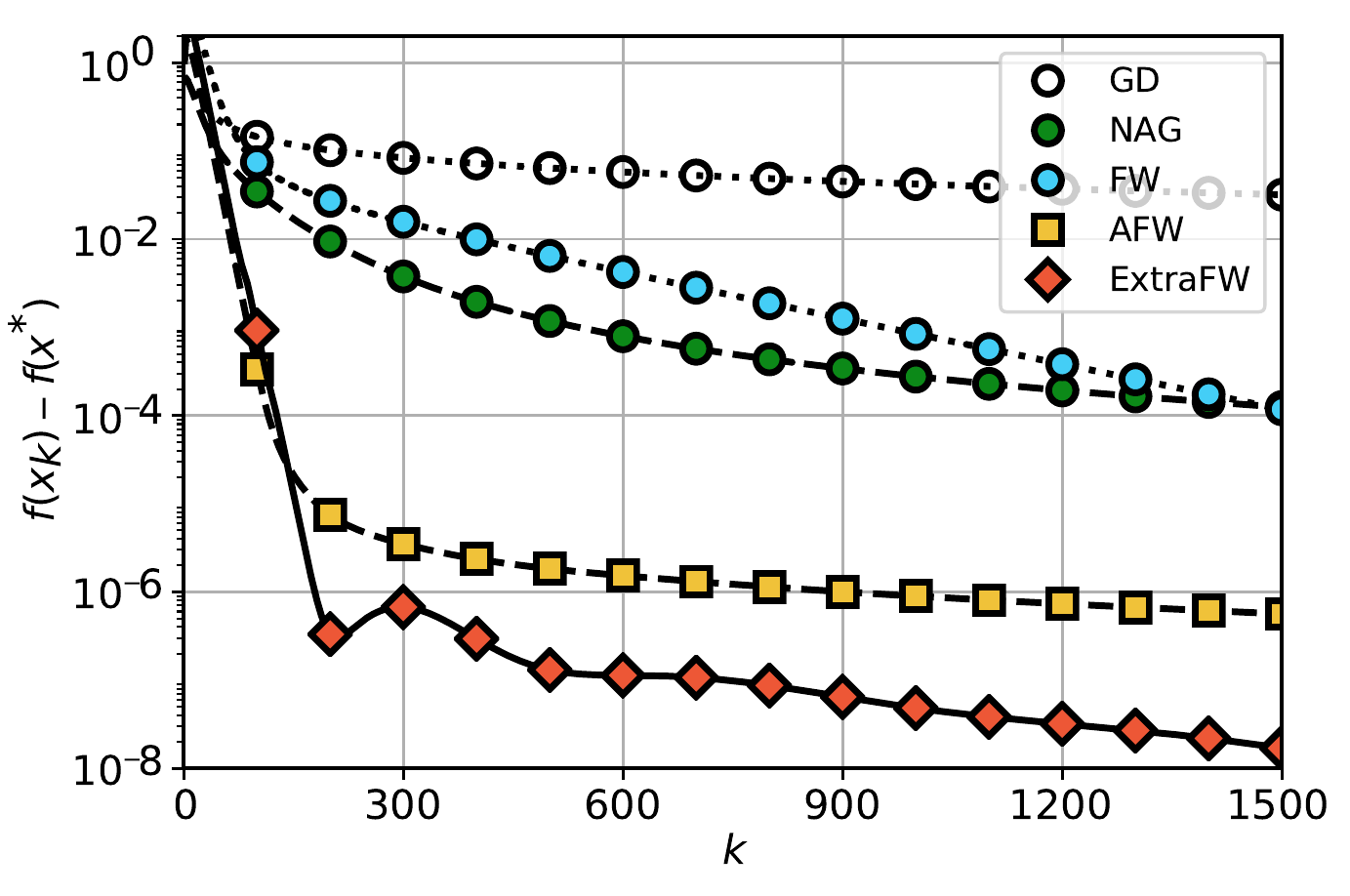}
		\\ (a)  & (b)	&  (c) & (d) 
	\end{tabular}
	\caption{Performance of ExtraFW for binary classification with an $\ell_2$ norm ball constraint on datasets: (a) \textit{mnist}, (b) \textit{w7a}, (c) \textit{realsim}, and, (d) \textit{mushroom}.}
	 \label{fig.l2}
\end{figure*}

\begin{figure*}[t]
	\vspace{-0.1cm}
	\centering
	\begin{tabular}{cccc}
		\hspace{-0.2cm}
		\includegraphics[width=.24\textwidth]{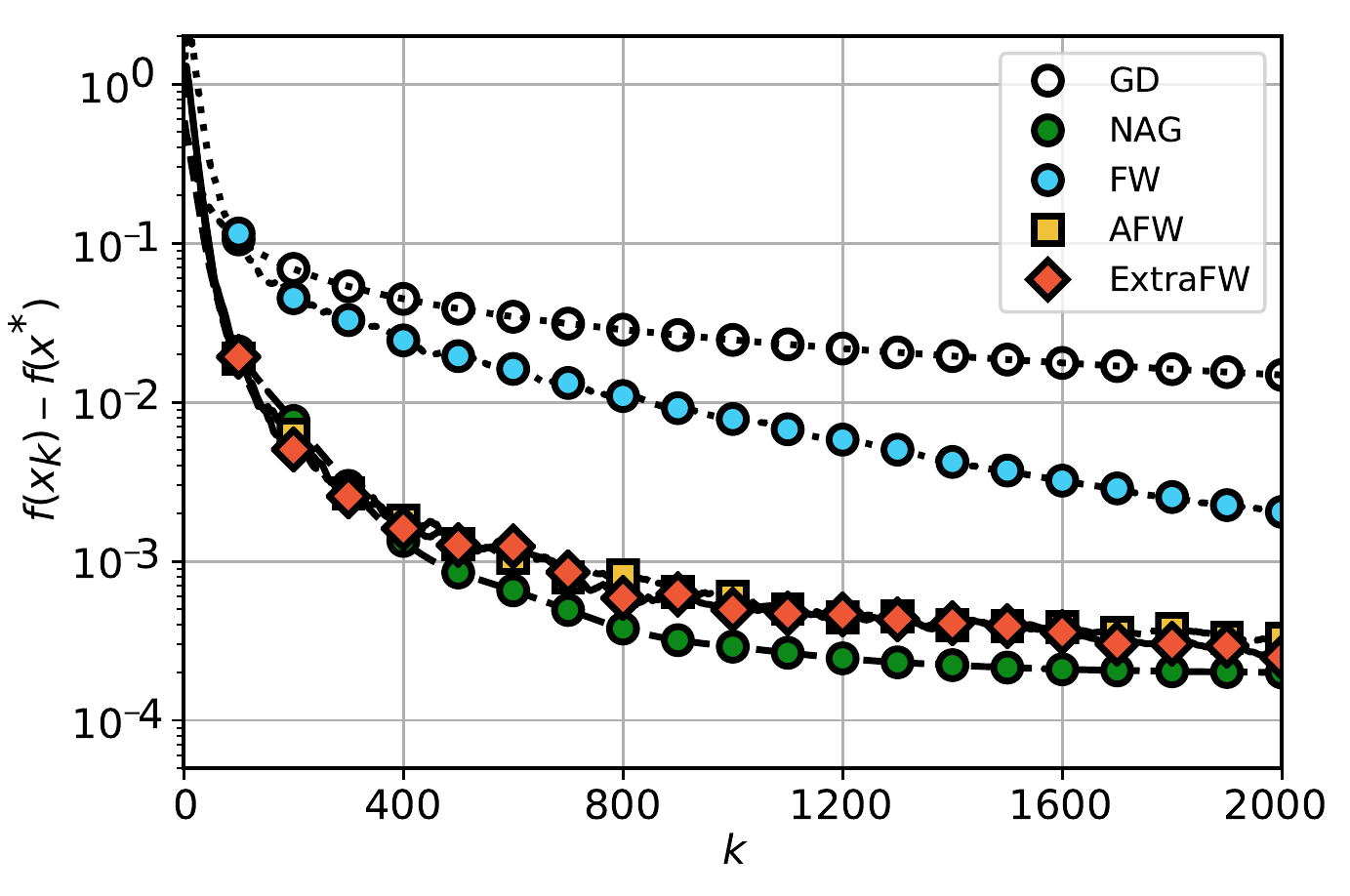}&
		\hspace{-0.4cm}
		\includegraphics[width=.24\textwidth]{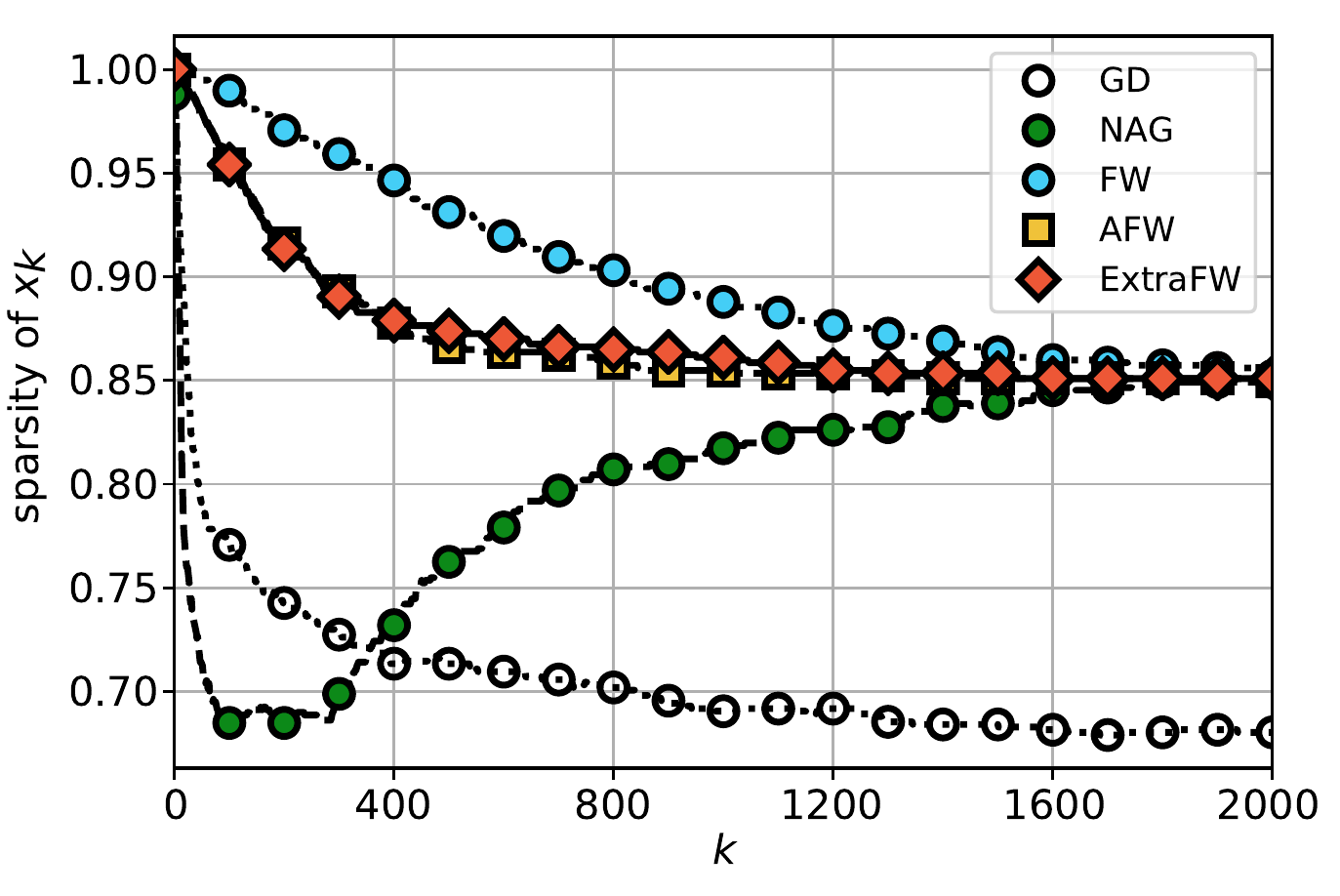}&
		\hspace{-0.4cm}
		\includegraphics[width=.24\textwidth]{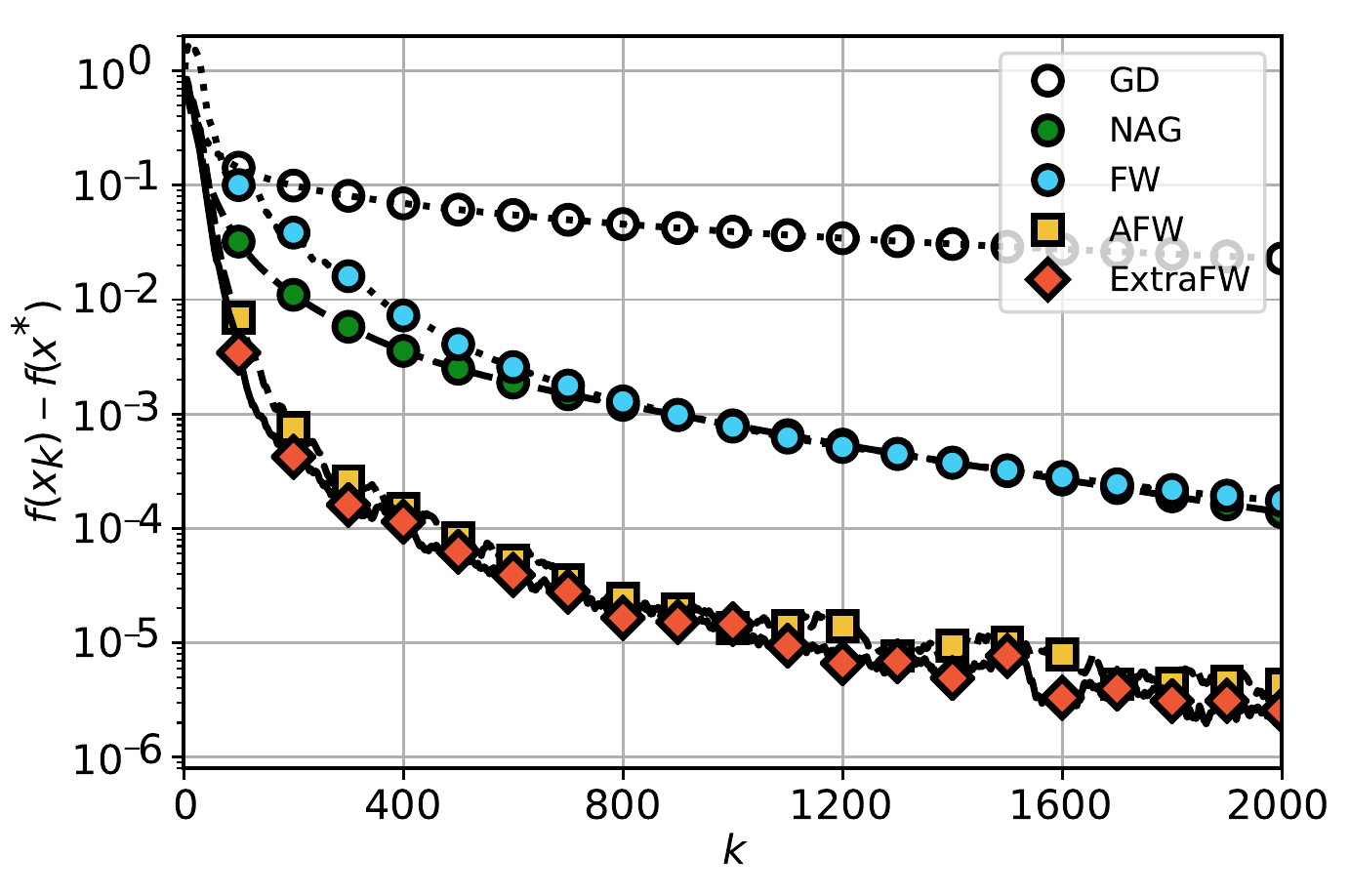}&
		\hspace{-0.4cm}
		\includegraphics[width=.24\textwidth]{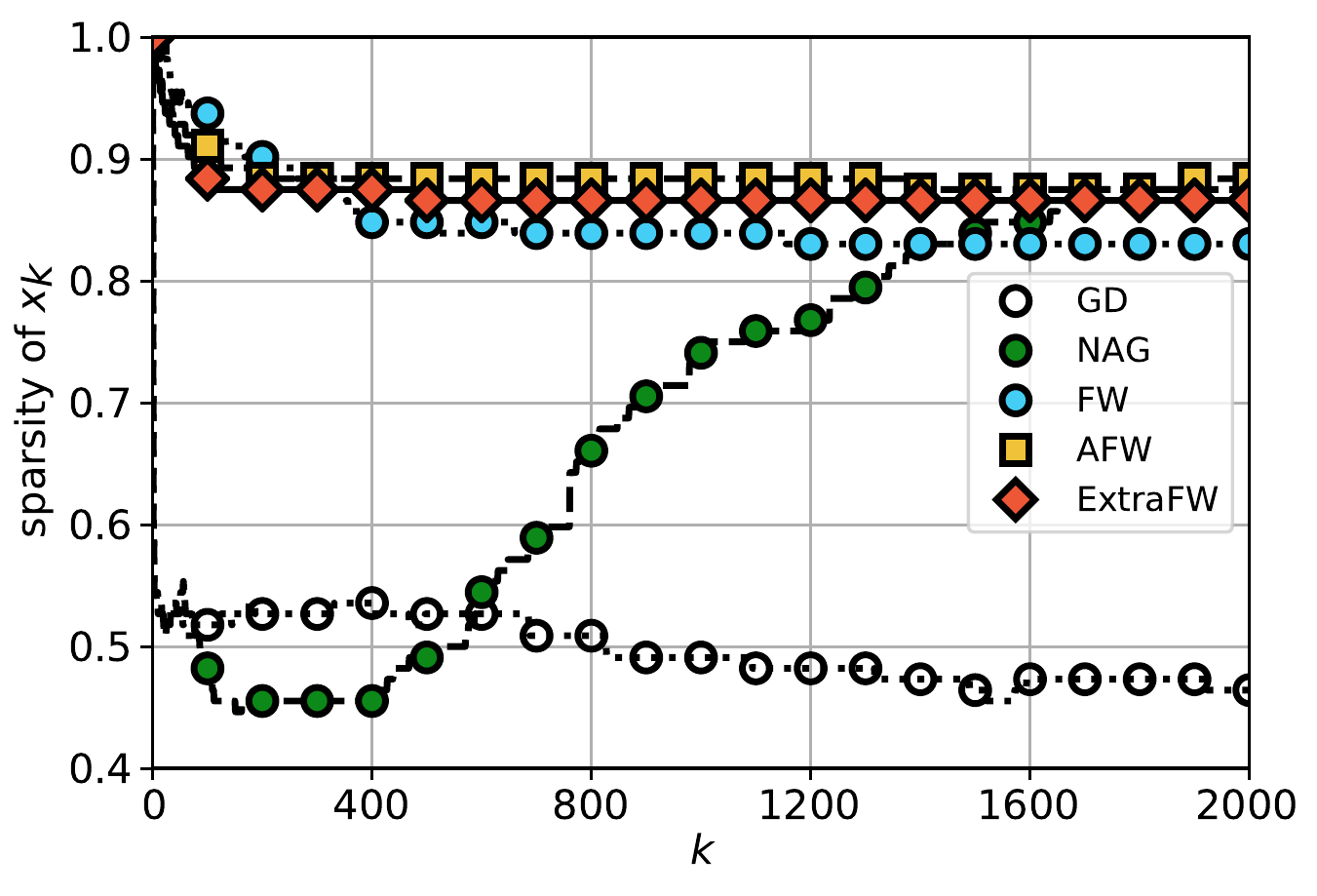}
		\\ (a1)  &  (a2) & (b1)  & (b2) 
	\end{tabular}
	\caption{Performance of ExtraFW for binary classification with an $\ell_1$ norm ball constraint: (a1) optimality error on \textit{mnist}, (a2) solution sparsity on \textit{mnist}, (b1) optimality error on \textit{mushroom}, and, (b2) solution sparsity on \textit{mushroom}.} 
	 \label{fig.l1}
\end{figure*}

\begin{figure*}[t]
	\vspace{-0.1cm}
	\centering
	\begin{tabular}{cccc}
		\hspace{-0.2cm}
		\includegraphics[width=.24\textwidth]{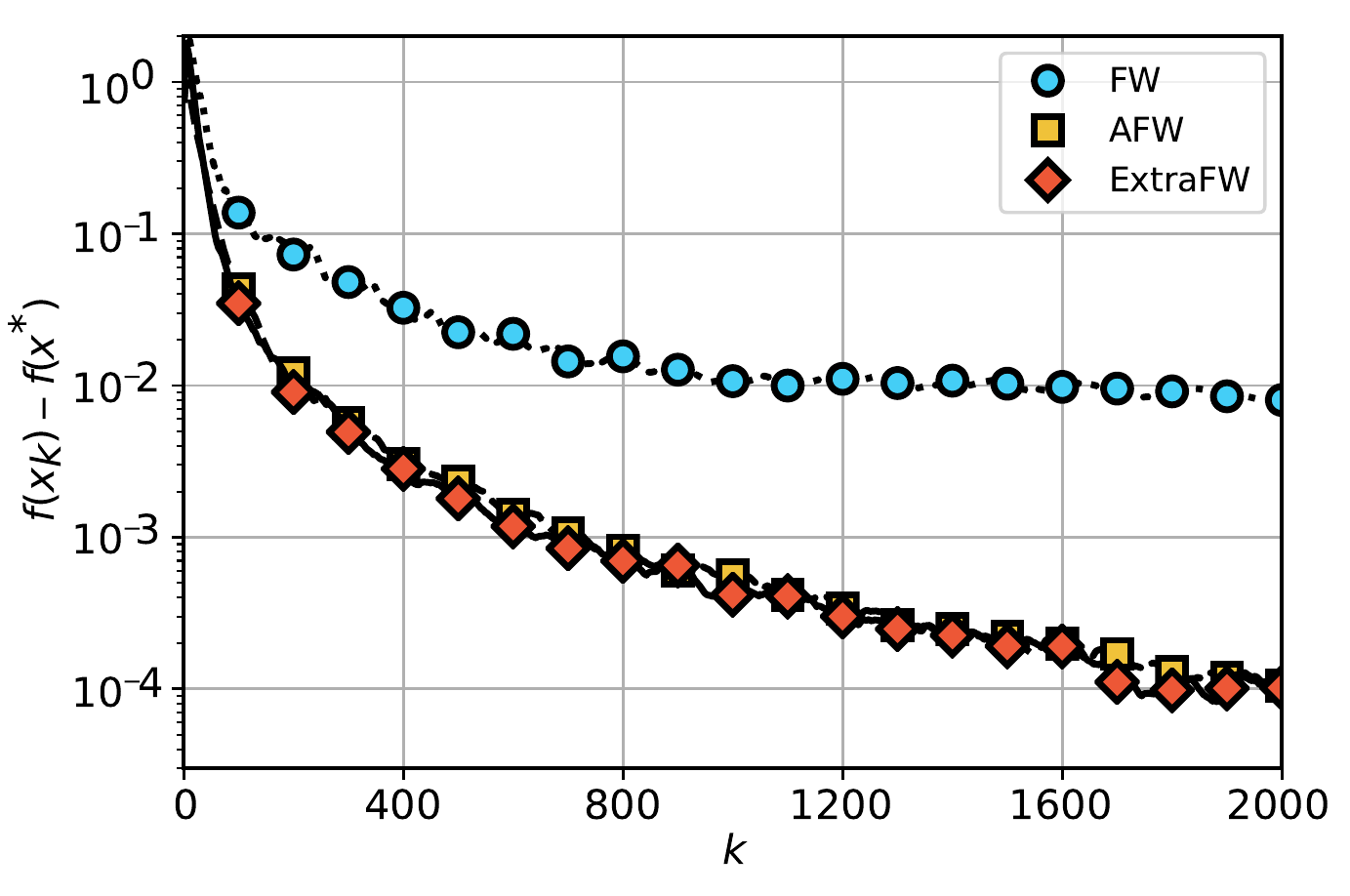}&
		\hspace{-0.4cm}
		\includegraphics[width=.24\textwidth]{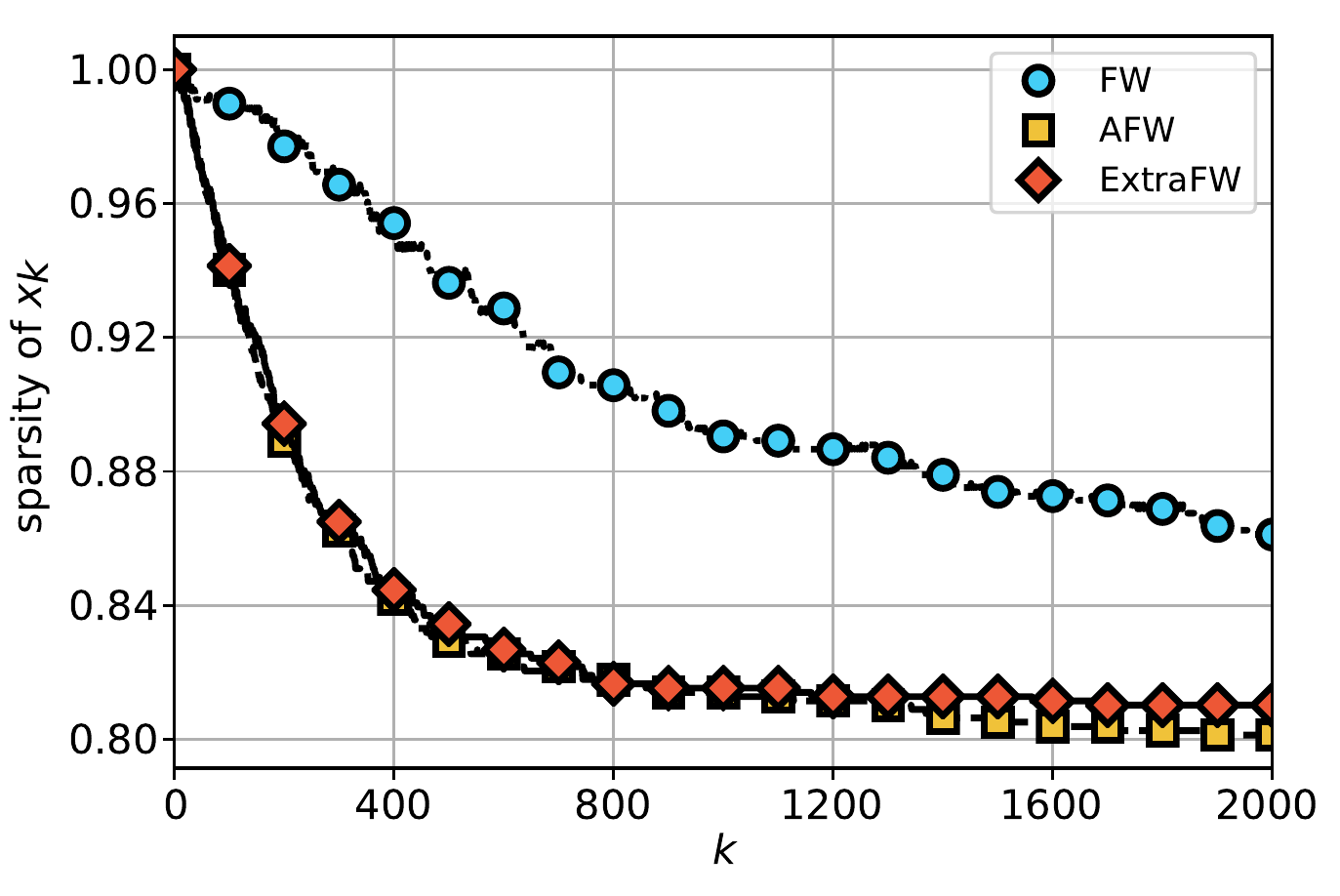}&
		\hspace{-0.4cm}
		\includegraphics[width=.24\textwidth]{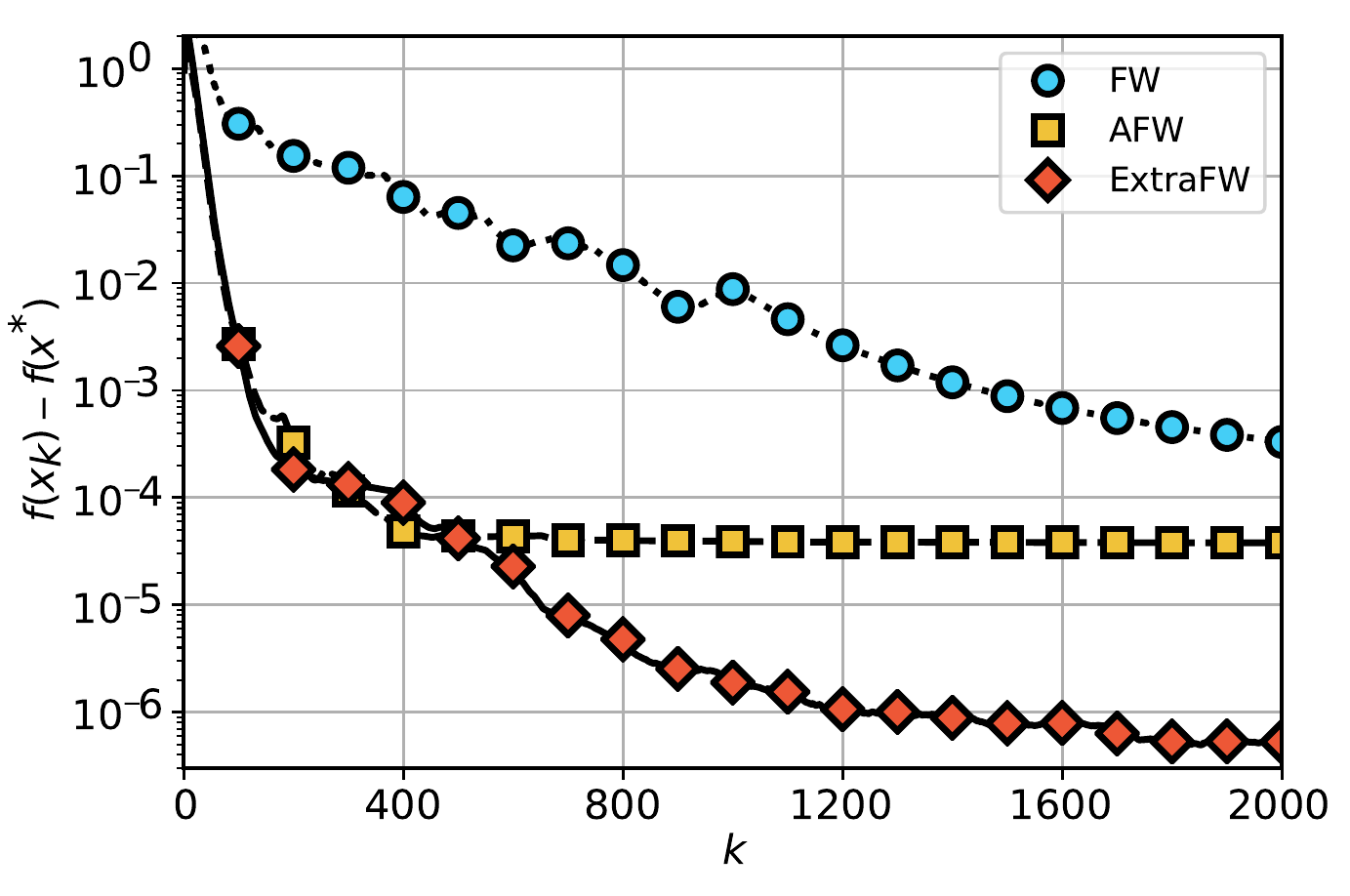}&
		\hspace{-0.4cm}
		\includegraphics[width=.24\textwidth]{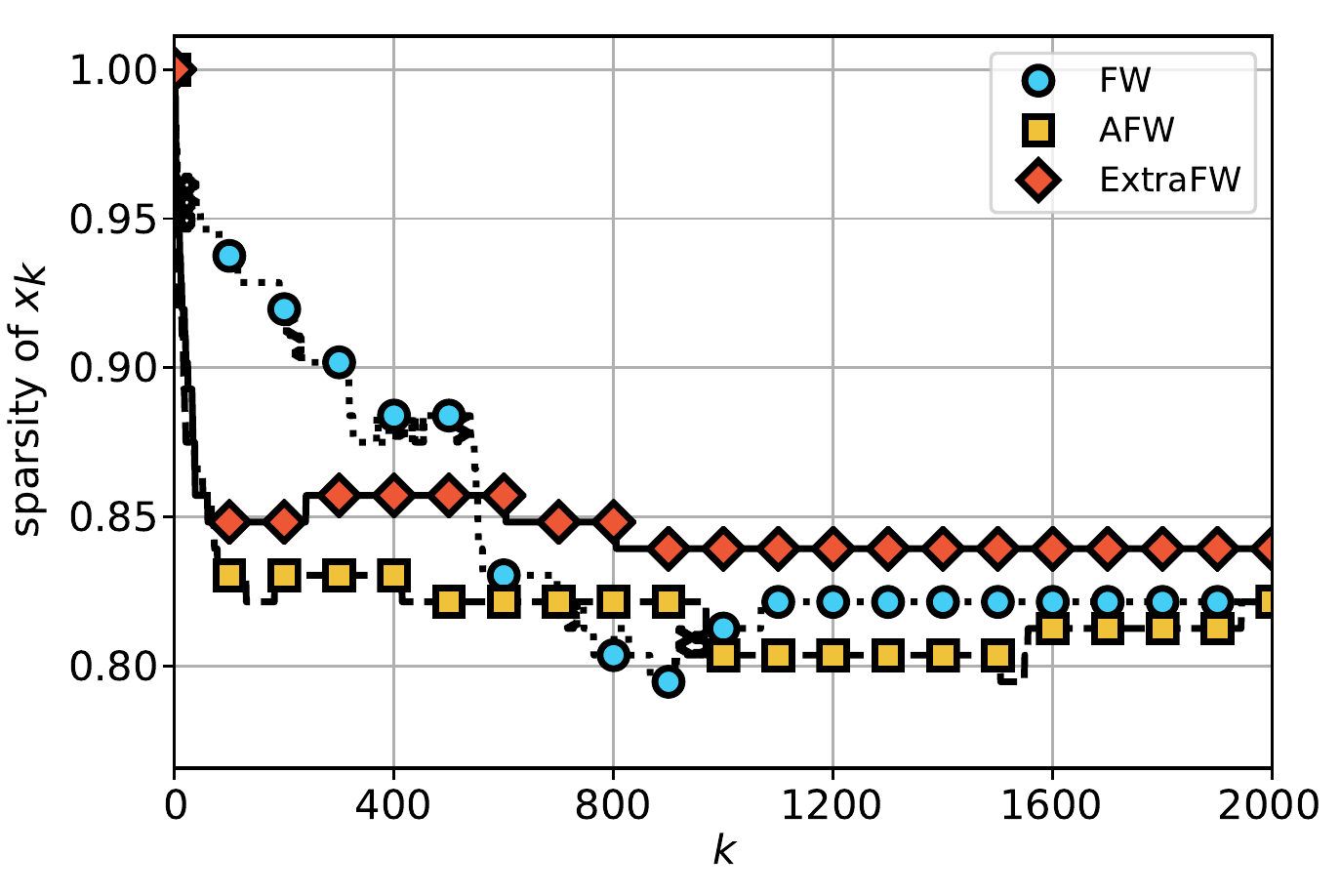}
		\\ (a1)  &  (a2)  & (b1)  & (b2)
	\end{tabular}
	\caption{Performance of ExtraFW for binary classification with an $n$-support norm ball constraint: (a1) optimality error on \textit{mnist}, (a2) solution sparsity on \textit{mnist}, (b1) optimality error on \textit{mushroom}, and, (b2) solution sparsity on \textit{mushroom}. } 
	 \label{fig.n_supp}
\end{figure*}

\textbf{Other constraints.}  Note that the faster rates for ExtraFW are not limited to the exemplified constraint sets. In principle, if i) certain structure such as sparsity is promoted by the constraint set so that $\mathbf{x}^*$ is likely to lie on the boundary of ${\cal X}$; and ii) one can ensure the uniqueness of $\mathbf{v}_k$ through either a closed-form solution or a specific implementation manner, the acceleration of ExtraFW is achievable. Discussions for faster rates on a simplex ${\cal X}$ can be found in Appendix \ref{apdx.ell1}. In addition, one can easily extend our results to the matrix case, where the constraint set is the Frobenius or the nuclear norm ball since they are $\ell_2$ and $\ell_1$ norms on the singular values of matrices, respectively.

\section{Numerical Tests}\label{sec.numerical}
This section deals with numerical tests of ExtraFW to showcase its effectiveness on different machine learning problems. Due to the space limitation, details of the datasets and implementation are deferred to Appendix \ref{apdx.numerical}. For comparison, the benchmarked algorithms are chosen as: i) GD with standard step size $\frac{1}{L}$; ii) Nesterov accelerated gradient (NAG) with step sizes in \citep{allen2014}; iii) FW with parameter-free step size $\frac{2}{k+2}$ \citep{jaggi2013}; and iv) AFW with step size $\frac{2}{k+3}$ \citep{li2020}.

\begin{figure*}[t]
	\vspace{-0.1cm}
	\centering
	\begin{tabular}{cccc}
		\hspace{-0.2cm}
		\includegraphics[width=.24\textwidth]{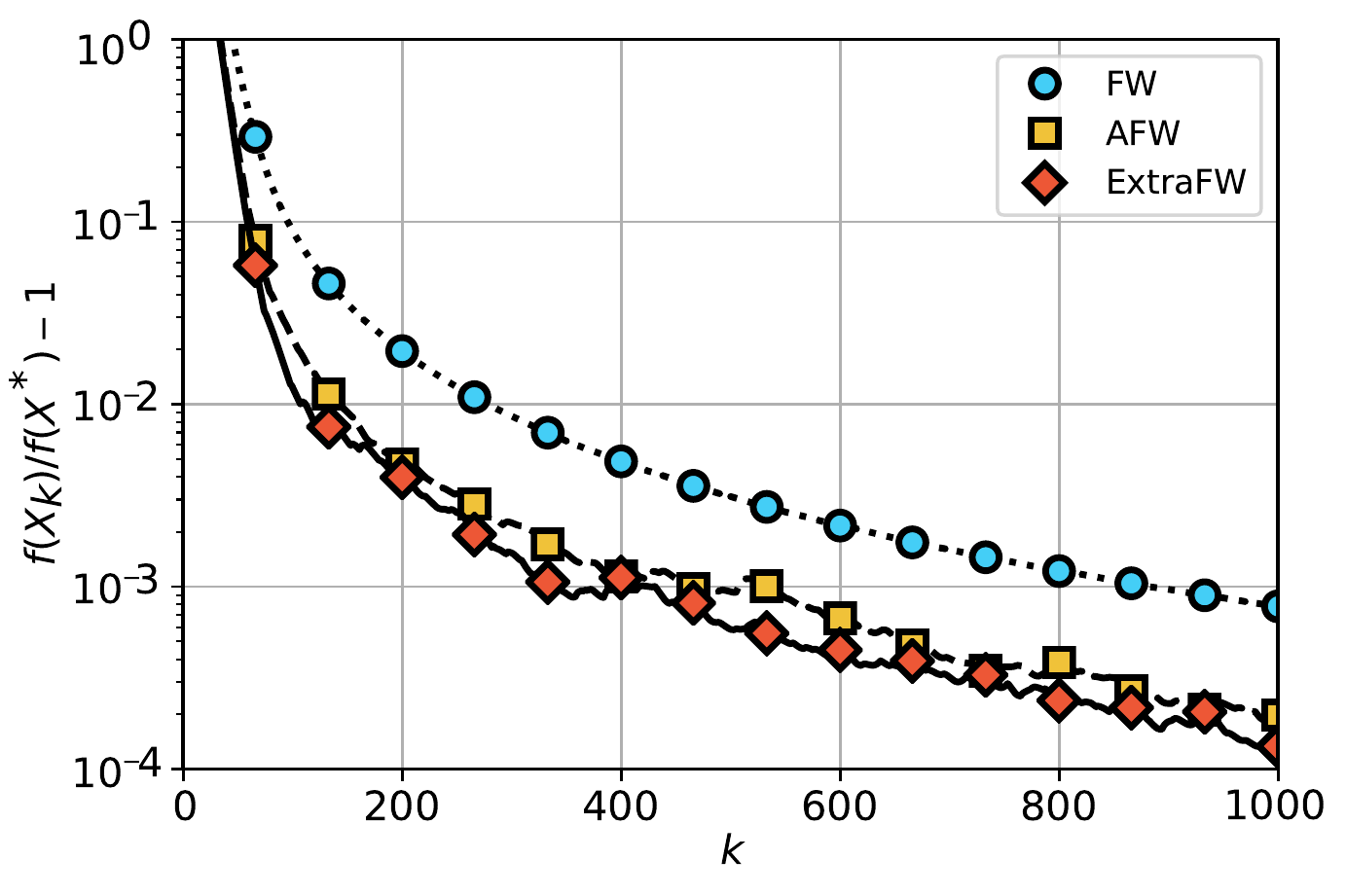}&
		\hspace{-0.4cm}
		\includegraphics[width=.235\textwidth]{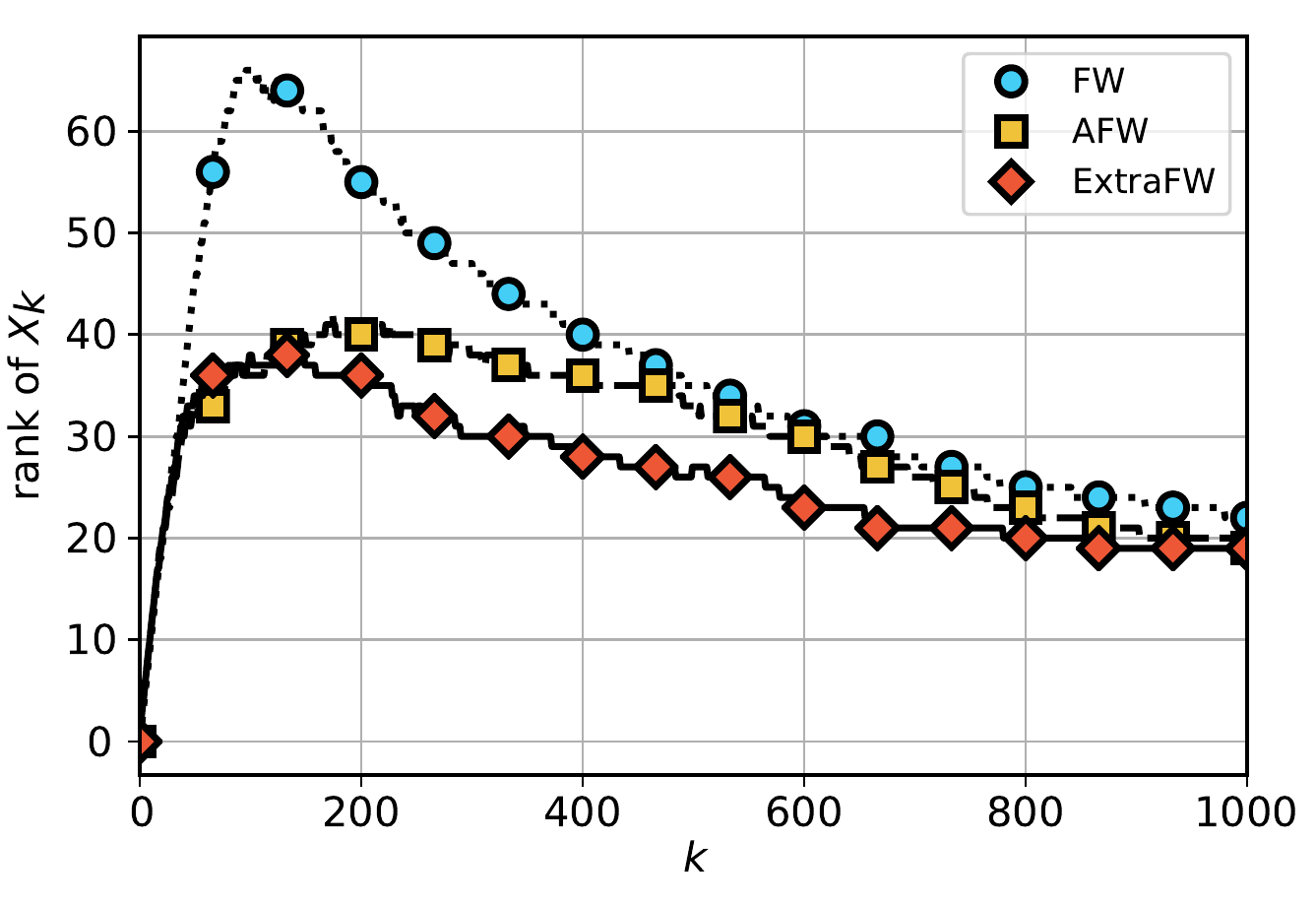}&
		\hspace{-0.4cm}
		\includegraphics[width=.24\textwidth]{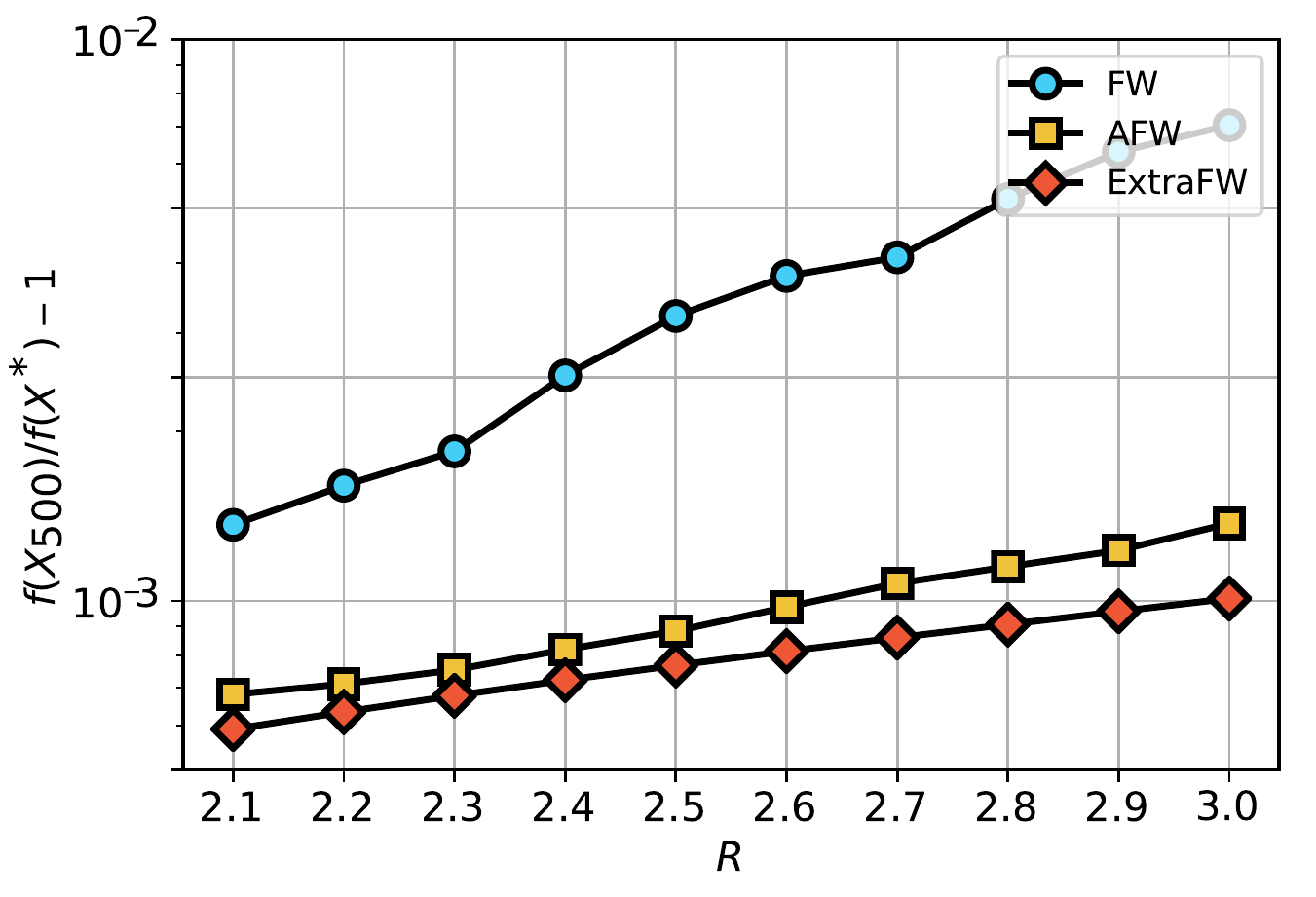} &
		\hspace{-0.4cm}
		\includegraphics[width=.235\textwidth]{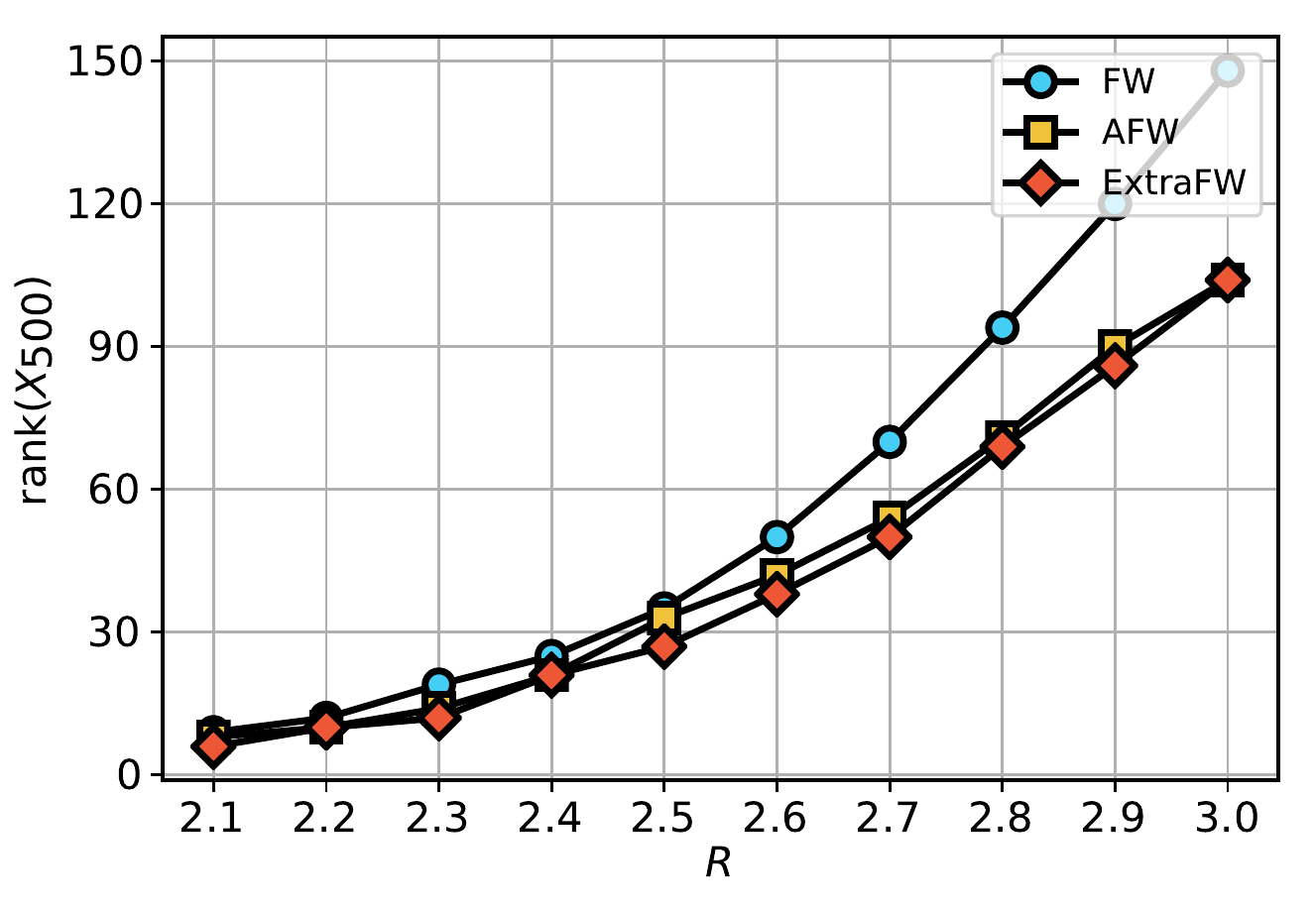}
		\\ (a)  & (b) &  (c) & (d) 
	\end{tabular}
	\caption{Performance of ExtraFW for matrix completion: (a) optimality vs $k$, (b) solution rank vs $k$,  (c) optimality at $k=500$ vs $R$, and,  (d) solution rank at $k=500$ vs $R$. }
	\label{fig.mtrx_comp}
\end{figure*}

\subsection{Binary Classification}
We first investigate the performance of ExtraFW on binary classification using logistic regression. The constraints considered include: i) $\ell_2$ norm ball for generalization merits; and, ii) $\ell_1$ and $n$-support norm ball for promoting a sparse solution. The objective function is
\begin{equation}\label{eq.test}
	f(\mathbf{x}) =\frac{1}{N} \sum_{i =1}^N \ln \big(1+ \exp(- b_i \langle \mathbf{a}_i, \mathbf{x} \rangle ) \big) 
\end{equation}
where $(\mathbf{a}_i, b_i)$ is the (feature, label) pair of datum $i$, and $N$ is the number of data. Datasets \textit{mnist} and those from LIBSVM\footnote{\url{http://yann.lecun.com/exdb/mnist/}, and \url{https://www.csie.ntu.edu.tw/~cjlin/libsvmtools/datasets/binary.html}.} are used in the numerical tests. Figures reporting test accuracy, and additional tests are postponed into Appendix.

\textbf{$\ell_2$ norm ball constraint.} We start with ${\cal X} = \{ \mathbf{x}| \| \mathbf{x} \|_2 \leq R \}$. The optimality error are plotted in Figure \ref{fig.l2}. On all tested datasets, ExtraFW outperforms AFW, NAG, FW and GD, demonstrating the ${\cal O}(\frac{1}{k^2})$ convergence rate established in Theorem \ref{thm.acc}. In addition, the simulation also suggests that $T$ is in general small for logistic loss. On dataset \textit{w7a} and \textit{mushroom}, ExtraFW is significantly faster than AFW. All these observations jointly confirm the usefulness of the extra gradient and the PC update.

\textbf{$\ell_1$ norm ball constraint.} Let ${\cal X} = \{ \mathbf{x}| \|\mathbf{x}\|_1 \leq R \}$ be the constraint set to promote sparsity on the solution. Note that FW type updates directly guarantee that $\mathbf{x}_k$ has at most $k$ non-zero entries when initialized at $\mathbf{x}_0 = \mathbf{0}$; see detailed discussions in Appendix \ref{apdx.classification}. In the simulation, $R$ is tuned to obtain a solution that is almost as sparse as the dataset itself. The numerical results on datasets \textit{mnist} and \textit{mushroom} including both optimality error and the sparsity level of the solution can be found in Figure \ref{fig.l1}. On dataset \textit{mnist}, ExtraFW slightly outperforms AFW but is not as fast as NAG. However, ExtraFW consistently finds solutions sparser than NAG. While on dataset \textit{mushroom}, it can be seen that both AFW and ExtraFW outperform NAG, with ExtraFW slightly faster than AFW. And ExtraFW finds sparser solutions than NAG.

\textbf{$n$-support norm ball constraint.} Effective projection onto such a constraint is unknown yet and hence GD and NAG are not included in the test. The performance of ExtraFW can be found in Figure \ref{fig.n_supp}. On dataset \textit{mnist}, both AFW and ExtraFW converge much faster than FW with ExtraFW slightly faster than AFW. However, FW trades the solution accuracy with its sparsity. On dataset \textit{mushroom}, ExtraFW converges much faster than AFW and FW, while finding the sparsest solution. 


\subsection{Matrix Completion}

We then consider matrix completion problems that are ubiquitous in recommender systems. Consider a matrix $\mathbf{A} \in \mathbb{R}^{m \times n}$ with partially observed entries, that is, entries $A_{ij}$ for $(i,j) \in {\cal K}$ are known, where ${\cal K} \subset \{1,\ldots,m \} \times \{1,\ldots,n \}$. Note that the observed entries can also be contaminated by noise. The task is to predict the unobserved entries of $\mathbf{A}$. Although this problem can be approached in several ways, within the scope of recommender systems, a commonly adopted empirical observation is that $\mathbf{A}$ is low rank \citep{bennett2007,bell2007,fazel2002}. Hence the objective boils down to
\begin{align}\label{eq.mtrx_complete_relax}
	\min_{\mathbf{X}} ~~ \frac{1}{2} \sum_{(i,j) \in {\cal K}} (X_{ij} - A_{ij})^2 ~~~	\text{s.t.} ~~ \| \mathbf{X} \|_{\rm nuc}\leq R
\end{align}
where $\| \cdot \|_{\rm nuc}$ denotes the nuclear norm. Problem \eqref{eq.mtrx_complete_relax} is difficult to be solved via GD or NAG because projection onto a nuclear norm ball requires to perform SVD, which has complexity ${\cal O}\big(mn (m \wedge n)\big)$. On the contrary, FW and its variants are more suitable for \eqref{eq.mtrx_complete_relax} given the facts: i) Assumptions \ref{as.1} -- \ref{as.3} are satisfied under nuclear norm \citep{freund2017}; ii) FW step can be solved easily with complexity at the same order as the number of nonzero entries; and iii) the update promotes low-rank solution directly \citep{freund2017}. More on ii) and iii) are discussed in Appendix \ref{apdx.mtrx_completion}.

We test ExtraFW on a widely used dataset, \textit{MovieLens100K}\footnote{\url{https://grouplens.org/datasets/movielens/100k/}}. The experiments follow the same steps in \citep{freund2017}. The numerical performance of ExtraFW, AFW, and FW can be found in Figure \ref{fig.mtrx_comp}. We plot the optimality error and rank versus $k$ choosing $R = 2.5$ in Figures~\ref{fig.mtrx_comp}(a) and~\ref{fig.mtrx_comp}(b). It is observed that ExtraFW exhibits the best performance in terms of both optimality error and solution rank. In particular, ExtraFW roughly achieves $2.5$x performance improvement compared with FW in terms of optimality error. We further compare the convergence of ExtraFW to AFW and FW at iteration $k=500$ under different choices of $R$ in Figures~\ref{fig.mtrx_comp}(c) and~\ref{fig.mtrx_comp}(d). ExtraFW still finds solutions with the lowest optimality error and rank. Moreover, the performance gap between ExtraFW and AFW increases with $R$, suggesting the inclined tendency of preferring ExtraFW over AFW and FW as $R$ grows.

\section{Conclusions}
A new parameter-free FW variant, ExtraFW, is introduced and analyzed
in this work. ExtraFW leverages two gradient evaluations per iteration to update in a ``prediction-correction'' manner. We show
that ExtraFW converges at ${\cal O}(\frac{1}{k})$ on general problems, while achieving a faster rate ${\cal O}(\frac{TLD^2}{k^2})$ on certain types of constraint sets including active $\ell_1$, $\ell_2$ and $n$-support norm balls.
Given the possibility of acceleration, ExtraFW is thus a competitive alternative to FW. The efficiency of ExtraFW is validated on tasks such as i) binary classification with different constraints, where ExtraFW can be even faster than NAG; and ii) matrix completion where ExtraFW finds solutions with lower optimality error and rank rapidly.

\subsection*{Acknowledgements}
BL and GG gratefully acknowledge the support from NSF grants 1711471, and 1901134. LW and ZZ are supported by Alfred P. Sloan Foundation.

\bibliographystyle{plainnat}
\bibliography{datactr_arxiv.bib}


\onecolumn
\appendix

\begin{center}
{\LARGE \bf Supplementary Document for \\
``Enhancing Parameter-Free Frank Wolfe \\with an Extra Subproblem'' }
\end{center}

\section{Missing Proofs in Section \ref{sec.conv_gen} }
\subsection{Proof of Lemma \ref{lemma.es_convergence}}
\begin{proof} 
	If $f(\mathbf{x}_k) \leq \min_{\mathbf{x} \in {\cal X}} \Phi_k(\mathbf{x}) + \xi_k$ holds, then we have
	\begin{align*}
		f(\mathbf{x}_k) \leq \min_{\mathbf{x} \in {\cal X} } \Phi_k(\mathbf{x}) + \xi_k \leq \Phi_k(\mathbf{x}^*)  + \xi_k \leq (1 - \lambda_k)	 f(\mathbf{x}^*) + \lambda_k \Phi_0 (\mathbf{x}^*) + \xi_k
	\end{align*}
	where the last inequality is because Definition \ref{def.es}. Subtracting $f(\mathbf{x}^*) $ on both sides, we arrive at
	\begin{align*}
		f(\mathbf{x}_k) - f(\mathbf{x}^*) \leq  \lambda_k \big( \Phi_0 (\mathbf{x}^*) - f(\mathbf{x}^*) \big) + \xi_k 
	\end{align*}
	which completes the proof.
\end{proof}

\subsection{Proof of Lemma \ref{lemma.es}}
\begin{proof}
We prove $\big( \{\Phi_k(\mathbf{x}) \}_{k=0}^\infty,$$ \{ \lambda_k \}_{k=0}^\infty \big)$ is an ES of $f$ by induction. Because $\lambda_0 =1$, it holds that $\Phi_0(\mathbf{x}) = (1 - \lambda_0) f(\mathbf{x}) + \lambda_0\Phi_0(\mathbf{x}) = \Phi_0(\mathbf{x})$. Suppose that $\Phi_k (\mathbf{x}) \leq (1 - \lambda_k)	 f(\mathbf{x}) + \lambda_k \Phi_0 (\mathbf{x})$ is true for some $k$. We have
	\begin{align*}
		\Phi_{k+1}(\mathbf{x}) & = (1- \delta_k) \Phi_k (\mathbf{x}) + \delta_k \Big[ f(\mathbf{x}_{k+1}) + \big\langle  \nabla f(\mathbf{x}_{k+1}), \mathbf{x} - \mathbf{x}_{k+1}  \big\rangle \Big] \\
		& \stackrel{(\text{a})}{\leq} (1- \delta_k) \Phi_k (\mathbf{x}) + \delta_k f(\mathbf{x}) \\
		& \leq (1- \delta_k) \Big[ (1 - \lambda_k)	 f(\mathbf{x}) + \lambda_k \Phi_0 (\mathbf{x}) \Big] + \delta_k f(\mathbf{x}) \\
		& = (1 - \lambda_{k+1})	 f(\mathbf{x}) + \lambda_{k+1} \Phi_0 (\mathbf{x})
	\end{align*}
	where (a) is because $f$ is convex; and the last equation is by definition of $\lambda_{k+1}$. Together with the fact that $\lim_{k\rightarrow\infty} \lambda_k = 0$, one can see that the tuple $\big( \{\Phi_k(\mathbf{x}) \}_{k=0}^\infty, \{ \lambda_k \}_{k=0}^\infty \big)$ is an ES of $f$. 
	
	Next we show $\big( \{\hat{\Phi}_k(\mathbf{x}) \}_{k=0}^\infty, \{ \lambda_k \}_{k=0}^\infty \big)$ is also an ES. Clearly $\hat{\Phi}_0(\mathbf{x}) = (1 - \lambda_0) f(\mathbf{x}) + \lambda_0\Phi_0(\mathbf{x}) = \hat{\Phi}_0(\mathbf{x})$. Next for $k \geq 0 $, using similar arguments,	we have
	\begin{align*}
		\hat{\Phi}_{k+1}(\mathbf{x}) & = (1- \delta_k) \Phi_k (\mathbf{x}) +  \delta_k \Big[ f(\mathbf{y}_k) +  \big\langle \nabla f(\mathbf{y}_k) , \mathbf{x}  - \mathbf{y}_k \big\rangle  \Big] \\
		& \leq (1- \delta_k) \Phi_k (\mathbf{x}) + \delta_k f(\mathbf{x}) \\
		& \leq (1- \delta_k) \Big[ (1 - \lambda_k)	 f(\mathbf{x}) + \lambda_k \Phi_0 (\mathbf{x}) \Big] + \delta_k f(\mathbf{x}) \\
		& = (1 - \lambda_{k+1})	 f(\mathbf{x}) + \lambda_{k+1} \Phi_0 (\mathbf{x}) \\
		& = (1 - \lambda_{k+1})	 f(\mathbf{x}) + \lambda_{k+1} \hat{\Phi}_0 (\mathbf{x}).
	\end{align*}
	The proof is thus completed.
\end{proof}

\subsection{Proof of Lemma \ref{lemma.vstar}}
\begin{proof}
	For convenience, denote $B_k (\mathbf{x}):= f(\mathbf{x}_k) + \langle \nabla f(\mathbf{x}_k), \mathbf{x} - \mathbf{x}_k \rangle$. We can unroll $\Phi_{k+1}(\mathbf{x})$ as
	\begin{align}\label{eq.stop_use1}
		\Phi_{k+1}(\mathbf{x}) &= (1 - \delta_k) \Phi_k(\mathbf{x}) + \delta_k B_{k+1}(\mathbf{x}) 	 \\
		& = (1 - \delta_k)(1 - \delta_{k-1}) \Phi_{k-1}(\mathbf{x}) + (1 - \delta_k)\delta_{k-1} B_k(\mathbf{x})  + \delta_k B_{k+1}(\mathbf{x}) \nonumber \\
		& = \Phi_0(\mathbf{x}) \prod_{\tau = 0}^k (1-\delta_\tau) + \sum_{\tau=0}^k \delta_\tau B_{\tau+1}(\mathbf{x}) \prod_{j = \tau+1}^k (1-\delta_j) \nonumber \\
		& = f(\mathbf{x}_0) \prod_{\tau = 0}^k (1-\delta_\tau) + \sum_{\tau=0}^k \delta_\tau B_{\tau+1}(\mathbf{x}) \prod_{j = \tau+1}^k (1-\delta_j) \nonumber.
	\end{align}
	Hence, the minimizer of $\Phi_{k+1}(\mathbf{x})$ can be rewritten as
	\begin{align}\label{eq.ccc1}
		\argmin_{\mathbf{x} \in {\cal X}} ~	\Phi_{k+1}(\mathbf{x}) & = \argmin_{\mathbf{x} \in {\cal X}} ~ f(\mathbf{x}_0) \prod_{\tau = 0}^k (1-\delta_\tau)  + \sum_{\tau=0}^k \delta_\tau B_{\tau+1}(\mathbf{x}) \prod_{j = \tau+1}^k (1-\delta_j)  \\
		&= \argmin_{\mathbf{x} \in {\cal X}} ~  \sum_{\tau=0}^k \delta_\tau \big[ f(\mathbf{x}_{\tau+1}) + \langle \nabla f(\mathbf{x}_{\tau+1}), \mathbf{x} - \mathbf{x}_{\tau+1} \rangle \big] \prod_{j = \tau+1}^k (1-\delta_j) \nonumber  \\
		& = \argmin_{\mathbf{x} \in {\cal X}} ~  \sum_{\tau=0}^k \delta_\tau  \langle \nabla f(\mathbf{x}_{\tau+1}), \mathbf{x}  \rangle  \prod_{j = \tau+1}^k (1-\delta_j) \nonumber \\
		& = \argmin_{\mathbf{x} \in {\cal X}} ~ \sum_{\tau=0}^k   \Big\langle \delta_\tau  \nabla f(\mathbf{x}_{\tau+1}) \prod_{j = \tau+1}^k (1-\delta_j), \mathbf{x}  \Big\rangle \nonumber \\
		& =   \argmin_{\mathbf{x} \in {\cal X}} ~  \langle \mathbf{g}_{k+1} , \mathbf{x} \rangle  \nonumber
	\end{align}
	where the last equation is because
	\begin{align*}
		\mathbf{g}_{k+1} & =  (1-\delta_k) \mathbf{g}_k +  \delta_k \nabla f(\mathbf{x}_{k+1})  \\
		& = (1-\delta_k)  (1-\delta_{k-1}) \mathbf{g}_{k-1} +  (1-\delta_k) \delta_{k-1} \nabla f(\mathbf{x}_k) +  \delta_k \nabla f(\mathbf{x}_{k+1}) \nonumber  \\
		& = \mathbf{g}_0 \prod_{\tau = 0}^k (1-\delta_\tau) + \sum_{\tau = 0}^k \delta_\tau \nabla f(\mathbf{x}_{\tau+1}) \prod_{j=\tau+1}^k (1-\delta_j) = \sum_{\tau = 0}^k \delta_\tau \nabla f(\mathbf{x}_{\tau+1}) \prod_{j=\tau+1}^k (1-\delta_j). \nonumber
	\end{align*}
	From \eqref{eq.ccc1} it is not hard to see $\mathbf{v}_{k+1}$ minimizes $\Phi_{k+1}(\mathbf{x})$.

	If we write $\hat{\mathbf{g}}_{k+1}$ explicitly, we can obtain 
	\begin{align*}
		\hat{\Phi}_{k+1}(\mathbf{x}) &= (1 - \delta_k) \Phi_k(\mathbf{x}) + \delta_k \Big[ f(\mathbf{y}_k) +  \big\langle \nabla f(\mathbf{y}_k) , \mathbf{x}  - \mathbf{y}_k \big\rangle  \Big] 	 \\
		& = f(\mathbf{x}_0) \prod_{\tau = 0}^k (1-\delta_\tau)+ \sum_{\tau=0}^{k-1} \delta_\tau B_{\tau+1}(\mathbf{x}) \prod_{j = \tau+1}^k (1-\delta_j)  + \delta_k \Big[ f(\mathbf{y}_k) +  \big\langle \nabla f(\mathbf{y}_k) , \mathbf{x}  - \mathbf{y}_k \big\rangle  \Big] .
	\end{align*}
 	Hence using similar arguments as above we have
 	\begin{align*}
 		\argmin_{\mathbf{x} \in {\cal X}} ~	\hat{\Phi}_{k+1}(\mathbf{x}) & = \argmin_{\mathbf{x} \in {\cal X}} ~	  \Big\langle \delta_k \nabla f(\mathbf{y}_k) + \sum_{\tau=0}^{k-1} \delta_\tau  \nabla f(\mathbf{x}_{\tau+1})  \prod_{j = \tau+1}^k (1-\delta_j) , \mathbf{x}  \Big\rangle  \\
 		&= \argmin_{\mathbf{x} \in {\cal X}} ~   \langle \hat{\mathbf{g}}_{k+1} , \mathbf{x} \rangle  = \hat{\mathbf{v}}_{k+1}
  	\end{align*}
  	which implies that $\hat{\mathbf{v}}_{k+1}$ is a minimizer of $\hat{\Phi}_{k+1}(\mathbf{x})$ over ${\cal X}$. The lemma is thus proved.
\end{proof}

\subsection{Proof of Lemma \ref{lemma_phistar}}

\begin{proof}
	We prove this lemma by induction. Since $\Phi_0(\mathbf{x}) \equiv f(\mathbf{x}_0)$ and $\xi_0 = 0$, it is clear that $f(\mathbf{x}_0) \leq \Phi_0^* + \xi_0 $. 
	
	Now suppose that $f(\mathbf{x}_k )\leq \Phi_k^*  + \xi_k$ holds for some $k >  0$, we will show $f(\mathbf{x}_{k+1} )\leq \Phi_{k+1}^* + \xi_{k+1}$. To start with, we have from Assumption \ref{as.1} that
	\begin{align}\label{eq.bbb3}
		f(\mathbf{x}_{k+1}) & \leq f(\mathbf{y}_k)	 + \big \langle \nabla f(\mathbf{y}_k), \mathbf{x}_{k+1} - \mathbf{y}_k \big \rangle + \frac{L}{2} \| \mathbf{x}_{k+1} - \mathbf{y}_k \|^2  \\
		& \stackrel{(\text{a})}{=}  f(\mathbf{y}_k)	 +  (1-\delta_k) \big \langle \nabla f(\mathbf{y}_k), \mathbf{x}_k - \mathbf{y}_k \big \rangle + \delta_k \big \langle \nabla f(\mathbf{y}_k), \hat{\mathbf{v}}_{k+1} - \mathbf{y}_k \big \rangle + \frac{L}{2} \| \mathbf{x}_{k+1} - \mathbf{y}_k \|^2 \nonumber \\
		&  \stackrel{(\text{b})}{=} f(\mathbf{y}_k)	 +  (1-\delta_k) \big \langle \nabla f(\mathbf{y}_k), \mathbf{x}_k - \mathbf{y}_k \big \rangle + \delta_k \big \langle \nabla f(\mathbf{y}_k), \hat{\mathbf{v}}_{k+1} - \mathbf{y}_k \big \rangle + \frac{L \delta_k^2}{2} \| \hat{\mathbf{v}}_{k+1} - \mathbf{v}_k \|^2  \nonumber \\
		& \stackrel{(\text{c})}{\leq} (1-\delta_k) f(\mathbf{x}_k) + \delta_k f(\mathbf{y}_k) + \delta_k \big \langle \nabla f(\mathbf{y}_k), \hat{\mathbf{v}}_{k+1} - \mathbf{y}_k \big \rangle + \frac{L \delta_k^2}{2} \| \hat{\mathbf{v}}_{k+1} - \mathbf{v}_k \|^2 \nonumber
	\end{align}
	where (a) is because $\mathbf{x}_{k+1} = (1-\delta_k) \mathbf{x}_k + \delta_k \hat{\mathbf{v}}_{k+1}$; (b) is by the choice of $\mathbf{x}_{k+1}$ and $\mathbf{y}_k$; and (c) is from convexity, that is, $ \langle \nabla f(\mathbf{y}_k), \mathbf{x}_k - \mathbf{y}_k \rangle \leq f(\mathbf{x}_k) - f(\mathbf{y}_k)$. For convenience we denote $\hat{\Phi}_k^* := \hat{\Phi}_k(\hat{\mathbf{v}}_k)$ as the minimum value of $\hat{\Phi}_k(\mathbf{x})$ over ${\cal X}$ (the equation here is the result of Lemma \ref{lemma.vstar}). Then we have 
	\begin{align*}
		\hat{\Phi}_{k+1}^*	& = \hat{\Phi}_{k+1}(\hat{\mathbf{v}}_{k+1}) \stackrel{(\text{d})}{=}  (1 - \delta_k) \Phi_k(\hat{\mathbf{v}}_{k+1}) + \delta_k \Big[ f(\mathbf{y}_k) +  \big\langle \nabla f(\mathbf{y}_k) , \hat{\mathbf{v}}_{k+1} - \mathbf{y}_k \big\rangle \Big]  \\
		& \stackrel{(\text{e})}{\geq} (1 - \delta_k) \Phi_k^* + \delta_k \Big[ f(\mathbf{y}_k) +  \big\langle \nabla f(\mathbf{y}_k) , \hat{\mathbf{v}}_{k+1} - \mathbf{y}_k \big\rangle \Big]  \nonumber \\
		&\stackrel{(\text{f})}{\geq} (1 - \delta_k) f(\mathbf{x}_k) + \delta_k \Big[ f(\mathbf{y}_k) +  \big\langle \nabla f(\mathbf{y}_k) , \hat{\mathbf{v}}_{k+1} - \mathbf{y}_k \big\rangle \Big] - (1 - \delta_k) \xi_k \nonumber \\
		& \stackrel{(\text{g})}{\geq} f(\mathbf{x}_{k+1}) - \frac{L \delta_k^2}{2}  \| \hat{\mathbf{v}}_{k+1} - \mathbf{v}_k \|^2 - (1 - \delta_k) \xi_k \nonumber \\
		& \geq f(\mathbf{x}_{k+1}) - \frac{L D^2 \delta_k^2}{2} - (1 - \delta_k) \xi_k \nonumber
	\end{align*}
	where (d) is by the definition of $ \hat{\Phi}_{k+1}(\mathbf{x}) $; (e) uses $\Phi_k(\hat{\mathbf{v}}_{k+1}) \geq \Phi_k^*$; (f) is by the induction hypothesis $f(\mathbf{x}_k )\leq \Phi_k^*  + \xi_k$; (g) is by plugging \eqref{eq.bbb3} in; and the last inequality is because of Assumption \ref{as.3}. Rearrange the terms, we have
	\begin{align}\label{eq.bbb1}
		f(\mathbf{x}_{k+1})	& \leq \hat{\Phi}_{k+1}^* + \frac{L D^2 \delta_k^2}{2} + (1 - \delta_k) \xi_k \\
		& = \Phi_{k+1}^* + ( \hat{\Phi}_{k+1}^* - \Phi_{k+1}^*) + \frac{L D^2 \delta_k^2}{2} + (1 - \delta_k) \xi_k. \nonumber
	\end{align}
	Then, we have from Lemma \ref{lemma.vstar} that
	\begin{align}\label{eq.bbb2}
		\hat{\Phi}_{k+1}^* - \Phi_{k+1}^* & =s \hat{\Phi}_{k+1}(\hat{\mathbf{v}}_{k+1}) - \Phi_{k+1}(\mathbf{v}_{k+1}) \\
		& = \hat{\Phi}_{k+1}(\hat{\mathbf{v}}_{k+1})- \hat{\Phi}_{k+1}(\mathbf{v}_{k+1}) + \hat{\Phi}_{k+1}(\mathbf{v}_{k+1}) - \Phi_{k+1}(\mathbf{v}_{k+1}) \nonumber \\
		& \stackrel{(\text{h})}{\leq}  \hat{\Phi}_{k+1}(\mathbf{v}_{k+1}) - \Phi_{k+1}(\mathbf{v}_{k+1}) \nonumber \\
		& \stackrel{(\text{i})}{ = } \delta_k \Big[ f(\mathbf{y}_k) +  \big\langle \nabla f(\mathbf{y}_k) , \mathbf{v}_{k+1}  - \mathbf{y}_k \big\rangle  \Big] -  \delta_k \Big[ f(\mathbf{x}_{k+1}) +  \big\langle \nabla f(\mathbf{x}_{k+1}) , \mathbf{v}_{k+1}  - \mathbf{x}_{k+1} \big\rangle  \Big] \nonumber \\
		& \stackrel{(\text{j})}{\leq} \delta_k \big\langle \nabla f(\mathbf{y}_k) - \nabla f(\mathbf{x}_{k+1}), \mathbf{v}_{k+1}  - \mathbf{x}_{k+1}   \big\rangle \nonumber \\ 
		& \leq \delta_k \big\| \nabla f(\mathbf{y}_k) - \nabla f(\mathbf{x}_{k+1})\big\|_* \big\|  \mathbf{v}_{k+1}  - \mathbf{x}_{k+1} \big\|  \nonumber \\
		& \stackrel{(\text{k})}{ \leq }  \delta_k  L \big\| \mathbf{y}_k - \mathbf{x}_{k+1}\big\| \big\|  \mathbf{v}_{k+1}  - \mathbf{x}_{k+1} \big\| \nonumber \\
		& \stackrel{(\text{l})}{ \leq }  \delta_k^2  L \big\| \mathbf{v}_k - \hat{\mathbf{v}}_{k+1}\big\| \big\|  \mathbf{v}_{k+1}  - \mathbf{x}_{k+1} \big\|  \leq \delta_k^2  L D^2 \nonumber
	\end{align}
	where (h) is because $ \hat{\Phi}_{k+1}(\hat{\mathbf{v}}_{k+1}) \leq  \hat{\Phi}_{k+1}(\mathbf{x}), \forall \mathbf{x} \in {\cal X}$ according to Lemma \ref{lemma.vstar}; (i) follows from \eqref{eq.phi}; (j) uses $f(\mathbf{y}_k) - f(\mathbf{x}_{k+1}) \leq \langle  \nabla f(\mathbf{y}_k), \mathbf{y}_k - \mathbf{x}_{k+1} \rangle$; (k) is because of Assumption \ref{as.1}; and (l) uses the choice of $\mathbf{y}_k$ and $\mathbf{x}_{k+1}$. Plugging \eqref{eq.bbb2} back into \eqref{eq.bbb1}, we have
	\begin{align*}
		f(\mathbf{x}_{k+1})	\leq \Phi_{k+1}^* + \frac{3 L D^2 \delta_k^2}{2} + (1 - \delta_k) \xi_k
	\end{align*}
	which completes the proof.
\end{proof}

\subsection{Proof of Theorem \ref{thm.general}}
\begin{proof}
	Given $\big( \{\Phi_k(\mathbf{x}) \}_{k=0}^\infty, \{ \lambda_k \}_{k=0}^\infty \big)$ is an ES as shown in Lemma \ref{lemma.es}, together with the fact $f(\mathbf{x}_k) \leq \min_{\mathbf{x} \in {\cal X}} \Phi_k(\mathbf{x}) + \xi_k, \forall k$ as shown in Lemma \ref{lemma_phistar}, one can directly apply Lemma \ref{lemma.es_convergence} to have
	\begin{align}\label{eq.new2}
		f(\mathbf{x}_k) - f(\mathbf{x}^*) & \leq  \lambda_k \big( f (\mathbf{x}_0) - f(\mathbf{x}^*) \big)  + \xi_k  = \frac{ 2\big( f(\mathbf{x}_0) - f(\mathbf{x}^*) \big) }{(k+1)(k+2)} + \xi_k
	\end{align}
	where $\xi_k$ is defined in Lemma \ref{lemma_phistar}. Clearly, $\xi_k \geq 0, \forall k$, and one can find an upper bound of it as
	\begin{align*}
		\xi_k & = (1- \delta_{k-1}) \xi_{k-1} + \frac{3 \delta_{k-1}^2}{2} LD^2 \nonumber \\
		& = \frac{3 LD^2}{2} \sum_{\tau=0}^{k-1} \delta_\tau^2 \bigg[ \prod_{j=\tau+1}^{k-1} (1 - \delta_j) \bigg] \\
		& = \frac{3LD^2}{2} \sum_{\tau=0}^{k-1} \frac{4}{(\tau + 3)^2} \frac{(\tau+2)(\tau+3)}{(k+1)(k+2)} \leq \frac{6 LD^2}{k+2}.
	\end{align*}
	Plugging $\xi_k$ into \eqref{eq.new2} completes the proof.
\end{proof}

\subsection{Stopping Criterion}\label{apdx.stop}
In this subsection we show that the value of $f(\mathbf{x}_k) - \Phi_k^*$ can be used to derive a stopping criterion (see \eqref{eq.stop_crit}). How to obtain the value of $\Phi_k^*$ iteratively (via \eqref{eq.stop_crit2} and \eqref{eq.stop_crit3}) is also discussed. 

First, as a consequence of Lemma \ref{lemma_phistar}, we have $f(\mathbf{x}_k) - \Phi_k^* \leq \xi_k = {\cal O}\big(\frac{LD^2}{k}\big)$. This means that the value of $f(\mathbf{x}_k) - \Phi_k^*$ converges to $0$ at the same rate of $ f(\mathbf{x}_k) - f(\mathbf{x}^*)$.

Next we show that how to estimate $f(\mathbf{x}_k) - f(\mathbf{x}^*) $ using $f(\mathbf{x}_k) - \Phi_k^*$. We have that
\begin{align*}
	f(\mathbf{x}_k) - \Phi_k^* & \stackrel{(\text{a})}{\geq} 	f(\mathbf{x}_k) -  \Phi_k(\mathbf{x}^*) \stackrel{(\text{b})}{\geq} f(\mathbf{x}_k) - (1 - \lambda_k)	 f(\mathbf{x}^*) - \lambda_k \Phi_0 (\mathbf{x}^*) \nonumber \\
	& \stackrel{(\text{c})}{=} (1 - \lambda_k)	\big[  f(\mathbf{x}_k) -f(\mathbf{x}^*) \big]  + \lambda_k 	\big[  f(\mathbf{x}_k) -f(\mathbf{x}_0) \big]
\end{align*}
where (a) is because of $\Phi_k^* = \min_{\mathbf{x}\in {\cal X}} \Phi_k(\mathbf{x})$; (b) is by the definition of ES; and (c) uses $\Phi_0(\mathbf{x}) \equiv f(\mathbf{x}_0)$. The inequality above implies that
\begin{align}\label{eq.stop_crit}
	f(\mathbf{x}_k) - f(\mathbf{x}^*) \leq  \frac{1}{1-\lambda_k} \bigg( f(\mathbf{x}_k) - \Phi_k^*  - \lambda_k 	\big[  f(\mathbf{x}_k) -f(\mathbf{x}_0)\big]  \bigg).
\end{align}
Notice that the RHS of \eqref{eq.stop_crit} goes to $0$ as $k$ increases, hence \eqref{eq.stop_crit} can be used as the stopping criterion. 

Finally we discuss how to update $\Phi_k^*$ efficiently. From \eqref{eq.stop_use1}, we have 
	\begin{align*}
		& \Phi_{k+1}(\mathbf{x})  =  f(\mathbf{x}_0)  \prod_{\tau = 0}^k (1-\delta_\tau) + \sum_{\tau=0}^k \delta_\tau \bigg[f(\mathbf{x}_{\tau+1}) + \langle \nabla f(\mathbf{x}_{\tau+1}), \mathbf{x} - \mathbf{x}_{\tau + 1} \rangle \bigg] \prod_{j = \tau+1}^k (1-\delta_j)  \nonumber \\
		& =  f(\mathbf{x}_0)  \prod_{\tau = 0}^k (1-\delta_\tau) + \sum_{\tau=0}^k \delta_\tau \bigg[f(\mathbf{x}_{\tau+1}) + \langle \nabla f(\mathbf{x}_{\tau+1}), \mathbf{x} - \mathbf{x}_{\tau + 1} \rangle \bigg] \prod_{j = \tau+1}^k (1-\delta_j) \nonumber \\
		& =  f(\mathbf{x}_0) \prod_{\tau = 0}^k (1-\delta_\tau) + \sum_{\tau=0}^k \delta_\tau \bigg[f(\mathbf{x}_{\tau+1}) - \langle \nabla f(\mathbf{x}_{\tau+1}), \mathbf{x}_{\tau + 1} \rangle \bigg] \prod_{j = \tau+1}^k (1-\delta_j) + \langle \mathbf{g}_{k+1}, \mathbf{x}\rangle
	\end{align*}
where the last equation uses the definition of $\mathbf{g}_{k+1}$. Hence, we can obtain $\Phi_{k+1}^* $ as
\begin{align}\label{eq.stop_crit2}
	\Phi_{k+1}^* = \Phi_{k+1}(\mathbf{v}_{k+1}) = V_{k+1} + \langle \mathbf{g}_{k+1}, \mathbf{v}_{k+1}\rangle
\end{align}
and $V_{k+1}$ can be updated as
\begin{align}\label{eq.stop_crit3}
	V_{k+1} = (1-\delta_k) V_k + \delta_k \Big[f(\mathbf{x}_{k+1}) - \langle \nabla f(\mathbf{x}_{k+1}), \mathbf{x}_{k + 1} \rangle \Big], ~~{\rm with}~~ V_0 = f(\mathbf{x}_0).
\end{align}


\section{Proof of Theorem \ref{thm.acc}}
Because we are dealing with an $\ell_2$ norm ball constraint in this section, we use $R:= \frac{D}{2}$ for convenience. And we will extend the domain of $f(\mathbf{x})$ slightly to $\tilde{\cal X}:= {\rm conv} \{  \mathbf{x} - \frac{1}{L}  \nabla f(\mathbf{x}),~ \forall \mathbf{x} \in {\cal X} \}$, i.e., $f: \tilde{\cal X} \rightarrow \mathbb{R}$. This is a very mild assumption since most of practically used loss functions have domain $\mathbb{R}^d$.

\begin{lemma}\label{lemma.nesterov}
	\citep[Theorem 2.1.5]{nesterov2004} If Assumptions \ref{as.1} and \ref{as.2} hold with the extended domain $\tilde{\cal X}$, then it is true that for any $\mathbf{x}, \mathbf{y} \in {\cal X}$
	\begin{align*}
		\frac{1}{2L} \| \nabla f(\mathbf{x}) - \nabla f(\mathbf{y}) \|_2^2 \leq f(\mathbf{y}) -  f(\mathbf{x})	- \langle \nabla f(\mathbf{x}), \mathbf{y} - \mathbf{x} \rangle.
	\end{align*}
\end{lemma}

\begin{lemma}\label{lemma.fw_y_grad}
	Choose $\delta_k = \frac{2}{k+3}$, then we have
	\begin{align*}
		 \| \nabla f(\mathbf{x}_k) - \nabla f(\mathbf{x}^*) \|_2  \leq \sqrt{ \frac{ 4L \big( f(\mathbf{x}_0) - f(\mathbf{x}^*) \big) }{(k+1)(k+2)} + \frac{12 L^2 D^2}{k+2} } \leq \frac{C_1}{\sqrt{k+2}}
	\end{align*}
	where $C_1 \leq \sqrt{12L^2D^2 + 4L \big( f(\mathbf{x}_0) - f(\mathbf{x}^*) \big)} $.
\end{lemma}
\begin{proof}
Using Lemma \ref{lemma.nesterov}, we have
	\begin{align*}
		\frac{1}{2L} \| \nabla f(\mathbf{x}_k) - \nabla f(\mathbf{x}^*) \|_2^2 & \leq f(\mathbf{x}_k) -  f(\mathbf{x}^*)	- \langle \nabla f(\mathbf{x}^*), \mathbf{x}_k - \mathbf{x}^* \rangle	 \stackrel{(\text{a})}{\leq} f(\mathbf{x}_k) - f(\mathbf{x}^*) \\
		& \stackrel{(\text{b})}{\leq } \frac{ 2\big( f(\mathbf{x}_0) - f(\mathbf{x}^*) \big) }{(k+1)(k+2)} + \frac{6 LD^2}{k+2}
	\end{align*}
	where (a) is by the optimality condition, that is, $\langle \nabla f(\mathbf{x}^*), \mathbf{x} - \mathbf{x}^* \rangle \geq 0, \forall \mathbf{x} \in {\cal X}$; and (b) is by Theorem \ref{thm.general}. This further implies 
	\begin{align*}
		 \| \nabla f(\mathbf{x}_k) - \nabla f(\mathbf{x}^*) \|_2  \leq \sqrt{ \frac{ 4L \big( f(\mathbf{x}_0) - f(\mathbf{x}^*) \big) }{(k+1)(k+2)} + \frac{12 L^2 D^2}{k+2} }.
	\end{align*}
	The proof is thus completed.
\end{proof}

\begin{lemma}\label{lemma.equal}
	If both $\mathbf{x}_1^*$ and $\mathbf{x}_2^*$ minimize $f(\mathbf{x})$ over ${\cal X}$, then we have $ \nabla f(\mathbf{x}_1^*) = \nabla f(\mathbf{x}_2^*) $.
\end{lemma}
\begin{proof}
	From Lemma \ref{lemma.nesterov}, we have
	\begin{align*}
		\frac{1}{2L} \| \nabla f(\mathbf{x}_2^*) - \nabla f(\mathbf{x}_1^*) \|_2^2 & \leq f(\mathbf{x}_2^*) -  f(\mathbf{x}_1^*)	- \langle \nabla f(\mathbf{x}_1^*), \mathbf{x}_2^* - \mathbf{x}_1^* \rangle	 \stackrel{(\text{a})}{\leq} f(\mathbf{x}_2^*) - f(\mathbf{x}_1^*) = 0
	\end{align*}
	where (a) is by the optimality condition, that is, $\langle \nabla f(\mathbf{x}_1^*), \mathbf{x} - \mathbf{x}_1^* \rangle \geq 0, \forall \mathbf{x} \in {\cal X}$. Hence we can only have $\nabla f(\mathbf{x}_2^*) = \nabla f(\mathbf{x}_1^*)$. This means that the value of $\nabla f(\mathbf{x}^*)$ is unique regardless of the uniqueness of $\mathbf{x}^*$.
\end{proof}

\begin{lemma}\label{lemma.theta_gradx*}
Let $\| \nabla f(\mathbf{x}^*) \|_2 = G^*$, (and $G^*$ is unique bacause of Lemma \ref{lemma.equal}) where $G^* \geq G$. Choose $\delta_k = \frac{2}{k+3}$, it is guaranteed to have
	\begin{align*}
		\|  \mathbf{g}_{k+1} - \nabla f(\mathbf{x}^*) \|_2 \leq \frac{4 C_1}{3(\sqrt{k+3}-1)}  + \frac{2 G^*}{(k+2)(k+3)}.
	\end{align*}
	In addition, there exists a constant $C_2 \leq \frac{4}{3}C_1 + \frac{2}{3(\sqrt{3}+1)} G^*$ such that 
	\begin{align*}
		\|  \mathbf{g}_{k+1} - \nabla f(\mathbf{x}^*) \|_2 \leq \frac{C_2}{\sqrt{k+3}-1} .
	\end{align*}
\end{lemma}
\begin{proof}
	First we have
	\begin{align}\label{eq.theta_avg}
		\mathbf{g}_{k+1} & = ( 1 - \delta_k) 	\mathbf{g}_k + \delta_k \nabla f(\mathbf{x}_{k+1})  = \sum_{\tau = 0}^k \delta_\tau \nabla f(\mathbf{x}_{\tau+1}) \bigg[ \prod_{j = \tau+1}^k (1 - \delta_j) \bigg] \\
		& = \sum_{\tau=0}^k \frac{2(\tau+2)}{(k+2)(k+3)} \nabla f(\mathbf{x}_{\tau+1}). \nonumber
	\end{align}
	Noticing that $2\sum_{\tau=0}^k (\tau+2) = (k+1)(k+4) = (k+2)(k+3) - 2$, we have
	\begin{align*}
		\|  \mathbf{g}_{k+1} - \nabla f(\mathbf{x}^*) \|_2 & = \bigg{\|} \sum_{\tau = 0}^k \frac{2(\tau+2)}{(k+2)(k+3)} \big[\nabla f(\mathbf{x}_{\tau+1}) - \nabla f(\mathbf{x}^*) \big] - \frac{2}{(k+2)(k+3)} \nabla f(\mathbf{x}^*) \bigg{\|}_2  \\
		& \leq \sum_{\tau = 0}^k \frac{2(\tau+2)}{(k+2)(k+3)}  \big{\|} \nabla f(\mathbf{x}_{\tau+1}) - \nabla f(\mathbf{x}^*) \big{\|}_2 + \frac{2}{(k+2)(k+3)} \big{\|}  \nabla f(\mathbf{x}^*) \big{\|}_2  \\
		& \stackrel{(\text{a})}{\leq} \sum_{\tau = 0}^k \frac{2(\tau+2)}{(k+2)(k+3)} \frac{C_1}{\sqrt{\tau+3}} + \frac{2 G^*}{(k+2)(k+3)}  \\
		& \leq \frac{2 C_1}{(k+2)(k+3)} \sum_{\tau = 0}^k \sqrt{\tau+2} + \frac{2 G^*}{(k+2)(k+3)} \\
		& \leq \frac{4 C_1}{3(k+2)(k+3)}  (k+3)^{3/2}  + \frac{2 G^*}{(k+2)(k+3)} \\
		& = \frac{4 C_1}{3(\sqrt{k+3}+1)(\sqrt{k+3}-1)}  \sqrt{k+3}  + \frac{2 G^*}{(k+2)(k+3)} \\ 
		& \leq \frac{4 C_1}{3(\sqrt{k+3}-1)}  + \frac{2 G^*}{(k+2)(k+3)} 
	\end{align*}
 	where (a) follows from Lemma \ref{lemma.fw_y_grad}. This completes the proof for the first part of this lemma. Next, to find $C_2$, we have 
 	\begin{align*}
 		\|  \mathbf{g}_{k+1} - \nabla f(\mathbf{x}^*) \|_2 & \leq \frac{4 C_1}{3(\sqrt{k+3}-1)}  + \frac{2 G^*}{(k+2)(k+3)} \\
 		& = \frac{4 C_1}{3(\sqrt{k+3}-1)}  + \frac{2 G^*}{(k+3)(\sqrt{k+3}+1)(\sqrt{k+3}-1) } \\
 		& \stackrel{(\text{b})}{\leq} \frac{4 C_1}{3(\sqrt{k+3}-1)}  + \frac{2 G^*}{3(\sqrt{3}+1)(\sqrt{k+3}-1) }
 	\end{align*}
 	where in (b) we use $k+3 \geq 3$ and  $\sqrt{k+3} +1 \geq \sqrt{3}+1$. The proof is thus completed.
\end{proof}

\begin{lemma}\label{lemma.last_lemma}
	There exists a constant $T_1 \leq  \big(\frac{2C_2}{G^*} + 1 \big)^2 - 3$, such that $\| \mathbf{g}_{k+1} \|_2 \geq \frac{G^*}{2}, \forall k \geq T_1$. 
\end{lemma}
\begin{proof}
	Consider a specific $\tilde{k}$ with $\| \mathbf{g}_{\tilde{k}+1} \|_2 < \frac{G^*}{2}$ satisfied. In this case we have 
	\begin{align*}
		\| \mathbf{g}_{\tilde{k}+1} - \nabla f(\mathbf{x}^*) \|_2 \geq \| \nabla f(\mathbf{x}^*) \|_2 - \| \mathbf{g}_{\tilde{k}+1} \|_2 > G^* - \frac{G^*}{2}  = \frac{G^* }{2}.
	\end{align*}
	From Lemma \ref{lemma.theta_gradx*}, we have
	\begin{align*}
		\frac{G^* }{2} < \|  \mathbf{g}_{\tilde{k}+1} - \nabla f(\mathbf{x}^*) \|_2 \leq \frac{C_2}{\sqrt{\tilde{k}+3}-1} .
	\end{align*}
	From this inequality we can observe that $\|\mathbf{g}_{\tilde{k}+1}\|_2$ can be less than $\frac{\sqrt{G}}{2}$ only when $\tilde{k} < T_1 = \big(\frac{2C_2}{G^*} + 1 \big)^2 - 3$. Hence, this lemma is proved.
\end{proof}

\begin{lemma}\label{lemma.k_greater_thres}
	Let $T:= \max\{T_1, T_2\}$, with $T_2 = \sqrt{\frac{8LD}{G^*}} - 3$. When $k \geq T+1$, it is guaranteed that
	\begin{align}
		\|  \mathbf{v}_{k+1} - \hat{\mathbf{v}}_{k+1} \|_2 \leq 	\frac{\delta_k^3 LDC_3}{\| \mathbf{g}_{k+1} \|_2 \| \mathbf{g}_k \|_2  } \leq \frac{4 \delta_k^3 LDC_3}{(G^*)^2  }
	\end{align}
	where $C_3:= LD^2 + \frac{DC_2}{\sqrt{2}-1} $.
\end{lemma}
\begin{proof}
First we show that when $k \geq T+1$, both $\|\mathbf{g}_k \|_2 > 0$ and $\|\hat{\mathbf{g}}_{k+1} \|_2 > 0$. First, because $k \geq T+1 \geq T_1+1$, through Lemma \ref{lemma.last_lemma} we have $\| \mathbf{g}_k\|_2 \geq \frac{G^*}{2} > 0$. Then we have
\begin{align*}
	\big\| \hat{\mathbf{g}}_{k+1} \big\|_2 & = \big\| (1-\delta_k)\mathbf{g}_k + \delta_k \nabla f(\mathbf{x}_{k+1}) - \delta_k \nabla f(\mathbf{x}_{k+1}) +\delta_k \nabla f(\mathbf{y}_k) \big\|_2 \\
	& \geq \big\| \mathbf{g}_{k+1} \big\|_2 - \delta_k \big\| \nabla f(\mathbf{x}_{k+1}) -\nabla f(\mathbf{y}_k) \big\|_2 \geq \frac{G^*}{2} - \delta_k^2 LD
\end{align*}
the last inequality holds when $k \geq T_1$. Hence when $k \geq \max\{ T_1 , T_2\}+1$, we must have both $\|\mathbf{g}_k \|_2 > 0$ and $\|\hat{\mathbf{g}}_{k+1} \|_2 > 0$. Then for any $k \geq T+1$, in view of  \eqref{eq.acc_opt_v}, we can write
	\begin{align}\label{eq.ccc...}
		& ~~~~~ \|  \mathbf{v}_{k+1} - \hat{\mathbf{v}}_{k+1} \|_2  = \bigg\| - \frac{R}{\| \mathbf{g}_{k+1} \|_2} \mathbf{g}_{k+1} + \frac{R}{\| \hat{\mathbf{g}}_{k+1} \|_2} \hat{\mathbf{g}}_{k+1}  \bigg\|_2   \\
		& = \frac{R}{\| \mathbf{g}_{k+1} \|_2  \| \hat{\mathbf{g}}_{k+1} \|_2 } \bigg\|  \big\| \hat{\mathbf{g}}_{k+1} \big\|_2  \mathbf{g}_{k+1} - \big\| \mathbf{g}_{k+1}\big\|_2   \hat{\mathbf{g}}_{k+1}  \bigg\|_2 \nonumber \\
		& = \frac{R}{\| \mathbf{g}_{k+1} \|_2  \| \hat{\mathbf{g}}_{k+1} \|_2 } \bigg\|  \big\| \hat{\mathbf{g}}_{k+1} \big\|_2  \mathbf{g}_{k+1} - \big\| \hat{\mathbf{g}}_{k+1} \big\|_2  \hat{\mathbf{g}}_{k+1} + \big\| \hat{\mathbf{g}}_{k+1} \big\|_2  \hat{\mathbf{g}}_{k+1} - \big\| \mathbf{g}_{k+1}\big\|_2   \hat{\mathbf{g}}_{k+1}  \bigg\|_2 \nonumber \\
		& \leq \frac{R}{\| \mathbf{g}_{k+1} \|_2 } \bigg\|  \mathbf{g}_{k+1} -   \hat{\mathbf{g}}_{k+1}\bigg\|_2  + \frac{R}{\| \mathbf{g}_{k+1} \|_2 }  \bigg| \big\| \hat{\mathbf{g}}_{k+1} \big\|_2  - \big\| \mathbf{g}_{k+1}\big\|_2   \bigg| \nonumber \\
		& \stackrel{(\text{a})}{\leq} \frac{2R}{\| \mathbf{g}_{k+1} \|_2 } \bigg\|  \mathbf{g}_{k+1} -   \hat{\mathbf{g}}_{k+1}\bigg\|_2 = \frac{2R \delta_k}{\| \mathbf{g}_{k+1} \|_2 } \bigg\|  \nabla f(\mathbf{x}_{k+1}) -  \nabla f(\mathbf{y}_k) \bigg\|_2 \nonumber \\
		& \stackrel{(\text{b})}{\leq} \frac{2R L \delta_k}{\| \mathbf{g}_{k+1} \|_2 } \bigg\|  \mathbf{x}_{k+1} -  \mathbf{y}_k \bigg\|_2 = \frac{D L \delta_k^2 }{\| \mathbf{g}_{k+1} \|_2 } \bigg\|  \hat{\mathbf{v}}_{k+1} -  \mathbf{v}_k \bigg\|_2 \nonumber
	\end{align}
	where (a) is by $ \big| \|\mathbf{a} \|_2 -  \|\mathbf{b} \|_2 \big| \leq \big\| \mathbf{a}-  \mathbf{b}  \big\|_2 $; and (b) is by Assumption \ref{as.1}. Then we will bound $\| \hat{\mathbf{v}}_{k+1} -  \mathbf{v}_k \|_2$. 		
	\begin{align*}
		& ~~~~~ \big\|  \hat{\mathbf{v}}_{k+1} -  \mathbf{v}_k \big\|_2  = \bigg\| - \frac{R}{\| \hat{\mathbf{g}}_{k+1} \|_2} \hat{\mathbf{g}}_{k+1} + \frac{R}{\| \mathbf{g}_k \|_2} \mathbf{g}_k  \bigg\|_2 \nonumber \\
		& = \frac{R}{\| \mathbf{g}_k \|_2  \| \hat{\mathbf{g}}_{k+1} \|_2 } \bigg\|  \big\| \mathbf{g}_k \big\|_2  \hat{\mathbf{g}}_{k+1} - \big\| \hat{\mathbf{g}}_{k+1} \big\|_2  \hat{\mathbf{g}}_{k+1} + \big\| \hat{\mathbf{g}}_{k+1} \big\|_2  \hat{\mathbf{g}}_{k+1} - \big\| \hat{\mathbf{g}}_{k+1}\big\|_2   \mathbf{g}_k  \bigg\|_2 \nonumber \\
		& \leq \frac{R}{\| \mathbf{g}_k \|_2  } \bigg| \big\| \mathbf{g}_k \big\|_2  -  \big\| \hat{\mathbf{g}}_{k+1} \big\|_2  \bigg| + \frac{R}{\| \mathbf{g}_k \|_2  } \bigg\| \hat{\mathbf{g}}_{k+1} - \mathbf{g}_k\bigg\|_2 \nonumber \\
		& \stackrel{(\text{c})}{\leq} \frac{D}{\| \mathbf{g}_k \|_2  } \bigg\| \hat{\mathbf{g}}_{k+1} - \mathbf{g}_k\bigg\|_2 = \frac{\delta_k D}{\| \mathbf{g}_k \|_2  } \bigg\| \nabla f(\mathbf{y}_k) - \mathbf{g}_k\bigg\|_2 \nonumber \\
		& \leq \frac{\delta_k D}{\| \mathbf{g}_k \|_2  } \big\| \nabla f(\mathbf{y}_k) - \nabla f(\mathbf{x}^*)\big\|_2  + \frac{\delta_k D}{\| \mathbf{g}_k \|_2  }  \big\|\nabla f(\mathbf{x}^*) -\mathbf{g}_k\big\|_2 \nonumber \\
		& \leq \frac{\delta_k LD^2}{\| \mathbf{g}_k \|_2  } + \frac{\delta_k D}{\| \mathbf{g}_k \|_2  }  \big\|\nabla f(\mathbf{x}^*) -\mathbf{g}_k\big\|_2  \nonumber \\
		& \leq \frac{\delta_k LD^2}{\| \mathbf{g}_k \|_2  } + \frac{\delta_k D}{\| \mathbf{g}_k \|_2  }  \frac{C_2}{\sqrt{k+2} -1 }	 \leq \frac{\delta_k \big( LD^2 + \frac{DC_2}{\sqrt{T+3}-1} \big)}{\| \mathbf{g}_k \|_2  } := \frac{\delta_k C_3}{\| \mathbf{g}_k \|_2 }
	\end{align*}
	where (c) again uses $ \big| \|\mathbf{a} \|_2 -  \|\mathbf{b} |_2 \big| \leq \big\| \mathbf{a}-  \mathbf{b}  \big\|_2 $; and the last inequality is because of Lemma \ref{lemma.fw_y_grad}. Plugging back to \eqref{eq.ccc...}, we arrive at
	\begin{align*}
		\|  \mathbf{v}_{k+1} - \hat{\mathbf{v}}_{k+1} \|_2 \leq 	\frac{D L \delta_k^2 }{\| \mathbf{g}_{k+1} \|_2 } \frac{\delta_k C_3}{\| \mathbf{g}_k \|_2 }=\frac{\delta_k^3 LDC_3}{\| \mathbf{g}_{k+1} \|_2 \| \mathbf{g}_k \|_2  } \leq \frac{4\delta_k^3 LDC_3}{(G^*)^2  }. 
	\end{align*}		
	The proof is thus completed.
\end{proof}

\begin{lemma}\label{lemma.cntpart}
Let $\xi_0 = 0$ and $T$ defined the same as in Lemma \ref{lemma.k_greater_thres}. Denote $\Phi_k^* := \Phi_k(\mathbf{v}_k)$ as the minimum value of $\Phi_k(\mathbf{x})$ over ${\cal X}$, then we have 
	\begin{align*}
		f(\mathbf{x}_k )\leq 	\Phi_k^* + \xi_k, \forall k\geq 0
	\end{align*}
	where for $k < T+1$, $\xi_{k+1} = (1-\delta_k)\xi_k + \frac{3LD^2}{2} \delta_k^2$, and $\xi_{k+1} = C_4 \delta_k^4 + (1 - \delta_k)  \xi_k$ for $k \geq T+1$ with $C_4 = \Big( \frac{C_1}{\sqrt{T+4}} + G^* \Big)\frac{4 LDC_3}{(G^*)^2  }$.
\end{lemma}
\begin{proof}
	The proof for $k < T+1$ is similar as that in Lemma \ref{lemma_phistar}, hence it is omitted here. We mainly focus on the case where $k \geq T+1$.
	\begin{align*}
		\Phi_{k+1}^*  & = \Phi_{k+1} (\mathbf{v}_{k+1}) = (1 - \delta_k) \Phi_k(\mathbf{v}_{k+1}) + \delta_k \Big[ f(\mathbf{x}_{k+1}) +  \big\langle \nabla f(\mathbf{x}_{k+1}) , \mathbf{v}_{k+1}  - \mathbf{x}_{k+1} \big\rangle  \Big] \\
		& \stackrel{(\text{a})}{\geq} (1 - \delta_k) \Phi_k(\mathbf{v}_k) + \delta_k \Big[ f(\mathbf{x}_{k+1}) +  \big\langle \nabla f(\mathbf{x}_{k+1}) , \mathbf{v}_{k+1}  - \mathbf{x}_{k+1} \big\rangle  \Big] \\
		& \geq (1 - \delta_k) f(\mathbf{x}_k) + \delta_k \Big[ f(\mathbf{x}_{k+1}) +  \big\langle \nabla f(\mathbf{x}_{k+1}) , \mathbf{v}_{k+1}  - \mathbf{x}_{k+1} \big\rangle  \Big]  - (1 - \delta_k)  \xi_k \nonumber \\
		& = f(\mathbf{x}_{k+1}) + (1 - \delta_k) \big[ f(\mathbf{x}_k) - f(\mathbf{x}_{k+1})  \big] + \delta_k  \big\langle \nabla f(\mathbf{x}_{k+1}) , \mathbf{v}_{k+1}  - \mathbf{x}_{k+1} \big\rangle - (1 - \delta_k)  \xi_k \\
		& \stackrel{(\text{b})}{\geq}  f(\mathbf{x}_{k+1}) + (1 - \delta_k) \big\langle \nabla f(\mathbf{x}_{k+1}),  \mathbf{x}_k - \mathbf{x}_{k+1} \big\rangle + \delta_k  \big\langle \nabla f(\mathbf{x}_{k+1}) , \mathbf{v}_{k+1}  - \mathbf{x}_{k+1} \big\rangle - (1 - \delta_k)  \xi_k \\
		& =  f(\mathbf{x}_{k+1}) + \delta_k \big\langle \nabla f(\mathbf{x}_{k+1}),  \mathbf{v}_{k+1} - \hat{\mathbf{v}}_{k+1} \big\rangle  - (1 - \delta_k)  \xi_k  \\
		&  \stackrel{(\text{c})}{\geq} f(\mathbf{x}_{k+1}) - \delta_k  \| \nabla f(\mathbf{x}_{k+1}) \|_2 \|  \mathbf{v}_{k+1} - \hat{\mathbf{v}}_{k+1} \|_2 - (1 - \delta_k)  \xi_k  \\
		&  \stackrel{(\text{d})}{\geq} f(\mathbf{x}_{k+1}) -	 \| \nabla f(\mathbf{x}_{k+1}) \|_2 \frac{4 \delta_k^4 LDC_3}{(G^*)^2  }  - (1 - \delta_k)  \xi_k \\
		&  \stackrel{(\text{e})}{\geq} f(\mathbf{x}_{k+1}) -	\Big( \frac{C_1}{\sqrt{T+4}} + G^* \Big)\frac{4 \delta_k^4 LDC_3}{(G^*)^2  }  - (1 - \delta_k)  \xi_k 
	\end{align*}
	where (a) is because $\mathbf{v}_k$ minimizes $\Phi_k(\mathbf{x})$ shown in Lemma \ref{lemma.vstar}; (b) is by $f(\mathbf{x}_{k+1}) - f(\mathbf{x}_k) \leq \langle \nabla f(\mathbf{x}_{k+1}), \mathbf{x}_{k+1} - \mathbf{x}_k \rangle$; (c) uses Cauchy-Schwarz inequality; (d) uses Lemma \ref{lemma.k_greater_thres}, and (e) uses the following inequality.
	\begin{align*}
		\| \nabla f(\mathbf{x}_{k+1}) \|_2 & = \| \nabla f(\mathbf{x}_{k+1}) - \nabla f(\mathbf{x}^*) + \nabla f(\mathbf{x}^*) \|_2 \nonumber \\
		& \leq \| \nabla f(\mathbf{x}_{k+1}) - \nabla f(\mathbf{x}^*)\|_2 + \| \nabla f(\mathbf{x}^*) \|_2 \nonumber \\
		& \leq \frac{C_1}{\sqrt{k+3}} + G^* \leq \frac{C_1}{\sqrt{T+4}} + G^*.
	\end{align*}
	where the last line uses Lemma \ref{lemma.fw_y_grad}.
\end{proof}

\noindent\textbf{Proof of Theorem \ref{thm.acc}}
\begin{proof}
Let $T$ be defined the same as in Lemma \ref{lemma.last_lemma}. For convenience denote $\xi_{k+1} = (1-\delta_k)\xi_k + \theta_k$. When $k < T +1$, we have $\theta_k = \frac{3LD^2}{2} \delta_k^2$; when $k \geq T+1$, we have $\theta_k = C_4 \delta_k^4$.
	
Then we can write 
	\begin{align}\label{eq.apdx.xi_thm2}
		\xi_{k+1} & = (1-\delta_k) \xi_k +  \theta_k = \sum_{\tau=0}^k \theta_\tau \prod_{j=\tau+1}^k (1-\delta_j) \nonumber \\
		& =  \sum_{\tau=0}^k \theta_\tau \frac{(\tau+2)(\tau+3)}{(k+2)(k+3)}  \nonumber \\
		& = \sum_{\tau=0}^{T} \frac{3LD^2}{2} \delta_{\tau}^2 \frac{(\tau+2)(\tau+3)}{(k+2)(k+3)} + \sum_{\tau=T+1}^{k} C_4 \delta_\tau^4 \frac{(\tau+2)(\tau+3)}{(k+2)(k+3)} \nonumber \\
		& = \frac{6LD^2 (T+1)}{(k+2)(k+3)} + {\cal O}\bigg( \frac{C_4}{k^3} \bigg).
	\end{align}

	Again note that $T < {\cal O} \big( \max\{ \sqrt{\frac{LD}{G}}, \frac{L^2D^2}{G^2} \} \big)$ is a constant independent of $k$. Finally, applying Lemma \ref{lemma.es_convergence}, we have
	\begin{align}\label{eq.apdx.conv_thm2}
		f(\mathbf{x}_k) - f(\mathbf{x}^*) \leq \frac{2 \big[ f(\mathbf{x}_0) - f(\mathbf{x}^*)\big]}{(k+1)(k+2)} + \xi_k.
	\end{align}
	Plugging the expression of $\xi_k$, i.e., \eqref{eq.apdx.xi_thm2}, into \eqref{eq.apdx.conv_thm2} completes the proof.
\end{proof}

\section{Discussions for Other Constraints}
\subsection{$\ell_1$ norm ball}\label{apdx.ell1}
In this subsection we focus on the convergence of ExtraFW for $\ell_1$ norm ball constraint under the assumption that $\argmax_j \big{|}[ \nabla f(\mathbf{x}^*) ]_j \big{|}$ has cardinality $1$ (which is also known as \textit{strict complementarity} \citep{garber2020}, and it naturally implies that the constraint is active). Note that in this case Lemma \ref{lemma.equal} still holds hence the value of $ \nabla f(\mathbf{x}^*)$ is unique regardless the uniqueness of $\mathbf{x}^*$. This assumption directly leads to $\argmax_j \big{|}[ \nabla f(\mathbf{x}^*) ]_j \big{|} - | [ \nabla f(\mathbf{x}^*) ]_i| \geq \lambda, \forall i$ for some $\lambda>0$. 


The closed-form solution of $\mathbf{v}_{k+1}$ is given in \eqref{eq.fw_l1}. The constants required in the proof is summarized below for clearance. The norm considered in this subsection for defining $L$ and $D$ is $\|\cdot \|_1$, that is, $\| \nabla f(\mathbf{x}) - \nabla f(\mathbf{y}) \|_\infty \leq L \| \mathbf{x} - \mathbf{y} \|_1$, and $\| \mathbf{x} - \mathbf{y} \|_1 \leq D, \forall \mathbf{x},\forall \mathbf{y} \in \tilde{\cal X}$. Using equivalences of norms, we also assume $\| \nabla f(\mathbf{x}) - \nabla f(\mathbf{y}) \|_2\leq L_2 \| \mathbf{x} - \mathbf{y} \|_2 ,\forall \mathbf{x},\mathbf{y} \in \tilde{\cal X}$ and $\| \mathbf{x} - \mathbf{y} \|_2 \leq D_2, \forall \mathbf{x},\forall \mathbf{y} \in {\cal X}$.

\begin{lemma}\label{lemma.122}
	There exists a constant $T$ (which is irreverent with $k$), whenever $k \geq T$, it is guaranteed to have 
	\begin{align*}
			\|  \mathbf{v}_{k+1} - \hat{\mathbf{v}}_{k+1} \|_1 = 0
	\end{align*}
\end{lemma}
\begin{proof}
	In the proof, we denote $i = \argmax_{j} | [\nabla f(\mathbf{x}^*)]_j|$ for convenience. With $\| \nabla f(\mathbf{x}^*) \|_2 = G^*$ Lemma \ref{lemma.theta_gradx*} still holds. 
	
	We first show that there exist $T_1 = (\frac{3C_2}{\lambda}+1)^2 - 3$, such that for all $k \geq T_1$, we have $\argmax_j | [\mathbf{g}_{k+1}]_j | = i$, which further implies only the $i$-th entry of $\mathbf{v}_{k+1}$ is non-zero. Since Lemma \ref{lemma.theta_gradx*} holds, one can see whenever $k \geq T_1$, it is guaranteed to have $\| \mathbf{g}_{k+1} - \nabla f(\mathbf{x}^*) \|_2 \leq \frac{\lambda}{3}$. Therefore, one must have $\big| | [\mathbf{g}_{k+1}]_j | - | [\nabla f(\mathbf{x}^*)]_j | \big| \leq \frac{\lambda}{3}, \forall j$. Then it is easy to see that $ | [\mathbf{g}_{k+1}]_i | - | [\mathbf{g}_{k+1}]_j |  \geq \frac{\lambda}{3}, \forall j$. Hence, we have $\argmax_j | [\mathbf{g}_{k+1}]_j | = i$. 
	
	Next we show that there exists another constant $T =  \max\{ T_1,  (\frac{3C_5}{\lambda})^2 - 3 \}$, such that $\argmax_j | [\hat{\mathbf{g}}_{k+1}]_j |  = i, \forall k \geq T$, which further indicates only the $i$-th entry of $\hat{\mathbf{v}}_{k+1}$ is non-zero. In this case, in view of Lemma \ref{lemma.theta_gradx*}, we have
	\begin{align*}
		& ~~~~ \big\| \hat{\mathbf{g}}_{k+1} - \nabla f(\mathbf{x}^*) \big\|_2 = \big\| (1-\delta_k)\mathbf{g}_k + \delta_k \nabla f(\mathbf{x}_{k+1}) - \delta_k \nabla f(\mathbf{x}_{k+1}) +\delta_k \nabla f(\mathbf{y}_k) - \nabla f(\mathbf{x}^*)  \big\|_2 \\
		&  \leq \big\| \mathbf{g}_{k+1} - \nabla f(\mathbf{x}^*)\|_2 + \delta_k  \|  \nabla f(\mathbf{x}_{k+1}) - \nabla f(\mathbf{y}_k)   \big\|_2  \leq \big\| \mathbf{g}_{k+1} - \nabla f(\mathbf{x}^*)\|_2 + \delta_k^2 L_2D_2  \\
		& \leq \frac{C_2}{\sqrt{k+3}-1} + \frac{4L_2D_2}{(k+3)^2} \leq \frac{C_5}{\sqrt{k+3}-1}, \forall k \geq T_1
	\end{align*}
	where $C_5 \leq C_2 + \frac{4L_2D_2}{(\sqrt{T_1+3} -1 )^3}$.

	Hence whenever $k \geq \max\{ T_1,  (\frac{3C_5}{\lambda}+1)^2 - 3 \}$, it is guaranteed to have $\|  \hat{\mathbf{g}}_{k+1} - \nabla f(\mathbf{x}^*) \|_2 \leq \frac{\lambda}{3}$. Therefore, one must have $\big| | [\hat{\mathbf{g}}_{k+1}]_j | - | [\nabla f(\mathbf{x}^*)]_j | \big| \leq \frac{\lambda}{3}, \forall j$. It is thus straightforward to see that $ | [\hat{\mathbf{g}}_{k+1}]_i | - | [\hat{\mathbf{g}}_{k+1}]_j |  \geq \frac{\lambda}{3}, \forall j$. Hence, it is clear that $\argmax_j | [\hat{\mathbf{g}}_{k+1}]_j | = i$. 
	
	Then one can see that when $k \geq T$, we have $\mathbf{v}_{k+1} - \hat{\mathbf{v}}_{k+1} = \mathbf{0}$.
\end{proof}

Next, we modify Lemma \ref{lemma.cntpart} to cope with the $\ell_1$ norm ball constraint.
\begin{lemma}\label{lemma.sssss}
Let $\xi_0 = 0$ and $T$ be the same as in Lemma \ref{lemma.122}. Denote $\Phi_k^* := \Phi_k(\mathbf{v}_k)$ as the minimum value of $\Phi_k(\mathbf{x})$ over ${\cal X}$, then we have 
	\begin{align*}
		f(\mathbf{x}_k )\leq \Phi_k(\mathbf{v}_k) = 	\Phi_k^* + \xi_k, \forall k\geq 0
	\end{align*}
	where for $k < T$, $\xi_{k+1} = (1-\delta_k)\xi_k + \frac{3LD^2}{2} \delta_k^2$, and $\xi_{k+1} = (1 - \delta_k)  \xi_k$ for $k \geq T$.
\end{lemma}
\begin{proof}
	The proof for $k < T$ is similar as that in Lemma \ref{lemma_phistar}, hence it is omitted here. We mainly focus on the case where $k \geq T$. Using similar argument as in Lemma \ref{lemma.cntpart}, we have
	\begin{align*}
		\Phi_{k+1}^* 
		& \geq f(\mathbf{x}_{k+1}) + \delta_k \big\langle \nabla f(\mathbf{x}_{k+1}),  \mathbf{v}_{k+1} - \hat{\mathbf{v}}_{k+1} \big\rangle  - (1 - \delta_k)  \xi_k  \\
		&  = f(\mathbf{x}_{k+1}) - (1 - \delta_k)  \xi_k 
	\end{align*}
	where the last inequality is because of Lemma \ref{lemma.122}.
\end{proof}

\begin{theorem}\label{thm.l1}
	Consider ${\cal X}$ is an $\ell_1$ norm ball. If $\argmax_j \big{|}[ \nabla f(\mathbf{x}^*) ]_j \big{|}$ has cardinality $1$, and Assumptions \ref{as.1} - \ref{as.3} are satisfied, ExtraFW guarantees that
	\begin{align*}
		f(\mathbf{x}_k) - f(\mathbf{x}^*) = {\cal O}\Big(\frac{1}{k^2}\Big).
	\end{align*}
\end{theorem}
\begin{proof}
Let $T$ be defined the same as in Lemma \ref{lemma.122}. For convenience denote $\xi_{k+1} = (1-\delta_k)\xi_k + \theta_k$. When $k < T$, we have $\theta_k = \frac{3LD^2}{2} \delta_k^2$; when $k \geq T$, we have $\theta_k = 0$. Then we can write 
	\begin{align}\label{eq.apdx.xi_thm3}
		\xi_{k+1} & = (1-\delta_k) \xi_k +  \theta_k = \sum_{\tau=0}^k \theta_\tau \prod_{j=\tau+1}^k (1-\delta_j)  =  \sum_{\tau=0}^k \theta_\tau \frac{(\tau+2)(\tau+3)}{(k+2)(k+3)}  \nonumber \\
		& = \sum_{\tau=0}^{T-1} \frac{3LD^2}{2} \delta_{\tau}^2 \frac{(\tau+2)(\tau+3)}{(k+2)(k+3)}   = \frac{6LD^2 T}{(k+2)(k+3)} .
	\end{align}

	Finally, applying Lemma \ref{lemma.es_convergence}, we have
	\begin{align}\label{eq.apdx.conv_thm3}
		f(\mathbf{x}_k) - f(\mathbf{x}^*) \leq \frac{2 \big[ f(\mathbf{x}_0) - f(\mathbf{x}^*)\big]}{(k+1)(k+2)} + \xi_k.
	\end{align}
Plugging the expression of $\xi_k$, i.e., \eqref{eq.apdx.xi_thm3} into \eqref{eq.apdx.conv_thm3} completes the proof.
\end{proof}

\begin{figure}[t]
\centering
	\begin{tabular}{cc}
		\hspace{-0.4cm}
		\includegraphics[width=.4\textwidth]{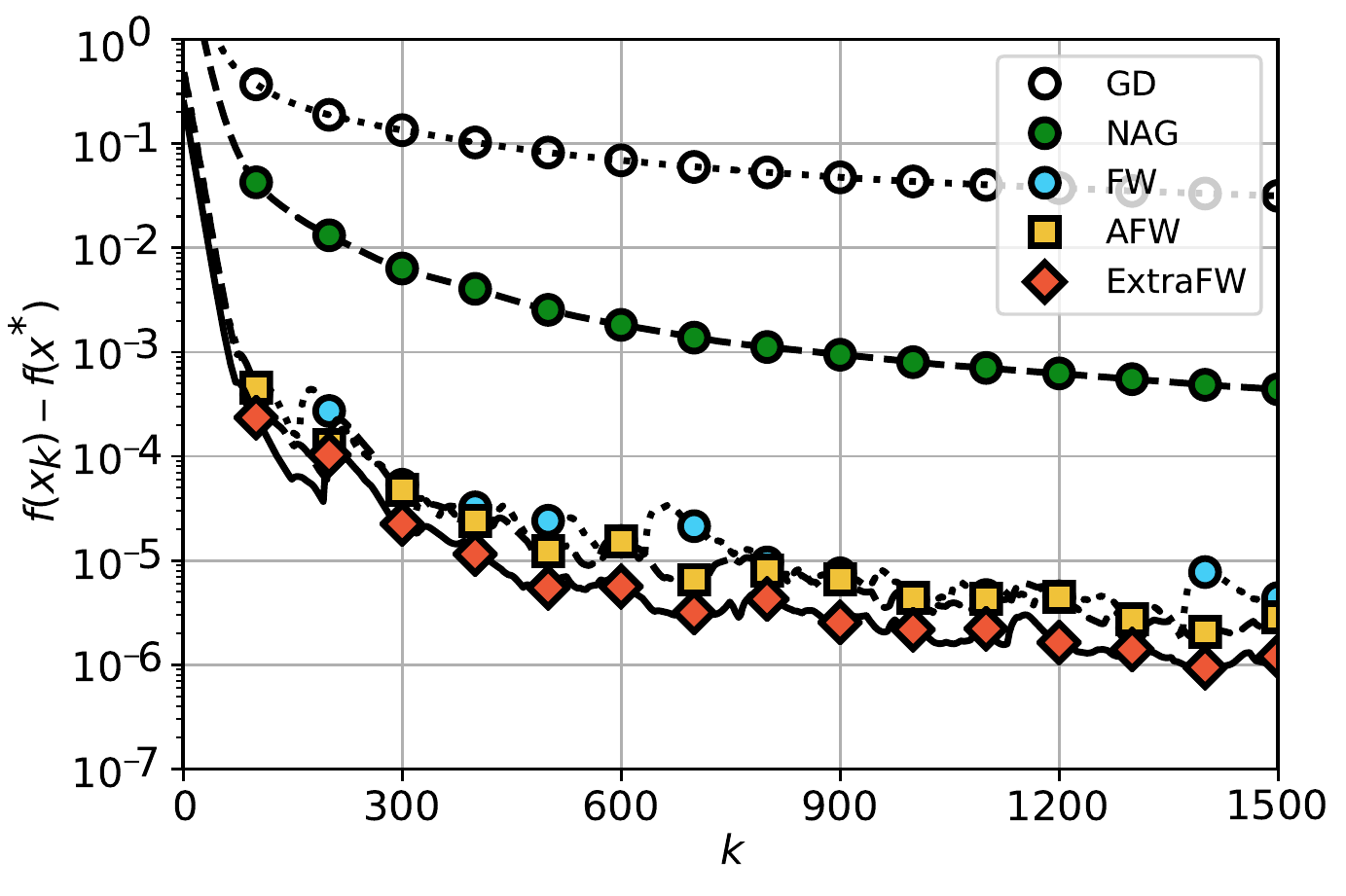}&
		\hspace{-0.3cm}
		\includegraphics[width=.4\textwidth]{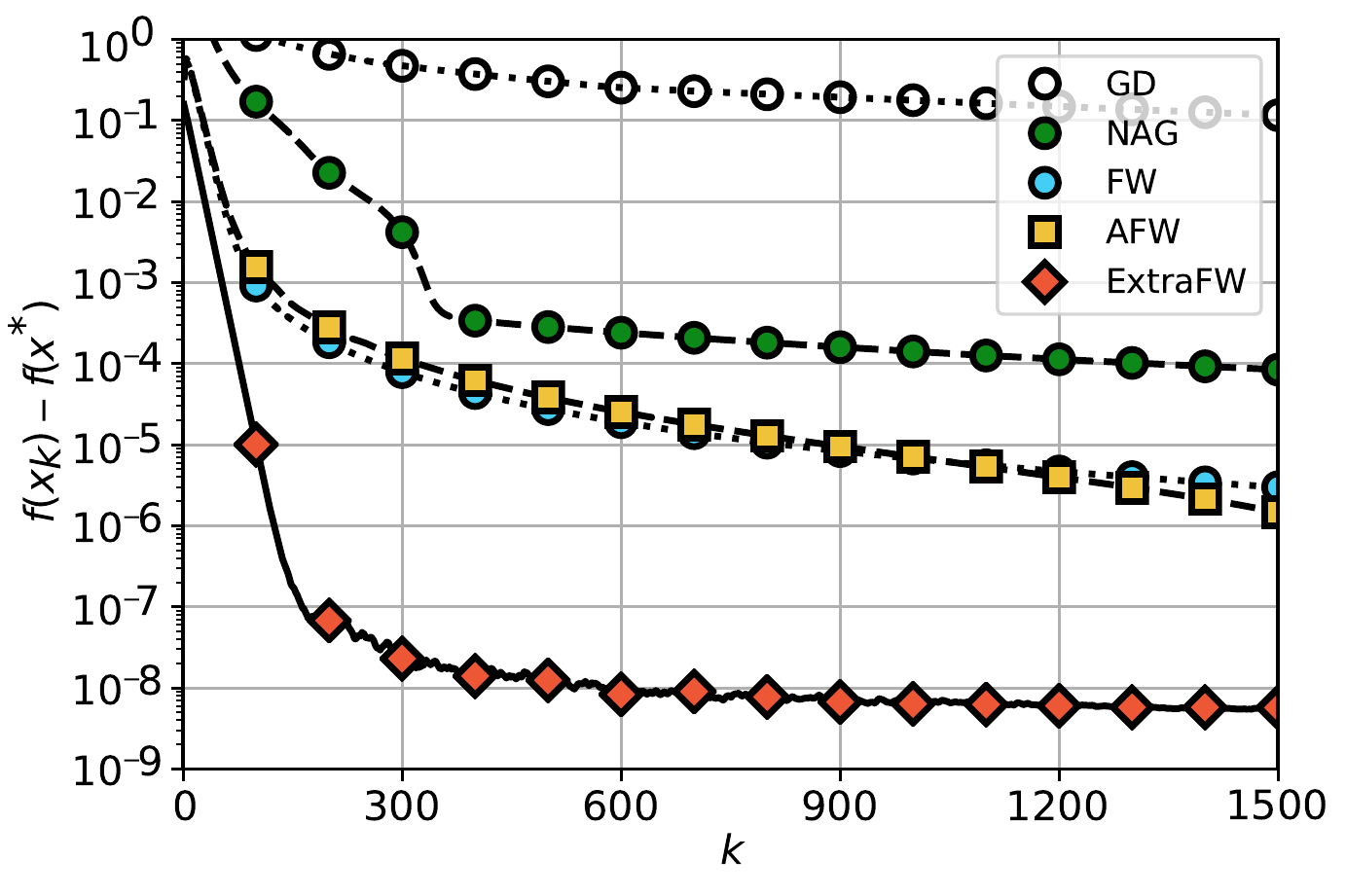}
		\\ (a) \textit{mnist}  & (b) \textit{mushroom} 
	\end{tabular}
	\caption{ExtraFW guarantees an ${\cal O}(\frac{1}{k^2})$ rate on simplex.} 
	 \label{fig.splx}
\end{figure}

\textbf{Beyond $\ell_1$ norm ball.} The ${\cal O}(\frac{T}{k^2})$ rate in Theorem \ref{thm.l1} can be generalized in a straightforward manner to simplex, that is, ${\cal X}:=\{\mathbf{x}| \mathbf{x} \geq \mathbf{0}, \langle \mathbf{1}, \mathbf{x} \rangle = R\}$ for some $R> 0$. A minor assumption needed is that the cardinality of $\argmin_j [ \nabla f(\mathbf{x}^*) ]_j $ is $1$. In this case, the FW steps in ExtraFW admit closed-form solutions. Again taking $\mathbf{v}_{k+1}$ as an example, we have $\mathbf{v}_{k+1} = [0,\ldots, 0, R,0,\ldots, 0]$, where the only non-zero is the $i=\argmin_{j} [\mathbf{g}_{k+1}]_j$-th entry. The proof is similar to the $\ell_1$ norm ball case, i.e., first show that both $\mathbf{g}_{k+1}$ and $\hat{\mathbf{g}}_{k+1}$ converge to $\nabla f(\mathbf{x}^*)$ so that $\mathbf{v}_{k+1} = \hat{\mathbf{v}}_{k+1},\forall k \geq T$, where $T$ is some constant depending on the difference of the smallest and the second smallest entry of $\nabla f(\mathbf{x}^*)$. Then one can follow similar steps of Lemma \ref{lemma.sssss} to obtain the ${\cal O}(\frac{T}{k^2})$ rate. Numerical evidences using logistic regression as objective function can be found in Figure \ref{fig.splx}. Note that in this case however, FW itself converges fast enough.

\subsection{$n$-support norm ball}\label{apdx.nsupp}

When ${\cal X}$ is an $n$-support norm ball, ExtraFW guarantees that $f(\mathbf{x}_k) - f(\mathbf{x}^*) = {\cal O} \big( \frac{T}{k^2} \big) $. The proof is just a combination of Theorem \ref{thm.acc} and \ref{thm.l1}, therefore, we highlight the general idea rather than repeat the proofs step by step. 

The norm considered in this section for defining $L$ and $D$ is $\|\cdot \|_2$, that is, $\| \nabla f(\mathbf{x}) - \nabla f(\mathbf{y}) \|_2 \leq L \| \mathbf{x} - \mathbf{y} \|_2, \forall \mathbf{x},\mathbf{y} \in \tilde{\cal X}$, and $\| \mathbf{x} - \mathbf{y} \|_2 \leq D, \forall \mathbf{x},\mathbf{y} \in {\cal X}$. Besides Assumptions \ref{as.1} - \ref{as.3}, the extra regularity condition we need is that: the $n$-th largest entry of $| [\nabla f(\mathbf{x}^*)]| $ is strictly larger than the $(n+1)$-th largest entry of $| [\nabla f(\mathbf{x}^*)] | $ by $\lambda$. Note that this condition is similar to what we used for the $\ell_1$ norm ball constraint. In addition, this extra assumption directly implies $\| \nabla f(\mathbf{x}^*) \|_2:= G^* >0$. In the proof one may find the constant $G_n^*:= \| {\rm top}_n(\nabla f(\mathbf{x}^*))\|_2$ helpful. Clearly, $G^* \geq G_n^* \geq \sqrt{\frac{n}{d}} G^*$.

\begin{theorem}\label{thm.nsupp}
	Consider ${\cal X}$ is an $n$-support norm ball. If the $n$-th largest entry of $| [\nabla f(\mathbf{x}^*)]| $ is strictly larger than the $(n+1)$-th largest entry of $| [\nabla f(\mathbf{x}^*)] | $, and Assumptions \ref{as.1} - \ref{as.3} are satisfied, ExtraFW guarantees that there exists a constant $T$ such that
	\begin{align*}
		f(\mathbf{x}_k) - f(\mathbf{x}^*) = {\cal O}\Big(\frac{T}{k^2}\Big).
	\end{align*}
\end{theorem}
\begin{proof}

First by using the regularity condition and similar arguments of Lemma \ref{lemma.122}, one can show that there exists a constant $T_1$ (depending on $\lambda$, $L$, $D$, and, $G$) such that the indices of the non-zero entries of $\mathbf{v}_{k+1}$ and $\hat{\mathbf{v}}_{k+1}$ are the same for all $k \geq T_1$. 

Next, using similar arguments of Lemma \ref{lemma.last_lemma}, one can show that there exists a constant $\tilde{T}_2$ such that $\|\mathrm{top}_n (\mathbf{g}_{k+1})\|_2 \geq \frac{G_n^*}{2}$. 

Let $T_2 = \max\{ \tilde{T}_2, T_1\}$. It is clear that for any $k \geq T_2$, the indices of non-zero entries of $\mathbf{v}_{k+1}$ and $\hat{\mathbf{v}}_{k+1}$ are the same. Together with $\|\mathrm{top}_n (\mathbf{g}_{k+1})\|_2 \geq \frac{G_n^*}{2}, \forall k \geq T_2$, we can show that for any $k \geq T_2 +1 $, $ \| \mathbf{v}_{k+1} - \hat{\mathbf{v}}_{k+1} \|_2 = {\cal O}(\delta_k^3)$ holds through similar steps as Lemma \ref{lemma.k_greater_thres}.

Finally, using similar arguments of Lemma \ref{lemma.cntpart} with the aid of  $ \| \mathbf{v}_{k+1} - \hat{\mathbf{v}}_{k+1} \|_2 = {\cal O}(\delta_k^3)$, and applying Lemma \ref{lemma.es_convergence}, we can obtain $f(\mathbf{x}_k) - f(\mathbf{x}^*) = {\cal O} \big( \frac{T_2}{k^2}  \big)$.
\end{proof}






\section{Additional Numerical Results}\label{apdx.numerical}

All numerical experiments are performed using Python 3.7 on an Intel i7-4790CPU @3.60 GHz (32 GB RAM) desktop.

\subsection{Efficiency of ExtraFW: Case Study of $n$-support Norm Ball}\label{apdx.efficiency}
In this subsection we show that ExtraFW achieves fast convergence rate and low iteration cost simultaneously when the constraint set is an $n$-support norm ball. We compare algorithms that can solve the constrained formulation or its equivalent regularized formulation discussed in Section \ref{sec.acc}, that is
\begin{subequations}
\begin{align}
& \min_{\mathbf{x}} ~ f(\mathbf{x}) + \lambda (\| \mathbf{x} \|_{\rm n-sp})^2 \label{eq.prox__} \\ 
	 \Leftrightarrow ~~~& \min_{\mathbf{x}} ~ f(\mathbf{x}) ~~{\rm s.t.} ~~ \| \mathbf{x} \|_{\rm n-sp} \leq R \label{eq.fw__}
\end{align}
\end{subequations}
where $\|\cdot \|_{\rm n-sp}$ denotes the $n$-support norm \citep{argyriou2012}. 


Clearly, one can apply proximal NAG (Prox-NAG) to \eqref{eq.prox__}. The proximal operator per iteration has complexity ${\cal O}(d (n + \log d))$ \citep{argyriou2012}.

One can also apply ExtraFW for \eqref{eq.fw__}. From the Lagrangian duality of \eqref{eq.fw__} and \eqref{eq.prox__}, one can see that if $\lambda \neq 0$, one must have an optimal solution for \eqref{eq.fw__} lies on the boundary of its constraint set. Hence ExtraFW achieves acceleration in this case. Below we summarize the convergence rate and per iteration cost of different algorithms. A simple comparison among different algorithms illustrates the efficiency of ExtraFW. 

\begin{table}[h]
\centering 
\caption{A comparison of different algorithms for logistic regression with  $n$-support norm }\label{tab.n_supp}
 \begin{tabular}{ c*{2}{|c}} 
    \hline
Alg.  & convergence rate  & per iteration cost  \\ \hline
Prox-NAG for \eqref{eq.prox__}  & ${\cal O}(1/k^2)$  & proximal operator: ${\cal O}(d (n + \log d))$  \\ \hline
Projected NAG for \eqref{eq.fw__}& ${\cal O}(1/k^2)$  & projection is expensive \\ \hline
FW for \eqref{eq.fw__} &   ${\cal O}(1/k)$  & FW step: ${\cal O}(d \log n)$  \\ \hline
ExtraFW for \eqref{eq.fw__} &  ${\cal O}(T/k^2)$& FW step: ${\cal O}(d \log n)$  \\ \hline
\end{tabular} 
\end{table}

\subsection{Binary Classification}\label{apdx.classification}

\begin{table}[h]
\centering 
\caption{A summary of datasets used in numerical tests}\label{tab.dataset}
 \begin{tabular}{ c*{3}{|c}} 
    \hline
Dataset  & $d$  & $N$ (train)  & nonzeros  \\ \hline
\textit{w7a}  & $300$  & $24,692$ & $3.89\%$ \\ \hline
\textit{realsim} &   $20,958$  & $50,617$ & $0.24\%$ \\ \hline
\textit{news20} &   $19,996$  & $1,355,191$ & $0.033\%$ \\ \hline
\textit{mushromm} &  $122$ & $8,124$ & $18.75\%$ \\ \hline
\textit{mnist} (digit $4$) &   $784$ & $60,000$ & $12.4\%$ \\ \hline
\end{tabular} 
\end{table}

The datasets used for the tests are summarized in Table \ref{tab.dataset}.

\textbf{Sparsity promoting property of FW variants in $\ell_1$ norm ball constraint.} FW in Alg. \ref{alg.fw} directly promotes sparsity on the solution if it is initialized at $\mathbf{x}_0 = \mathbf{0}$. To see this, suppose that the $i$-th entry of $\nabla f(\mathbf{x}_k)$ has the largest absolute value, then we have $\mathbf{v}_{k+1} = [0, \ldots ,-{\rm sgn}\big( [\nabla f(\mathbf{x}_k)]_i \big) R, \ldots,0]^\top$ with the $i$-th entry being non-zero. Hence, $\mathbf{x}_k$ has at most $k$ non-zero entries given $k-1$ entries are non-zero in $\mathbf{x}_{k-1}$. This sparsity promoting property also holds for ExtraFW.

\begin{figure*}[t]
	\centering
	\begin{tabular}{ccc}
		\hspace{-0.3cm}
		\includegraphics[width=.3\textwidth]{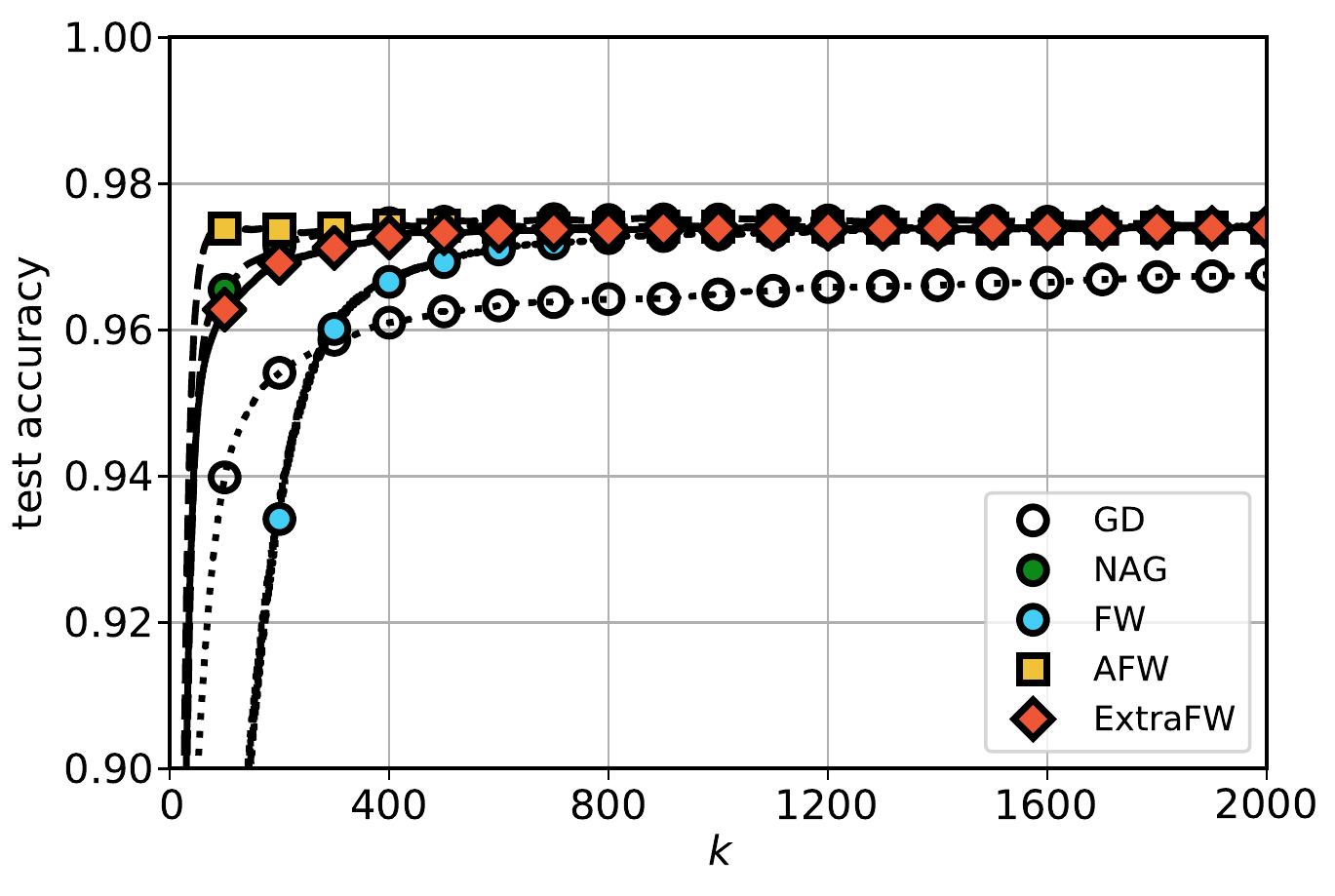}&
		\hspace{-0.3cm}
		\includegraphics[width=.3\textwidth]{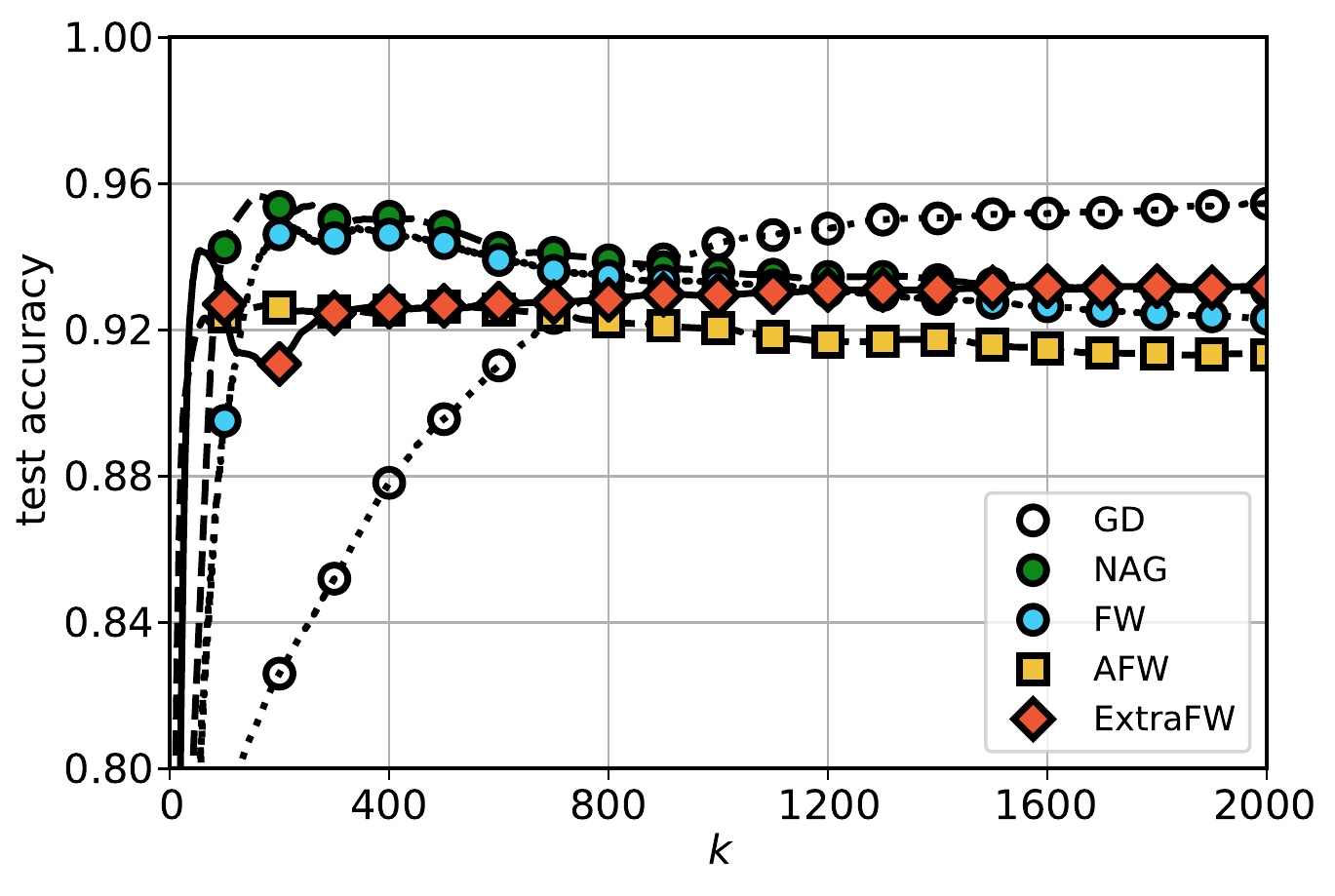}&
		\hspace{-0.3cm}
		\includegraphics[width=.3\textwidth]{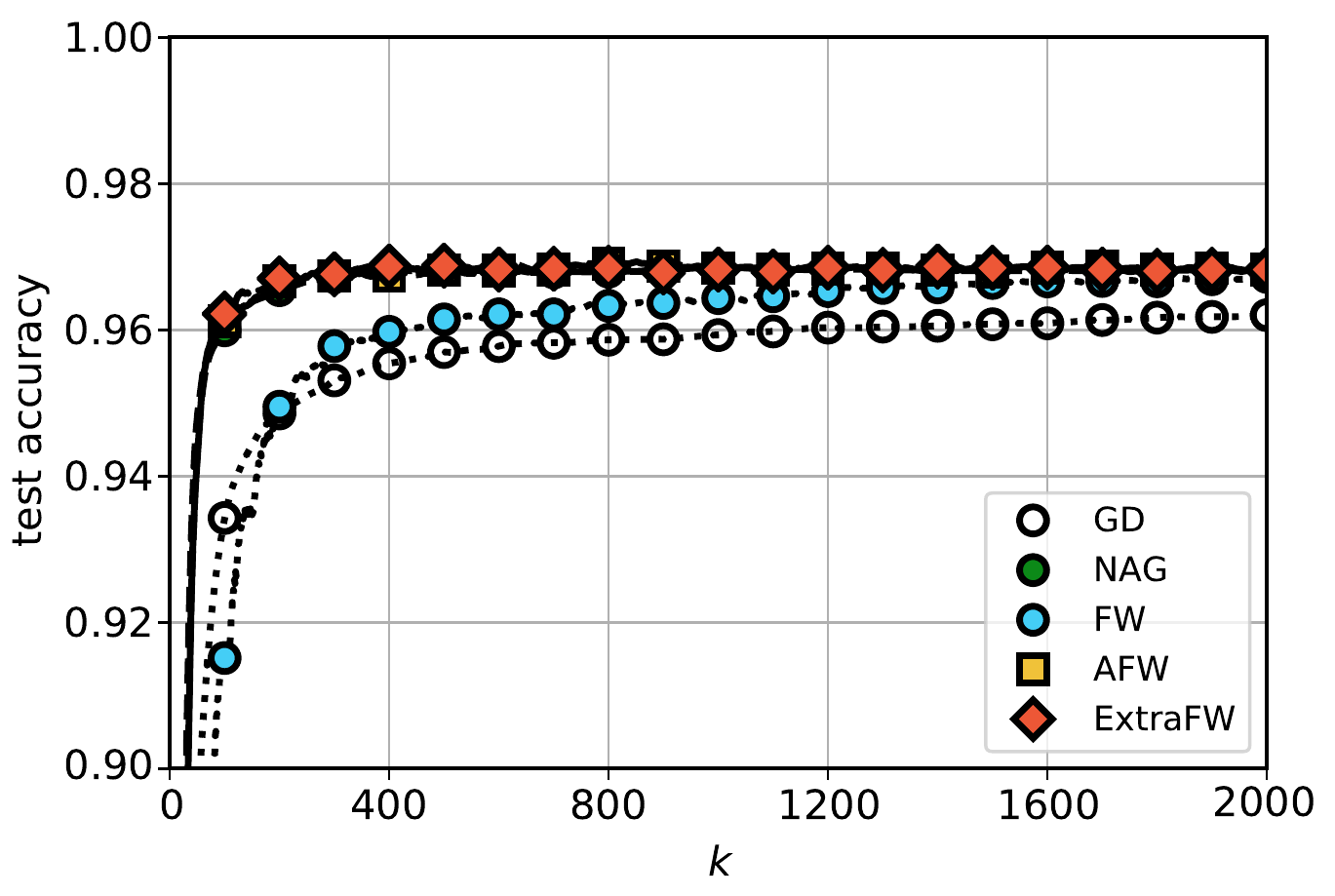}
		\\ (a1) \textit{mnist}, $\ell_2$ norm ball  & (a2) \textit{mushroom}, $\ell_2$ norm ball 	&  (b1) \textit{mnist}, $\ell_1$ norm ball \\
		\hspace{-0.3cm}
		\includegraphics[width=.3\textwidth]{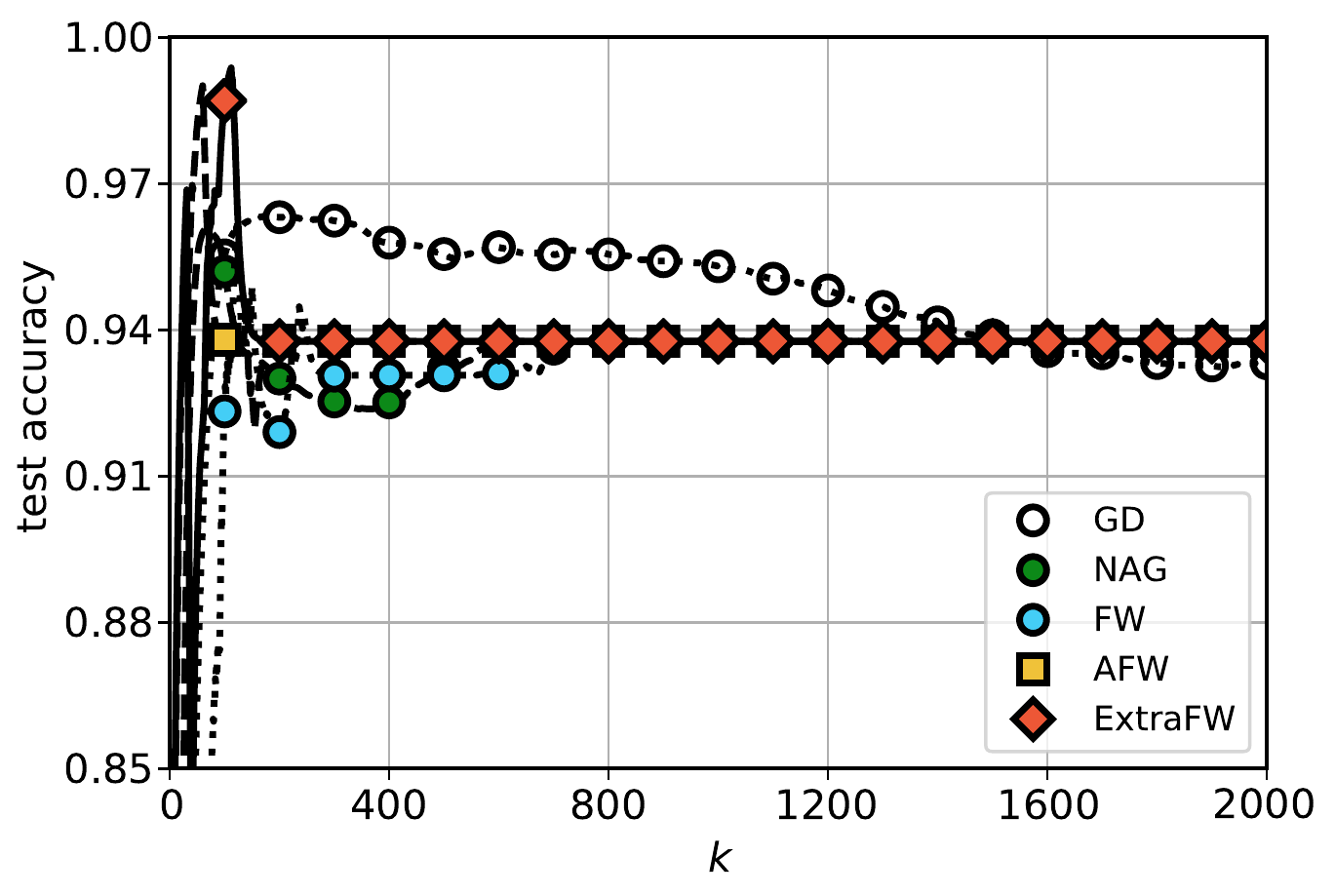}&
		\hspace{-0.3cm}
		\includegraphics[width=.3\textwidth]{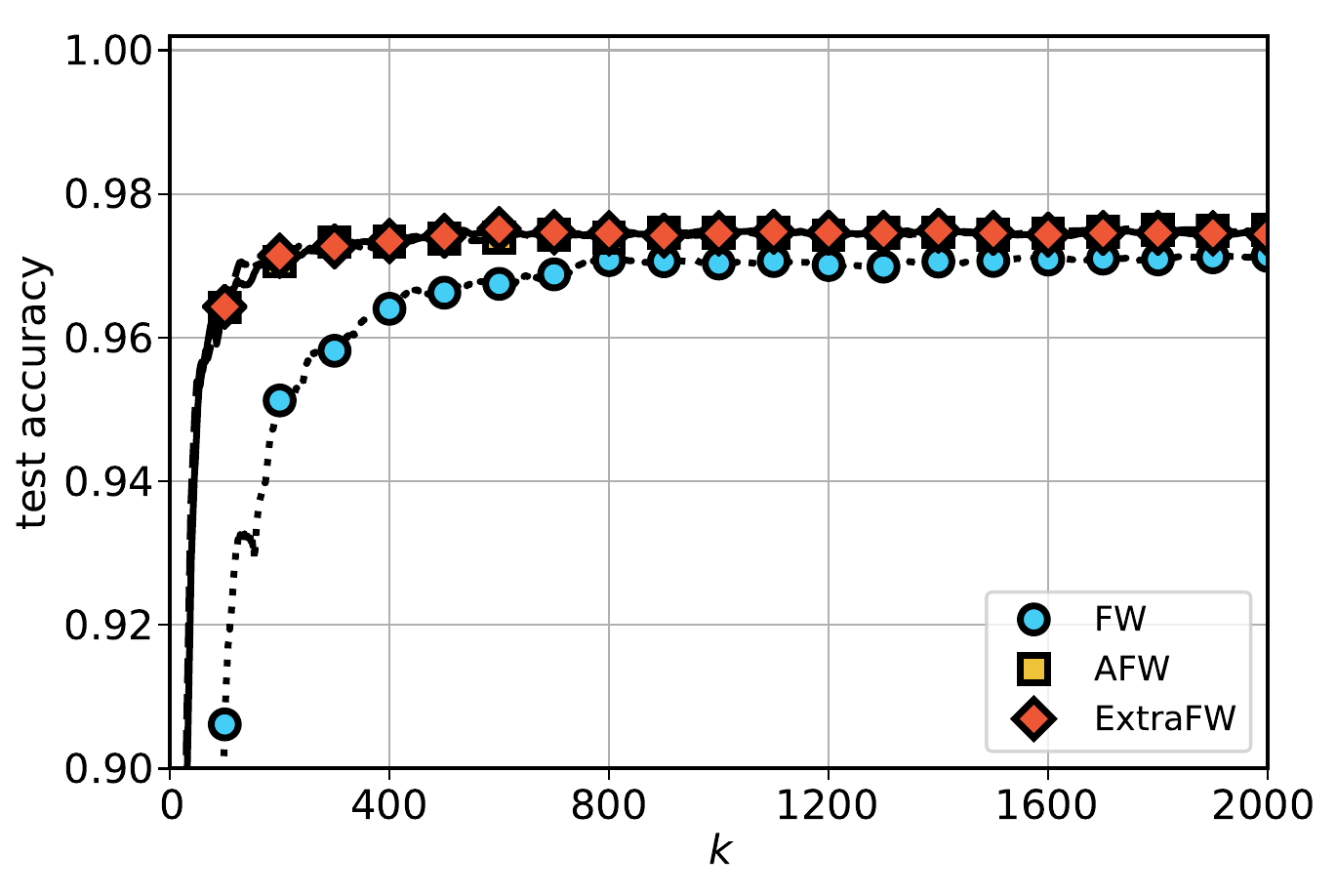}&
		\hspace{-0.3cm}
		\includegraphics[width=.3\textwidth]{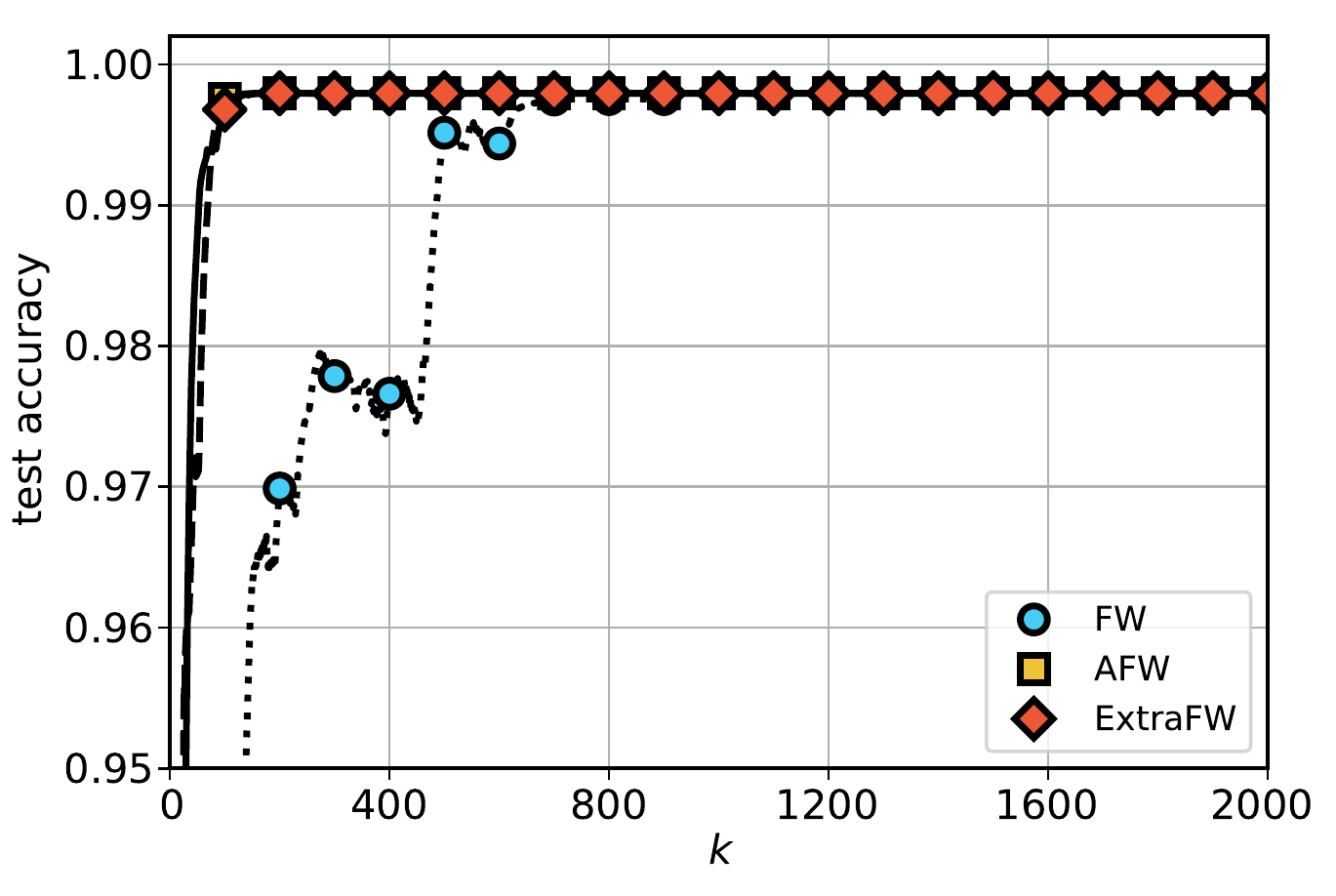}
		\\ (a2) \textit{mushroom}, $\ell_1$ norm ball  & (c1) \textit{mnist}, $n$-supp norm ball 	&  (c2) \textit{mushroom}, $n$-supp norm ball
	\end{tabular}
	\caption{Test accuracy of ExtraFW on different constraints.}
	 \label{fig.apdx_acc}
	 \vspace{-0.4cm}
\end{figure*}
\begin{figure}
\centering
	\begin{tabular}{cc}
		\includegraphics[width=.35\textwidth]{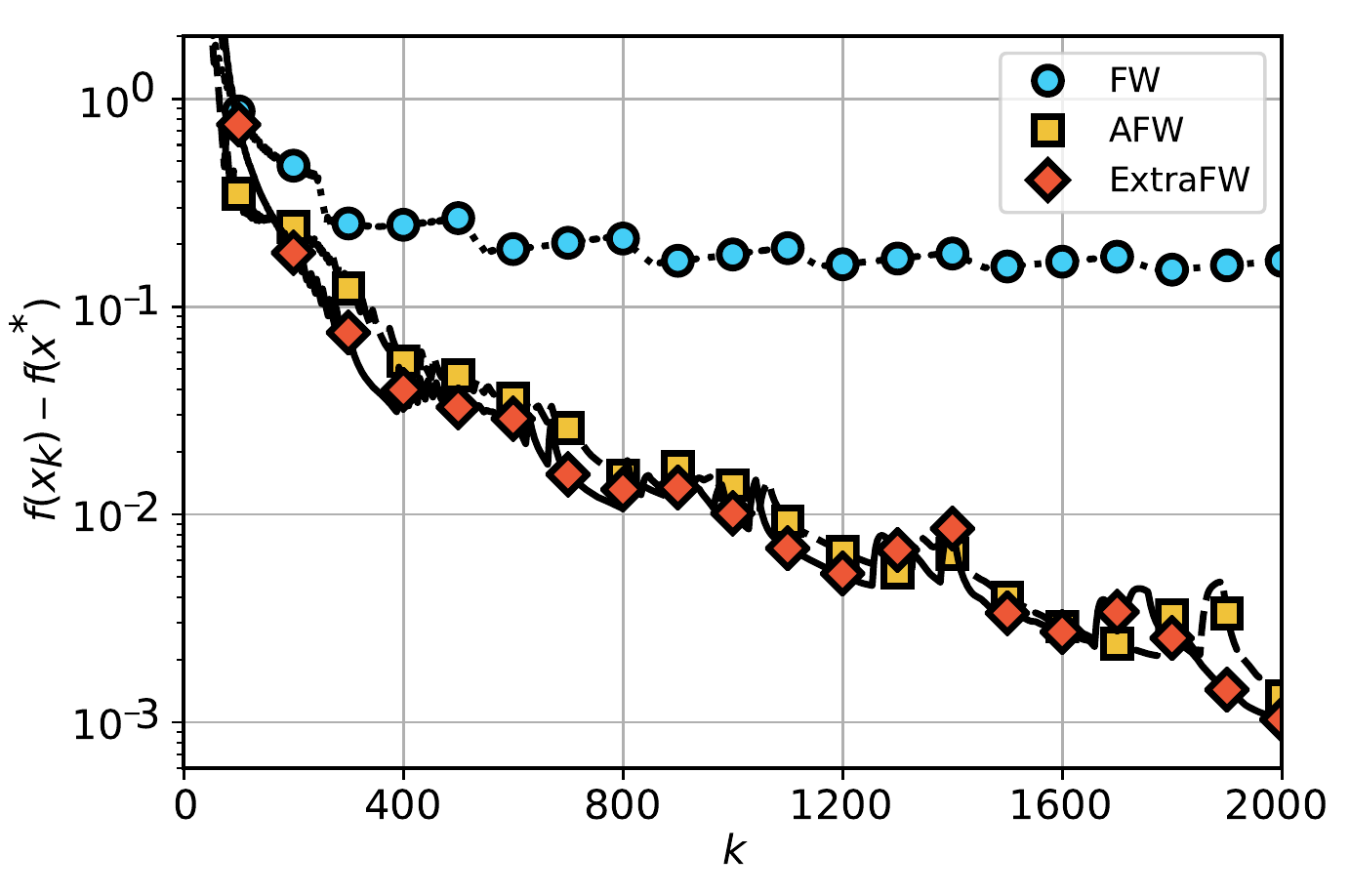}&
		\hspace{-0.3cm}
		\includegraphics[width=.35\textwidth]{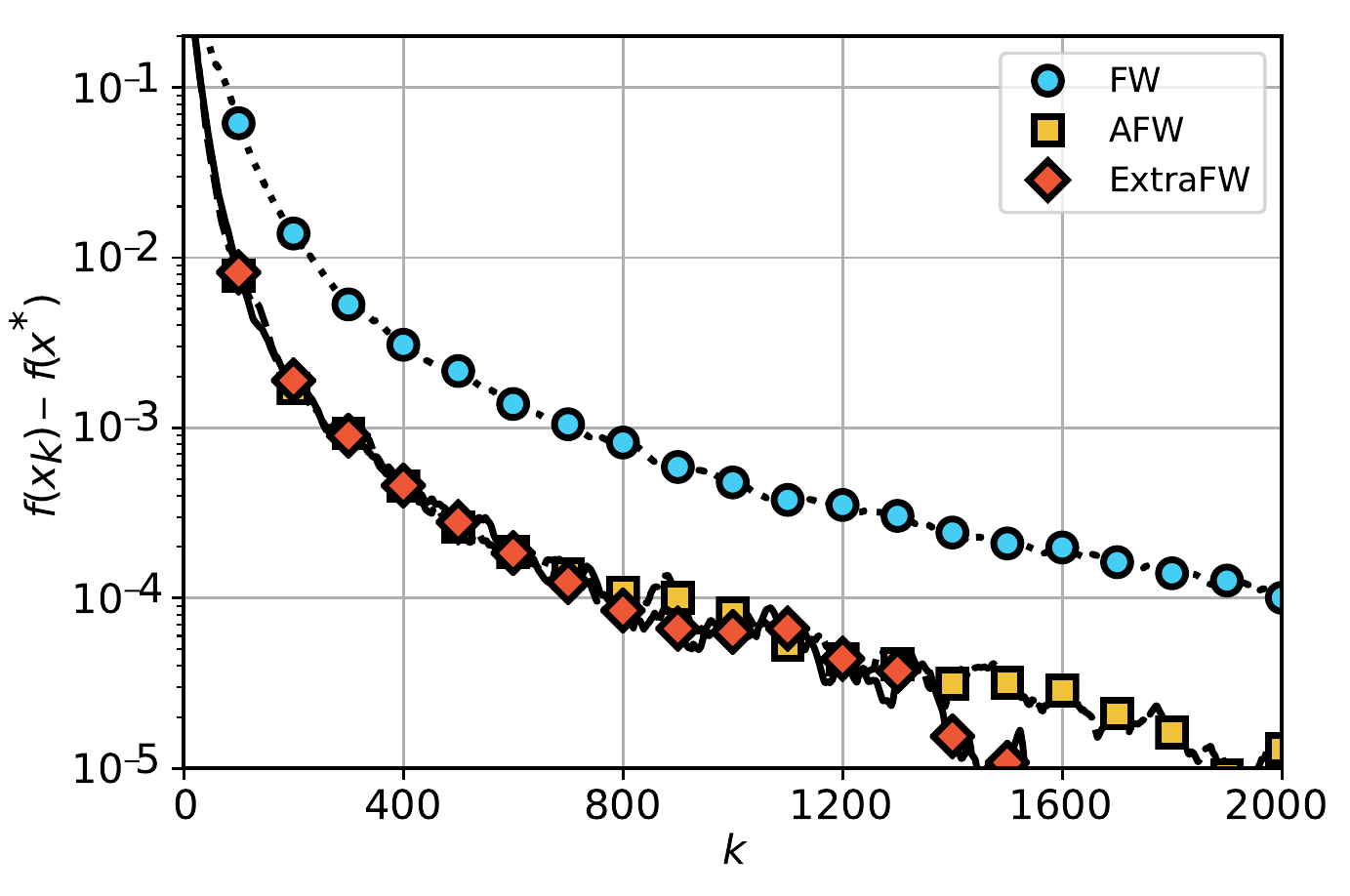}
		\\
		\includegraphics[width=.35\textwidth]{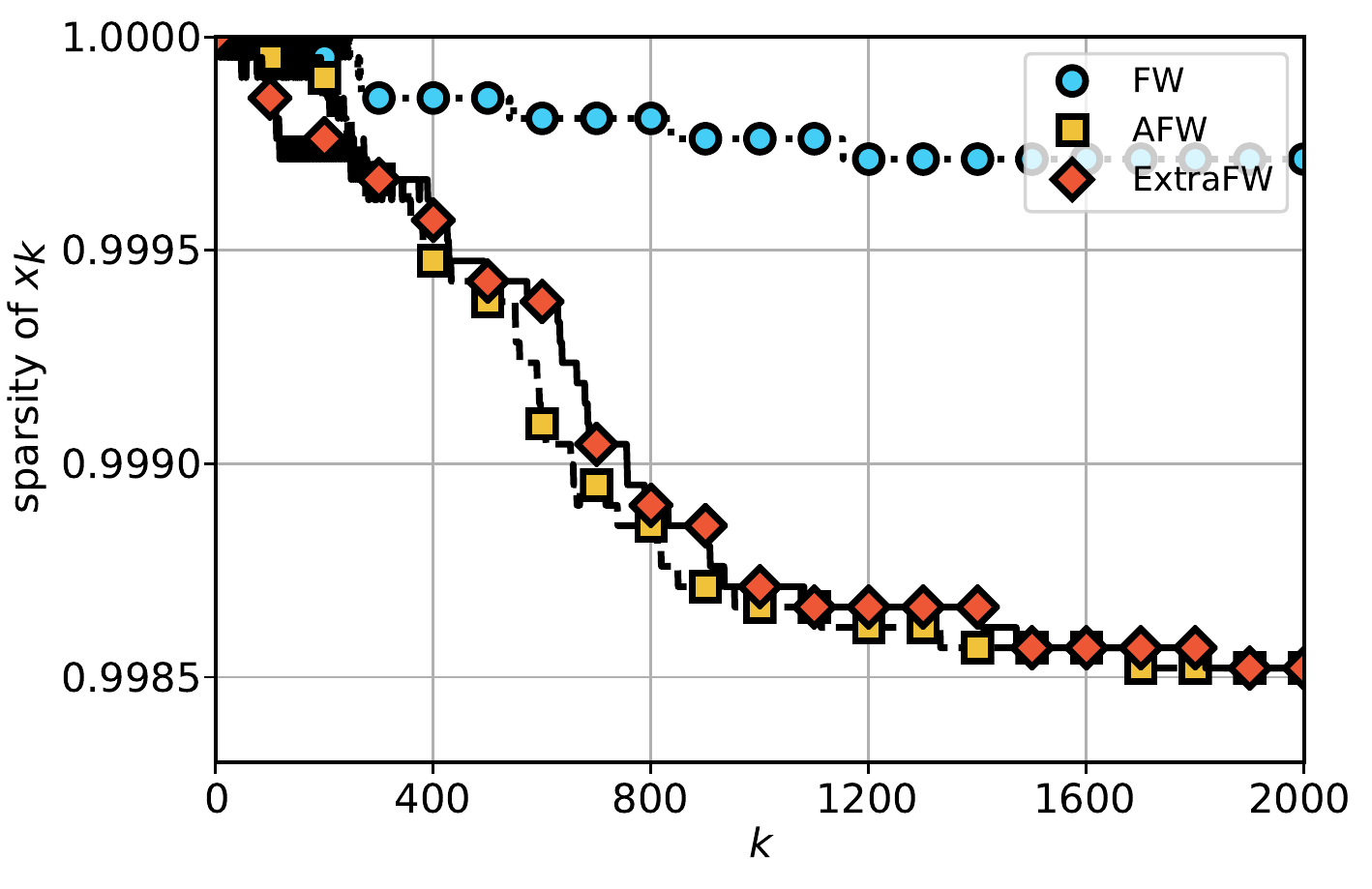}&
		\hspace{-0.3cm}
		\includegraphics[width=.35\textwidth]{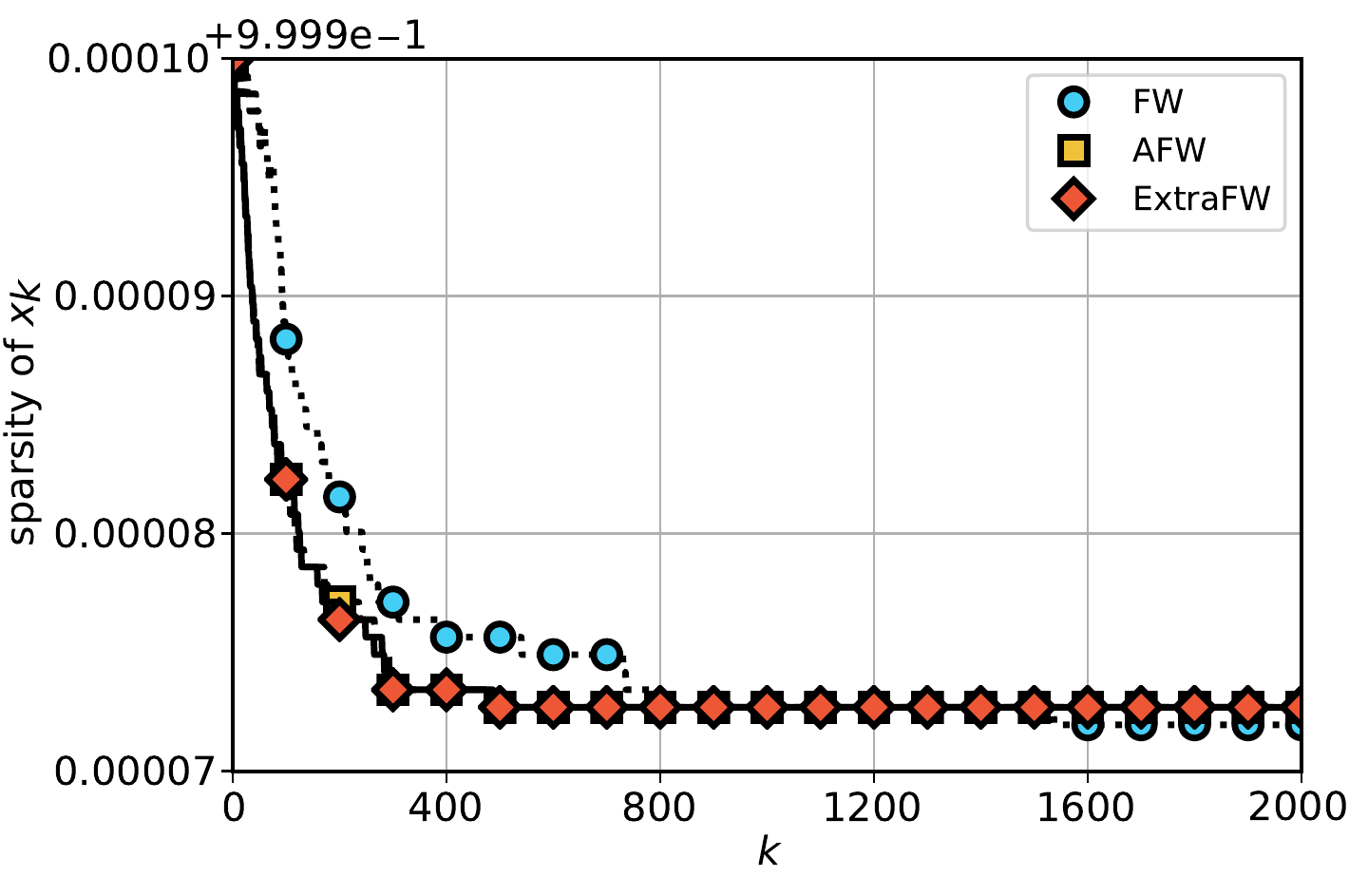}
		\\ (a) \textit{realsim}  & (b) \textit{news20} 
	\end{tabular}
	\caption{Additional tests of ExtraFW for classification with ${\cal X}$ being an $\ell_1$ norm ball.}  \label{fig.l1_supp}
\end{figure}

The test accuracy of different algorithms can be found in Figure \ref{fig.apdx_acc}. Additional numerical results for $\ell_1$ norm ball constraint can be found in Figure \ref{fig.l1_supp}. It can be seen that on dataset \textit{realsim}, ExtraFW has similar performance with AFW, both outperforming FW significantly. On dataset \textit{news20}, ExtraFW outperforms AFW in terms of optimality error.

\begin{figure}
\centering
	\begin{tabular}{cc}
		\includegraphics[width=.35\textwidth]{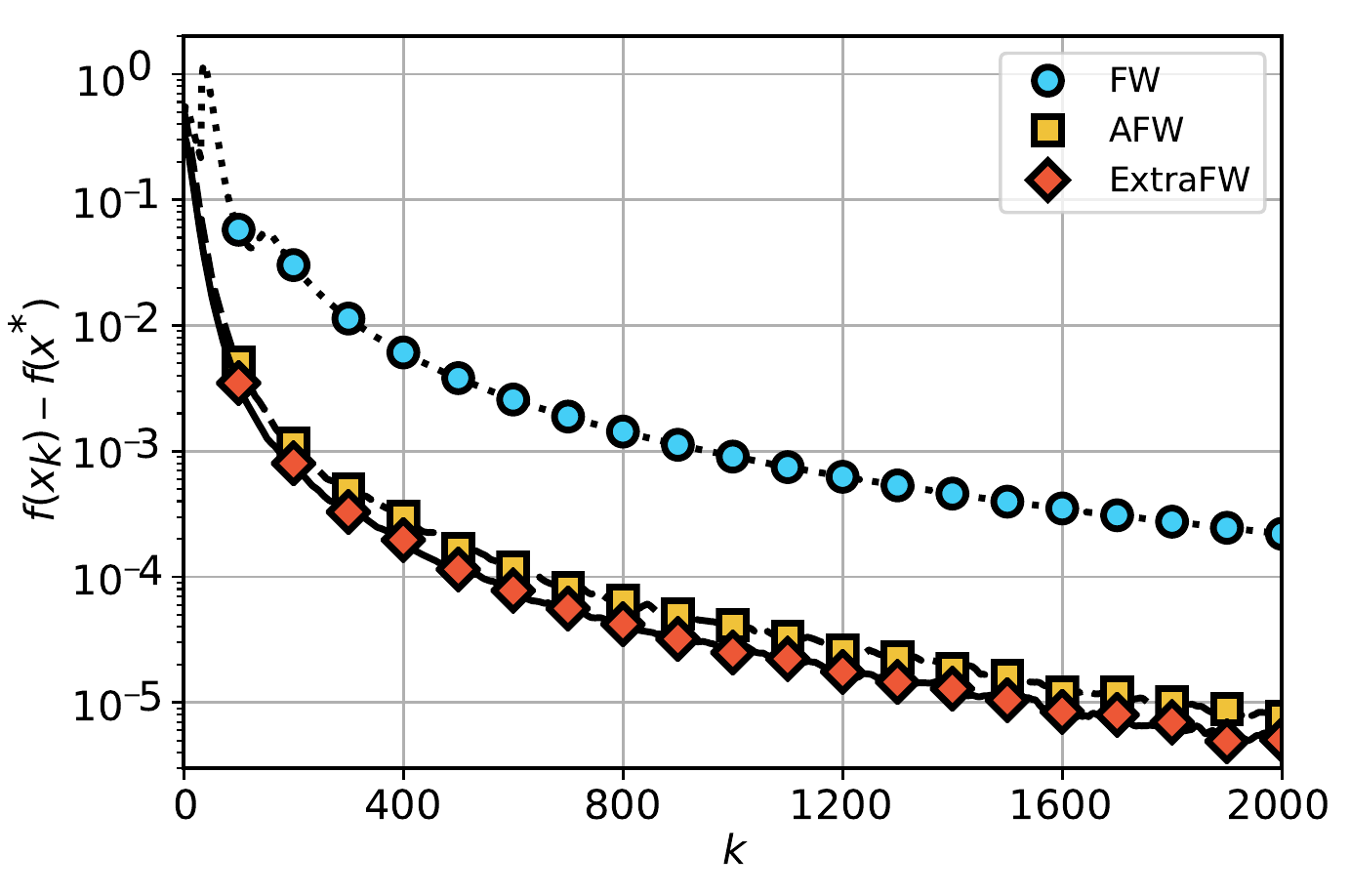}&
		\hspace{-0.3cm}
		\includegraphics[width=.35\textwidth]{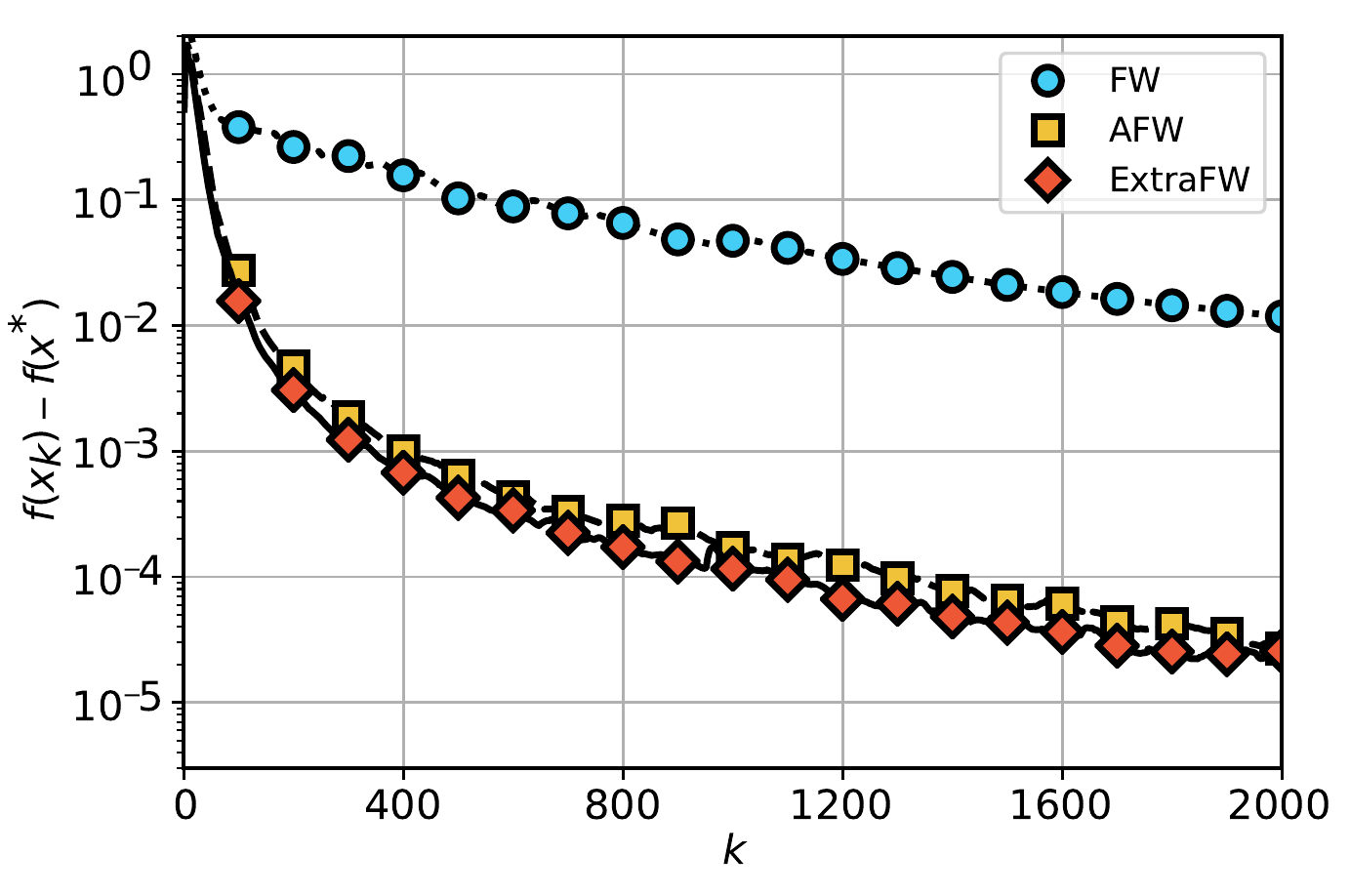}
		\\
		\includegraphics[width=.35\textwidth]{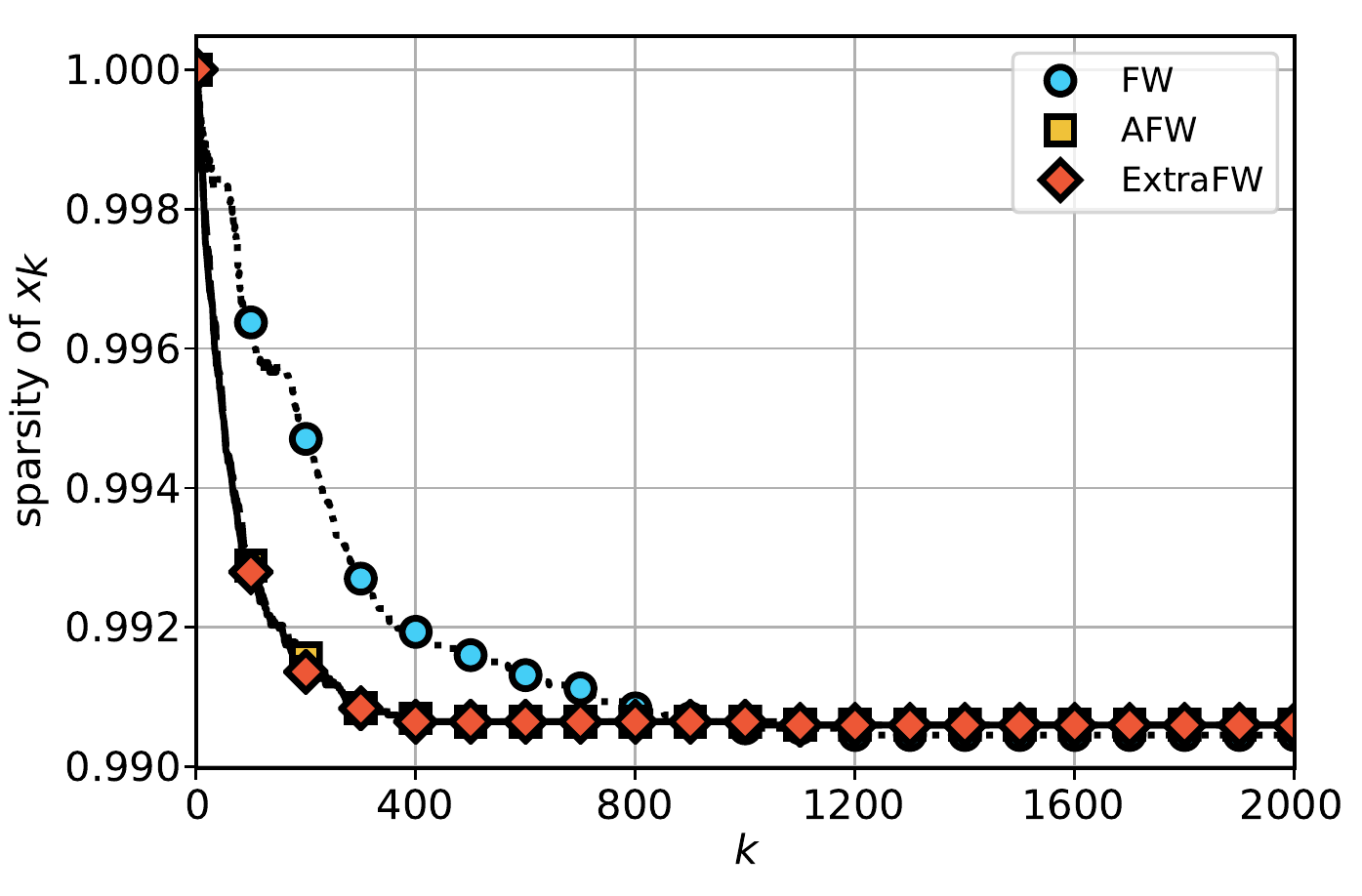}&
		\hspace{-0.3cm}
		\includegraphics[width=.35\textwidth]{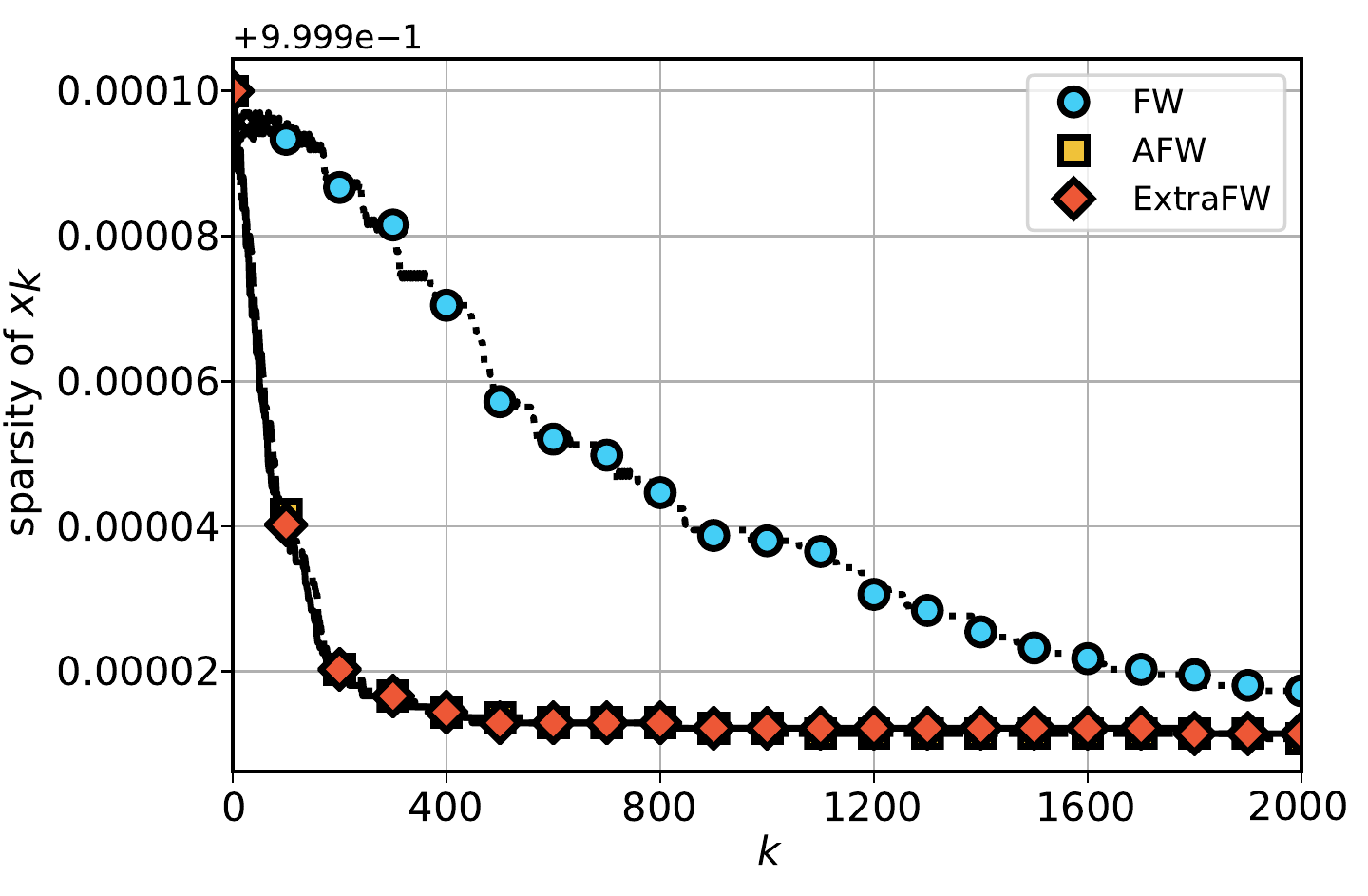}

		\\ (a) \textit{realsim}  & (b) \textit{news20} 
	\end{tabular}
	\caption{Additional tests of ExtraFW for classification with ${\cal X}$ being an $n$-support norm ball.} 
	 \label{fig.n_supp_supp}
\end{figure}
Additional tests for $n$-support norm ball constraint are listed in Figure \ref{fig.n_supp_supp}. The optimality error of ExtraFW is smaller than AFW on both \textit{realsim} and \textit{news20}. 

\subsection{Matrix Completion}\label{apdx.mtrx_completion}
Besides the projection-free property, FW and its variants are more suitable for problem \eqref{eq.mtrx_complete_relax} compared to GD/NAG because they also guarantee ${\rm rank }(\mathbf{X}_k) \leq k+1$ \citep{harchaoui2015,freund2017}. Take FW in Alg. \ref{alg.fw} for example. First it is clear that $\nabla f(\mathbf{X}_k) = (\mathbf{X}_k - \mathbf{A})_{\cal K}$. Suppose the SVD of $\nabla f(\mathbf{X}_k)$ is given by $\nabla f(\mathbf{X}_k) = \mathbf{P}_k\mathbf{\Sigma}_k\mathbf{Q}_k^\top$. Then the FW step can be solved easily by
\begin{align}\label{eq.fw_step_mtrx}
	\bm{V}_{k+1} = - R \mathbf{p}_k \mathbf{q}_k^\top
\end{align}
where $\mathbf{p}_k$ and $\mathbf{q}_k$ denote the left and right singular vectors corresponding to the largest singular value of $\nabla f(\mathbf{X}_k)$, respectively. Clearly $\bm{V}_{k+1}$ in \eqref{eq.fw_step_mtrx} has rank at most $1$. Hence it is easy to see that $\mathbf{X}_{k+1} = (1-\delta_k)\mathbf{X}_{k}+ \delta_k \bm{V}_{k+1}$ has rank at most $k+2$ if $\mathbf{X}_k$ is a rank-$(k+1)$ matrix (i.e., $\mathbf{X}_0$ has rank $1$).  
Using similar arguments, ExtraFW also ensures ${\rm rank }(\mathbf{X}_k) \leq k+1$. Therefore, the low rank structure is directly promoted by FW variants, and a faster convergence in this case implies a guaranteed lower rank $\mathbf{X}_k$.

The dataset used for the test is \textit{MovieLens100K}, where $1682$ movies are rated by $943$ users with $6.30\%$ percent ratings observed. And the initialization and data processing are the same as those used in \citep{freund2017}.

\end{document}